\documentclass[11pt]{article}
\usepackage{amsthm}
\usepackage{amsfonts}
\usepackage{amssymb}
\usepackage{amsmath}
\usepackage{oldgerm}
\numberwithin{equation}{section}
\usepackage{fullpage}
\usepackage{graphicx}%
\usepackage{subcaption}
\usepackage{float}
\allowdisplaybreaks[2]

\newcommand{\R}{\mathbb{R}}
\newcommand{\C}{\mathbb{C}}
\newcommand{\Z}{\mathbb{Z}}
\newcommand{\N}{\mathbb{N}}
\newcommand{\CalK}{\mathcal{K}}
\newcommand{\CalL}{\mathcal{L}}
\newcommand{\CalO}{\mathcal{O}}
\newcommand{\abs}[1]{\lvert#1\rvert}
\DeclareMathOperator{\Span}{span}

\renewcommand{\Re}{\operatorname{Re}}

\newtheorem{theorem}{Theorem}[section]

\theoremstyle{remark}
\newtheorem{remark}[theorem]{Remark}

\allowdisplaybreaks[3]

\title{Numerical Bifurcation and Spectral Stability of Wavetrains in Bidirectional Whitham Models}

\author{Kyle~M.~Claassen\thanks{Department of Mathematics, University of Kansas, 1460 Jayhawk Boulevard, 
Lawrence, KS 66045; kclaassen@ku.edu} \quad \&\quad Mathew~A.~Johnson\thanks{Department of Mathematics, University of Kansas, 1460 Jayhawk Boulevard, 
Lawrence, KS 66045; matjohn@ku.edu} }

\date{\today}

\begin{document}

\maketitle

\begin{abstract}
We consider several different bidirectional Whitham equations that have recently appeared in the literature.  Each of these models
combine the full two-way dispersion relation from the incompressible Euler equations with a canonical shallow water nonlinearity, providing
nonlocal model equations that may be expected to exhibit some of the interesting high-frequency phenomena present in the Euler equations
that standard ``long-wave" theories fail to capture.  Of particular interest here is the existence and stability of periodic traveling wave solutions
in such models.  Using numerical bifurcation techniques we construct global bifurcation diagrams for each system and compare the global
structure of branches, together with the possibility of bifurcation branches terminating in a ``highest" singular (peaked/cusped) wave.  We also numerically
approximate the stability spectrum along these bifurcation branches and compare the stability predictions of these models.  Our results confirm a number of analytical
results concerning the stability of asymptotically small waves in these models and provide new insights into the existence and stability of large amplitude waves.
\end{abstract}


\section{Introduction}  \label{S:intro}

While the Euler equations form the de-facto model for water waves, analysis of this very general model remains difficult, and simpler models have been proposed.  Using KdV as a starting point, which is known to be a good model for waves of large wavelength in shallow water, Whitham \cite{Whitham74} proposed the model
\[
u_t+\mathcal{M}u_{x}+uu_x=0,
\]
where the operator $\mathcal{M}$ is a Fourier multiplier with symbol $m(\xi)=\sqrt{\tanh(\xi)/\xi}$,
that incorporates the full unidirectional dispersion relation for the Euler equations in an effort to better capture mid- and high-frequency phenomena found in the Euler equations, such as wave steepening/peaking.  His equation, now referred to as the \emph{Whitham equation}, has received considerable attention in recent years, having been investigated analytically (\cite{EK09}, \cite{EK13}, \cite{EJ16}, \cite{HJ15}) as well as numerically (\cite{SKCK14}, \cite{EK09}, \cite{RKH17}).  These studies have found that the unidirectional Whitham equation not only captures wave steepening and derivative blow-up \cite{Hur1}, but also exhibits a highest cusped traveling wavetrain \cite{EW}, existence of smooth solitary traveling waves \cite{EW}, 
as well as the famous Benjamin-Feir, or modulational, instability for small amplitude wavetrains \cite{HJ15}.  Furthermore, generalizations of Whitham's equation accounting for constant
vorticity and surface tension have been studied \cite{HJ15b}.  However, Whitham's equation fails to capture
particular high-frequency (spectral) instabilities of small periodic wavetrains due to the unidirectionality of wave propagation; see \cite{DT1,DO11}.

Recently there has been interest in bidirectional Whitham models \cite{MDK15,EJ16,HP_shallow,Pandey2017,EPW17,Carter2017}
in an effort to discover models that capture more of the important qualitative properties of the full Euler equations, 
e.g. existence of smooth solitary waves and peaked traveling waves, and both high-frequency and modulational instabilities of small periodic wavetrains.  
These are often built as full-dispersion generalizations of the nonlinear shallow water equations
\begin{equation}\label{E:shallow}
\left\{\begin{aligned}
	u_t&=-\eta_x-uu_x\\
	\eta_t&=-u_x-(\eta u)_x.
	\end{aligned}\right.
\end{equation}
However, there does not seem to be a unique way to generalize \eqref{E:shallow} to incorporate the full dispersion relation for the Euler equations. Indeed, several different models can be found in the literature, all of which claim to be bidirectional Whitham equations.  
One of the main goals of this paper is to numerically investigate the existence and stability of periodic wavetrains of both small and large amplitude, including the existence of highest peaked / cusped traveling waves,
providing numerical confirmation of known analytical results as well as offering new conjectures for large amplitude waves in hopes of spurring further analytical study.

\

In this paper, we will study three fully-dispersive shallow water wave models that have been presented in the literature.  In \cite{MDK15}, the following bidirectional Whitham equation (written here in non-dimensional form) is derived from the Euler equations:
\begin{equation} \label{E:EJ}
\left\{
\begin{array}{rcl}
u_t &=& -\eta_x - uu_x \\
\eta_t &=& -\CalK u_x - (\eta u)_x
\end{array}
\right.
\end{equation}
where the operator $\CalK$ incorporates the full Euler-dispersion and is defined via its symbol
\begin{equation}\label{K}
\widehat{\CalK f}(n) := \frac{\tanh(n)}{n} \, \widehat{f}(n), \quad n \in \Z.
\end{equation}
Here, $\eta$ represents the fluid height and $u=\phi_x$ with $\phi(x,t)=\phi(x,\eta(x,t),t)$ denoting the trace of the velocity potential at the free surface.  An easy calculation
shows the dispersion relation of \eqref{E:EJ} agrees \emph{exactly} with that of the full Euler equations.  The model \eqref{E:EJ} also appeared in \cite[Section 6.3]{SWX} where it was described
as a linearly well-posed regularization of the linearly ill-posed, yet completely integrable, Kaup system.
While the well-posedness of this system is not fully understood, a recent result \cite[Theorem 1]{EPW17} establishes 
that \eqref{E:EJ} is locally well-posed for all initial data $\eta(x,0)=\psi(x)$, $u(x,0)=\phi(x)$ with $\inf \psi > 0$.  
Though this result does not prove that the system is ill-posed for initial data having negative infimum, numerical evidence 
suggests that this is indeed the case: see Section \ref{SSS:ill-posedness} for further details.  Furthermore, in \cite{EJ16} the existence of $2\pi$-periodic traveling wave
solutions of \eqref{E:EJ} were constructed, and it was proven that there exists a highest wave with a logarithimically-cusped singularity.
Concerning the stability of wavetrains, it has recently been shown that there exists a critical wavenumber $\kappa_c \approx 1.008$ 
such that all $2\pi/\kappa$-periodic wavetrains of \eqref{E:EJ} of sufficiently small amplitude are modulationally unstable provided $\kappa>\kappa_c$: see \cite{Pandey2017}.
Our studies numerically confirm this rigorous modulational instability result, and we investigate the presence of such instabilities for large amplitude wavetrains as well.

Alternatively, \cite{HP_shallow} proposes the non-dimensional, full dispersion shallow water model
\begin{equation} \label{E:HPShallow_model}
\left\{
\begin{array}{rcl}
u_t &=& -\CalK \eta_x - uu_x \\
\eta_t &=& -u_x - (\eta u)_x,
\end{array}
\right.
\end{equation}
where $\mathcal{K}$, $u$, and $\eta$ are as in \eqref{E:EJ}.  Moreover, \cite{HP_shallow} establishes the local well-posedness of \eqref{E:HPShallow_model},
the existence of even, small-amplitude, $2\pi/\kappa$-periodic wavetrains for all $\kappa > 0$, and
that this system possesses a critical wavenumber $\kappa_c \approx 1.610$ such that all periodic waves of sufficiently small amplitude are modulationally unstable
for $\kappa>\kappa_c$.
In Section \ref{SS:analysis_HP} below, we numerically study the existence and stability of wavetrain solutions of \eqref{E:HPShallow_model}, including those of large amplitudes.
We note that, in contrast with \eqref{E:EJ}, it appears branches of smooth periodic traveling waves in \eqref{E:HPShallow_model} 
bifurcating from a zero-amplitude state do not appear to peak/cusp; rather, we posit that such waves of arbitrarily large amplitude are smooth.

Lastly, in \cite[Section 5]{DT1}, a bidirectional \emph{Boussinesq-Whitham} model is proposed:
\begin{equation} \label{E:BoussinesqWhitham_equation}
u_{tt} = \partial_x^2 (u^2 + \CalK u),
\end{equation}
which we write as a first order system:
\begin{equation} \label{E:BoussinesqWhitham_model}
\left\{
\begin{array}{rcl}
u_t &=& \eta \\
\eta_t &=& \partial_x^2 (u^2 + \CalK u).
\end{array}
\right.
\end{equation}
This model is based on the so-called ``bad Boussinesq" equation
\[
u_{tt}=\partial_x^2\left(u^2+u+u_{xx}\right)
\]
from shallow water theory.  Although \eqref{E:BoussinesqWhitham_model} is known to be ill-posed, see \cite[Appendix C]{HP_shallow}, it has nevertheless been shown to exhibit high-frequency instabilities of small amplitude periodic wavetrains, which are known to exist in the full Euler equations; see \cite{Mc82,DO11}.  In Section \ref{SS:BW} below, we numerically 
study the existence and stability of large-amplitude periodic wavetrains of \eqref{E:BoussinesqWhitham_model}, and provide numerical evidence that the bifurcation branch
of such such solutions terminates with a peaked or cusped wave.

\

To streamline our discussion and to set up a framework amenable to ``black-box" computations, we will view each of the systems \eqref{E:EJ}, \eqref{E:HPShallow_model}, and \eqref{E:BoussinesqWhitham_model} in a very general framework of the form
\begin{equation} \label{E:general_model}
\left\{
\begin{array}{rcl}
u_t &=& F(\vec{u},\vec{\eta}) \\
\eta_t &=& G(\vec{u},\vec{\eta}),
\end{array}
\right.
\end{equation}
where
\[
\vec{u} := (u, \partial_x u, \partial_x^2 u, \ldots, \partial_x^n u), \quad \vec{\eta} := (\eta, \partial_x \eta, \partial_x^2 \eta, \ldots, \partial_x^n \eta).
\]
For the models being considered, one of $F$ or $G$ will incorporate the Euler dispersion operator $\CalK$.
Upon transforming $x \mapsto x-ct$ in \eqref{E:general_model}, traveling wave solutions with wavespeed $c$ satisfy (with a slight abuse of notation)
\begin{equation} \label{E:general_model_travwave}
\left\{
\begin{array}{rcl}
u_t &=& cu_x + F(\vec{u},\vec{\eta}) =: F(\vec{u},\vec{\eta}; c) \\
\eta_t &=& c\eta_x + G(\vec{u},\vec{\eta}) =: G(\vec{u},\vec{\eta}; c), \\
\end{array}
\right.
\end{equation}
and equilibrium solutions $(u(x,t),\eta(x,t)) = (\phi(x),\psi(x))$ of the above traveling wave system satisfy the profile equations
\begin{equation} \label{E:general_profile_equations}
\left\{
\begin{array}{rcl}
F(\vec{\phi},\vec{\psi}; c) &=& 0 \\
G(\vec{\phi},\vec{\psi}; c) &=& 0.
\end{array}
\right.
\end{equation}
Assuming that one can resolve either $\psi(x)$ in terms of $\phi(x)$ or $\phi(x)$ in terms of $\psi(x)$ in \eqref{E:general_profile_equations}, which will fortunately be the case for our intended applications,
we can conveniently express a single profile equation for, say, $\phi$ in the form
\begin{equation} \label{E:profile_equation}
\CalK \phi = g(c,\phi),
\end{equation}
where $\mathcal{K}$ is as in \eqref{K} and $g$ implicitly relates the wavespeed $c$ and its corresponding profile $\phi$.  Using this framework, we implement numerical methods for each model to compute a global bifurcation of periodic solutions from a constant amplitude state, as well as the spectrum of the linearization for a sampling of waves of various heights.  See Section \ref{S:bifurcation_methods} and Appendix \ref{A:general_FFHM} for detailed discussions of these methods.

\

Through numerical experiments, we find that, from the perspective of bifurcation and stability of periodic wavetrains, the three models discussed above have many similar qualitative features.  
We find that periodic traveling waves in each of these models indeed exhibit both high-frequency and modulational instabilities for waves of sufficiently large amplitude.  
However, high-frequency instabilities for small amplitude waves as predicted by the full Euler equations are numerically imperceptible for the models considered here, suggesting that these models fail to capture high-frequency instabilities in arbitrarily small amplitude waves.  Furthermore, we find that while the global bifurcation branches of smooth periodic wavetrains for the models \eqref{E:EJ}
and \eqref{E:BoussinesqWhitham_equation} seem to terminate in a highest singular (peaked/cusped) wave, our numerical experiments lead us to conjecture that 
the corresponding branches for the model \eqref{E:HPShallow_model} extend to arbitrarily high amplitude smooth solutions.
Furthermore, while the model \eqref{E:HPShallow_model} is known to be locally well-posed \cite{HP_shallow}, the model \eqref{E:EJ} is only known to be \emph{conditionally}
locally well-posed \cite{EPW17}, requiring that the height profile $\eta$ of the initial data be uniformly positive.  In Section \ref{SSS:ill-posedness} below we provide strong numerical evidence
that the local evolution of \eqref{E:EJ} is in fact \emph{ill-posed} if the positivity condition on $\eta$ is removed.  This leads us to construct new branches of 
periodic wavetrains of \eqref{E:EJ} with uniformly positive $\eta$, which had previously been unstudied: see Section \ref{SSS:positive} below.

Finally, we emphasize that the stability/instability found here is only at the spectral level and hence may not correspond to actual stability or instability
in the full nonlinear dynamics.  For an interesting investigation on how nonlinear terms may stabilize a spectrally unstable wave, see \cite{DDC17}.
\

The outline of this paper is as follows.  In Section \ref{S:bifurcation_methods} we discuss the numerical bifurcation techniques used in our studies to construct
the global bifurcation curves of periodic wavetrain solutions with fixed period.  In Section \ref{S:results} we present our main results for each of the models
\eqref{E:EJ}, \eqref{E:HPShallow_model}, and \eqref{E:BoussinesqWhitham_equation}: see Sections \ref{SS:analysis_EJ}, \ref{SS:analysis_HP}, and \ref{SS:BW}, respectively.
Further, in Section \ref{S:bf0} describes in detail the numerical methods used to approximate the stability spectrum associated to periodic wavetrains of \eqref{E:EJ}.
Section \ref{S:conclusion} provides a concluding summary of the main observations of this paper.
In Appendix \ref{A:general_FFHM} we detail how these numerical stability methods extend to more general models, including \eqref{E:HPShallow_model} and \eqref{E:BoussinesqWhitham_equation}.  Throughout, the parameters of the numerical methods (e.g. number of Fourier modes used) are chosen by experiment to ensure that the behavior of the computed waves and their corresponding spectra are qualitatively correct.  See Appendix \ref{A:numerical_parameters} for further details.

\

\noindent
{\bf Acknowledgments:}  The authors would like to thank Mats Ehrnstr\"om and Vera Hur for several useful conversations regarding the models studied here, as well as John Carter for his help in confirming the numerical ill-posedness results in Section \ref{SSS:ill-posedness} and Morgan Rozman for assistance in implementing the sixth-order operator splitting method \cite{Yoshida1990} used for time evolution.  MJ was supported by the National Science Foundation under grant DMS-1614785.


\section{Numerical Bifurcation Methods}  \label{S:bifurcation_methods}

In this section we discuss the methods used to numerically approximate the even, $2\pi$-periodic solutions of the bidirectional Whitham models.  In summary, the approximations will be obtained by truncating an expansion of the profile $\phi$ in an appropriate Fourier basis and discretizing the integrals that yield the Fourier coefficients. The resulting expression will then be used to form a discretized version of the profile equation \eqref{E:profile_equation}, which will be enforced at a fixed number of \textit{collocation points} in order to yield a nonlinear system of equations which may be solved, for example, via Newton's method.  Using the local bifurcation theory for the model to generate an ``initial guess" for the nonlinear system solver, an algorithm known as the \textit{pseudo-arclength method} will be employed as a robust method for simultaneously continuing the values of the wavespeed and its corresponding approximate profiles along a global bifurcation branch.

\subsection{Cosine Collocation Method}  \label{SS:cosine_collocation_method}

Note that even, $2\pi/\kappa$-periodic functions $\phi=\phi(x)$ are naturally represented in a cosine series
\begin{equation} \label{E:cos_series}
\phi(x) = \sum_{n=0}^\infty \widehat{\phi}(n) \cos(n\kappa x),
\end{equation}
where
\begin{equation} \label{E:cos_series_coeff}
\widehat{\phi}(n) = \begin{cases}
\displaystyle \frac{\kappa}{\pi} \int_0^{\pi/\kappa} \phi(x) \, dx & \text{if $n=0$} \vspace{0.5em} \\
\displaystyle \frac{2\kappa}{\pi} \int_0^{\pi/\kappa} \phi(x) \cos(n\kappa x) \, dx & \text{if $n \geq 1$.}
\end{cases}
\end{equation}
Partitioning the interval $[0,\pi/\kappa]$ into $N$ subintervals, for each $n \in \N_0$ we discretize the integrals in \eqref{E:cos_series_coeff} by the midpoint method as
\begin{equation} \label{E:fourier_integral_discretization}
\int_0^{\pi/\kappa} \phi(x) \cos(n\kappa x) \, dx \approx \sum_{i=1}^N \phi(x_i) \cos(n\kappa x_i) \frac{\pi}{\kappa N},
\end{equation}
where the $x_i = \dfrac{(2i-1)\pi}{2\kappa N}$ for $i=1,2,\ldots,N$ are the so-called \emph{collocation} points on $[0,\pi/\kappa]$.  Using this approximation in \eqref{E:cos_series_coeff}, we obtain the following approximations of the series coefficients:
\begin{equation}  \label{E:approx_cos_series_coeff}
\widehat{\phi}(n) \approx \widehat{\phi}_N(n) := w(n) \sum_{i=1}^N \phi(x_i) \cos(n\kappa x_i),
\end{equation}
where 
\[
w(n) := \begin{cases}
1/N & \text{if $n=0$} \\
2/N & \text{if $n \geq 1$.}
\end{cases}
\]
Using this approximation for the series coefficients in \eqref{E:cos_series}, and truncating the series to $N$ terms, we obtain
\begin{align}
\phi(x) \approx \phi_N(x) &:= \sum_{n=0}^{N-1} \left[ w(n)\sum_{i=1}^N \phi(x_i) \cos(n\kappa x_i) \right] \cos(n\kappa x) \nonumber \\
&= \sum_{i=1}^N \left[ w(n) \sum_{n=0}^{N-1} w(n)\cos(n\kappa x_i)\cos(n\kappa x) \right] \phi(x_i).  \label{E:phi_N}
\end{align}
Moreover, note that since $\phi$ is even, so is $\CalK\phi$, where again $\mathcal{K}$ is defined as in \eqref{K} .  It follows 
we can approximate $\CalK\phi(x)$ similarly by truncating cosine series expansion and again approximating coefficients with \eqref{E:approx_cos_series_coeff}:
\begin{align}
\CalK\phi(x) &= \sum_{n=0}^\infty \widehat{\CalK\phi}(n) \cos(n\kappa x) \nonumber \\
&\approx \sum_{n=0}^{N-1} \widehat{\CalK\phi_N}(n) \cos(n\kappa x) \nonumber \\
&= \sum_{n=0}^{N-1} \widehat{\CalK}(n) \widehat{\phi}_N(n) \cos(n\kappa x) \nonumber \\
&= \sum_{i=1}^N \left[ \sum_{n=0}^{N-1} w(n) \frac{\tanh(\kappa n)}{\kappa n} \cos(n\kappa x_i) \cos(n\kappa x) \right] \phi(x_i) =: (\CalK\phi)_N(x).
\end{align}
Thus we arrive at an approximate profile equation:
\begin{equation} \label{E:approximate_profile_equation}
(\CalK\phi)_N(x) = g(c,\phi_N(x)).
\end{equation}
For notational convenience, going forward we denote evaluation of $\phi_N$ and $(\CalK\phi)_N$ at each collocation point $x_i$ via a superscript:
\[
\phi_N^i := \phi_N(x_i), \quad \text{and} \quad (\CalK\phi)_N^i := (\CalK\phi)_N(x_i).
\]
Enforcing \eqref{E:approximate_profile_equation} at each $x_i$ yields a nonlinear system
\begin{equation}  \label{E:discretized_profile_system}
f_i(c,\phi_N^1,\phi_N^2,\ldots,\phi_N^N) := g(c,\phi_N^i) - (\CalK\phi)_N^i = 0, \quad i=1,\ldots,N,
\end{equation}
which we will endeavor to solve for wavespeed $c$ and approximate points $\phi_N^1, \ldots, \phi_N^N$ on the corresponding profile.
\begin{remark}
Note that the system defined by \eqref{E:discretized_profile_system} is actually underdetermined, as it contains $N+1$ unknowns but only $N$ equations.  However, as discussed below, the numerical continuation algorithm will impose an additional condition as these solutions are computed along the global bifurcation branch, at each step yielding a well-defined, full-rank system.
\end{remark}

\subsection{Numerical Continuation by the Pseudo-Arclength Method} \label{SS:pseudoarclength_method}

To solve the discretized profile system \eqref{E:discretized_profile_system}, we use the \emph{pseudo-arclength} method, a numerical continuation algorithm which is well-known to the numerical bifurcation community.  Nevertheless, we provide a summary of the method here.

\begin{figure}[th] 
\centering
\includegraphics[scale=0.8,clip=true,trim=1.5cm 1cm 1.5cm 0.5cm]{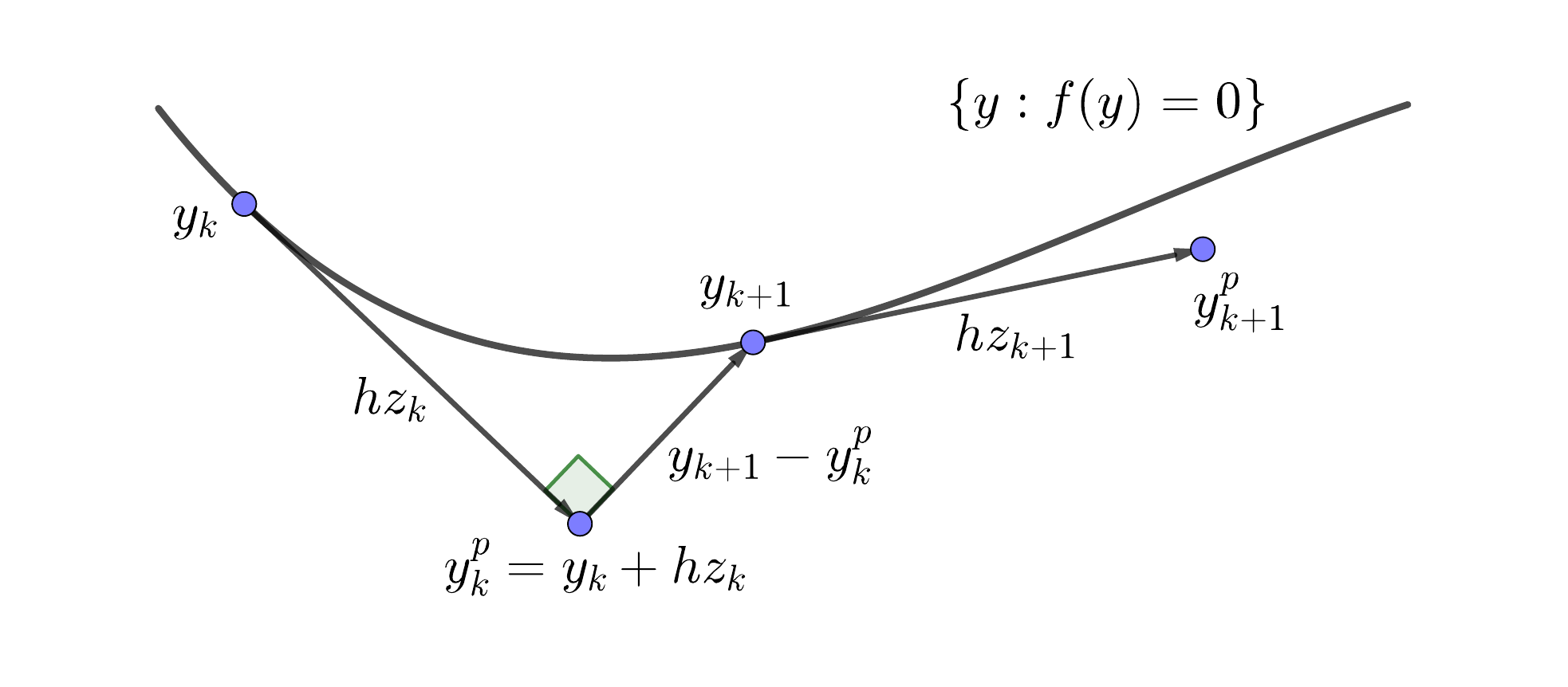}
\caption{ \small
Illustration of the pseudo-arclength method: via a predictor-corrector scheme, for a given point $y_k$ such that $f(y_k)=0$ the method computes $y_{k+1}$ such that $f(y_{k+1})=0$ and a consistently-oriented tangent direction $z_{k+1}$ at $y_{k+1}$.
}\label{F:pseudoarclength_method}
\end{figure}

For further notational convenience, let $y := (c,\phi_N^1,\ldots,\phi_N^N) \in \R^{N+1}$ and define $f : \R^{N+1} \to \R^N$ by  $f(y) := (f_1(y),\ldots,f_N(y))$.  Then to solve \eqref{E:discretized_profile_system}, we seek $y$ such that
\begin{equation} \label{E:compact_nonlinear_equation}
f(y) = 0.
\end{equation}
The pseudo-arclength method is a \emph{predictor-corrector} method; given a point on the solution curve $f(y)=0$, another solution is found by first applying an extrapolation (\emph{predictor}) from the known solution, followed by a \emph{corrector} process that projects the extrapolation back onto the solution curve.  More specifically, the pseudo-arclength method follows the following program:
\begin{enumerate}
\item Given a point $y_0 \in \R^{N+1}$ with $f(y_0) = 0$, a unit tangent direction $z_0$ to the solution curve at $y_0$, and a step size $h$, form the predictor $y_1^p$ by extrapolating $h$ units along $z_0$, i.e. $y_1^p := y_0 + hz_0$.

\item Project $y_0^p$ back onto the solution curve $f(y)=0$ along the direction orthogonal to $z_0$.  That is, solve for $y_1$ in the augmented nonlinear system
\[
\left\{
\begin{array}{rcl}
f(y_1) &=& 0 \\
z_0 \cdot (y_1 - y_1^p) &=& 0.
\end{array}
\right.
\]

\item For the next step of the method, a suitable tangent vector $z_1$ to the solution curve at the point $y_1$ can be found by solving for $z_1$ in the following system of $N+1$ equations and $N+1$ unknowns:
\begin{equation}  \label{E:pseudoarclen_tanvecsys}
\left\{
\begin{array}{rcl}
Df(y_1)z_1 &=& 0 \\
z_0 \cdot z_1 &=& 1.
\end{array}
\right.
\end{equation}
\end{enumerate}
Note that the $z_1$ solved for in \eqref{E:pseudoarclen_tanvecsys} will not necessarily be of unit length and should therefore be normalized before further iteration.  Steps (1)-(3) above can be iterated to continue from a solution $y_k$ to another point $y_{k+1}$ such that $f(y_{k+1})=0$.  See Figure \ref{F:pseudoarclength_method} for a graphical illustration of the method.

\begin{remark}
The first equation in \eqref{E:pseudoarclen_tanvecsys} is the \emph{tangency condition} ensuring that $z_1$ is tangent to the solution curve at $y_1$.  The second equation is an \emph{orientation condition} guaranteeing that the angle between $z_0$ and $z_1$ is acute.  Hence the tangent vectors will be consistently oriented from one iteration to the next so that the method always makes ``forward progress", not backtracking toward the previous solution point.  Together, these conditions form a full-rank system that can also be solved using a nonlinear equation solver.  In theory, any positive number could be used on the right-hand side of the orientation condition  of \eqref{E:pseudoarclen_tanvecsys}.  The only difference would be the magnitude of the solution $z_1$, which is inconsequential as this initial solution is scaled to unit length before further iteration.
\end{remark}

As is evident from Step (1) of the pseudo-arclength algorithm given above, it is necessary to start the method with an initial point $y_0$ such that $f(y_0)=0$.  We can find such a point near the start of the bifurcation branch by using the local bifurcation theory for the model, which yields curves $c=c(\varepsilon)$ and $\phi=\phi(x,\varepsilon)$ that trace the bifurcation branch for $\varepsilon > 0$ small.  Indeed, for small $\varepsilon_0 > 0$ fixed, consider
\begin{equation}  \label{E:y0_star}
y_0^* := \left( c(\varepsilon_0), \phi(x_1;\varepsilon_0), \phi(x_2;\varepsilon_0), \ldots, \phi(x_N;\varepsilon_0) \right).
\end{equation}
This $y_0^*$ does not necessarily solve $f(y_0^*) = 0$, but by virtue of $\varepsilon_0$ being small, $y_0^*$ should be \emph{close} to the solution curve.  Thus we may use $y_0^*$ as a good ``initial guess" for a nonlinear system solver applied to finding $\left( \phi_{N,0}^1, \phi_{N,0}^2, \ldots, \phi_{N,0}^N \right)$ such that
\[
f\left( c(\varepsilon_0), \phi_{N,0}^1, \phi_{N,0}^2, \ldots, \phi_{N,0}^N \right) = 0.
\]
We use this solution to form the initial point on the solution curve:
\begin{equation}  \label{E:initial_y0}
y_0 := \left( c(\varepsilon_0), \phi_{N,0}^1, \phi_{N,0}^2, \ldots, \phi_{N,0}^N \right).
\end{equation}
Moreover, since the local bifurcation curve is parameterized by $\varepsilon$ for $\varepsilon > 0$ small, we compute the tangent direction to the local bifurcation curve at $y_0^*$ by differentiating with respect to $\varepsilon$ at $\varepsilon_0$ and normalizing to unit length:
\[
z_0^* := \left( c'(\varepsilon_0), \frac{\partial\phi}{\partial \varepsilon}(x_1; \varepsilon_0), \frac{\partial\phi}{\partial \varepsilon}(x_2; \varepsilon_0), \ldots, \frac{\partial\phi}{\partial \varepsilon}(x_N; \varepsilon_0) \right), \quad z_0 := \frac{z_0^*}{\abs{z_0^*}},
\]
Now, $z_0^*$ is not necessarily tangent to the solution curve, but since $\varepsilon_0 > 0$ is small, $z_0^*$ is very nearly tangent to the solution curve at $y_0$ as found above, and we take $z_0 := z_0^* / \abs{z_0^*}$ as the initial tangent direction.

\begin{remark}
In practice it typically does not matter that $z_0$ is not exactly tangent to the solution curve at $y_0$.  For that matter, $y_0^*$ as defined in \eqref{E:y0_star} can typically be used as the initial starting point without the preliminary application of the solver discussed above to find a $y_0$ as defined in \eqref{E:initial_y0}.  The accuracy of $y_0^*$ and $z_0^*$ affects the quality of the first predictor $y_1^p$, but for small $\varepsilon_0$ these initial inaccuracies are mitigated by the Newton solver in steps (2) and (3) of the method to obtain $y_1$ and $z_1$.  Then, fortunately, all future iterations will have highly accurate starting points and tangent directions from which to form their predictors.
\end{remark}


\section{Numerical Results}  \label{S:results}

In this section, we present our main numerical results for each of the bidirectional Whitham models \eqref{E:EJ}, \eqref{E:HPShallow_model}, and \eqref{E:BoussinesqWhitham_model}.
Precisely, we implement the program outlined by the numerical bifurcation methods described in Section \ref{S:bifurcation_methods} to compute small and large amplitude traveling
periodic profiles, and we produce the global bifurcation diagrams of periodic traveling waves with fixed period.  
Moreover, we study the stability of these waves to localized, i.e. integrable on $\mathbb{R}$, perturbations by 
using a Fourier-Floquet-Hill method \cite{DK05} to numerically compute the spectrum of the associated linearized operators.  We begin our analysis by demonstrating the relevant
details of the computations for  the model \eqref{E:EJ}.  
Since the details for the models \eqref{E:HPShallow_model} and \eqref{E:BoussinesqWhitham_model} are highly similar, 
only a summary of the relevant formulas will be provided for these other models.

\subsection{Analysis of System \eqref{E:EJ}}\label{SS:analysis_EJ}

We begin our analysis by considering the existence of periodic traveling waves to the bidirectional Whitham model \eqref{E:EJ}:
\[
\left\{
\begin{array}{rcl}
u_t &=& -\eta_x - uu_x \\
\eta_t &=& -\CalK u_x -(\eta u)_x.
\end{array}
\right.
\]
Recall in the modeling of shallow water waves, $\eta$ represents the fluid height and $u$ is the trace of the velocity potential at the free surface.  
Transforming to traveling wave coordinates $x \mapsto x-ct$, such solutions are seen to be stationary, spatially periodic solutions $(u(x,t),\eta(x,t))=(\phi(x),\psi(x))$  
of the evolutionary equation
\begin{equation} \label{E:EJ_traveling_wave_model}
\left\{
\begin{array}{rcl}
u_t &=& -\eta_x - uu_x + cu_x =: F(u,u_x,\eta,\eta_x; c) \\
\eta_t &=&  -\CalK u_x - (\eta u)_x + c\eta_x =: G(u,u_x,\eta,\eta_x; c),
\end{array}
\right.
\end{equation}
and hence the profiles $\phi,\psi$ satisfy the system
\begin{equation} \label{E:EJ_traveling_profile_system}
\left\{
\begin{array}{rcl}
0 &=& -\psi' - \phi\phi' + c\phi' \\
0 &=& -\CalK\phi' - (\psi\phi)' + c\psi'.
\end{array}
\right.
\end{equation}
Upon integrating \eqref{E:EJ_traveling_profile_system} and setting constants of integration to zero\footnote{Note by Galilean invariance, one of the two integration constants can always be set to zero.  
For our analysis, we set the other integration constant to zero for simplicity.  See \cite{EJ16} for details.}, we can resolve $\psi$ in terms of $\phi$ in \eqref{E:EJ_traveling_profile_system} to obtain
\begin{align}
\psi &= c\phi - \frac{1}{2}\phi^2 \nonumber \\
\CalK\phi &= \frac{1}{2}\phi^3 - \frac{3}{2}c\phi^2 + c^2 \phi =: g(c,\phi).\label{E:EJ_profile}
\end{align}
As explained in \cite{EJ16}, one should not expect \eqref{E:EJ_profile} to admit smooth solutions of arbitrary amplitude.  Indeed, notice if $\phi$ is an $H^1$ solution of \eqref{E:EJ_profile}
then differentiating \eqref{E:EJ_profile} yields
\begin{equation} \label{E:EJ_smoothness_argument}
\CalK \phi' = \left( \frac{3}{2}\phi^2 - 3c\phi + c^2 \right) \phi' \implies \phi' = \frac{\CalK \phi'}{\frac{3}{2}\phi^2 - 3c\phi + c^2}.
\end{equation}
The operator $\CalK$ improves the smoothness of its operand by exactly one order, hence by \eqref{E:EJ_smoothness_argument} we have that $\phi'$ is as smooth as $\phi$ so long as $\frac{3}{2}\phi^2 - 3c\phi + c^2 \neq 0$, in which case a bootstrap argument demonstrates that $\phi$ is in fact $C^\infty$.  
When continuing from constant solutions, i.e. solutions with height zero, the breakdown of smoothness occurs when $\frac{3}{2}\phi^2 - 3c\phi + c^2 = 0$ or, more precisely, 
when the bifurcation branch intersects the curve $\phi = c \left( 1-\frac{1}{\sqrt{3}} \right)$.  
Thus we expect (as is shown in \cite{EJ16}) that periodic solutions of \eqref{E:EJ_profile} along the bifurcation branch with fixed period will have a maximum amplitude of 
\[
\max\phi(x) = \gamma := c\left( 1-\frac{1}{\sqrt{3}} \right).
\]
See \cite{EJ16} for details of the above arguments.

\begin{figure}[t]
    \centering
    
    \begin{subfigure}[t]{0.5\linewidth}
    \centering
	\includegraphics[scale=0.4]{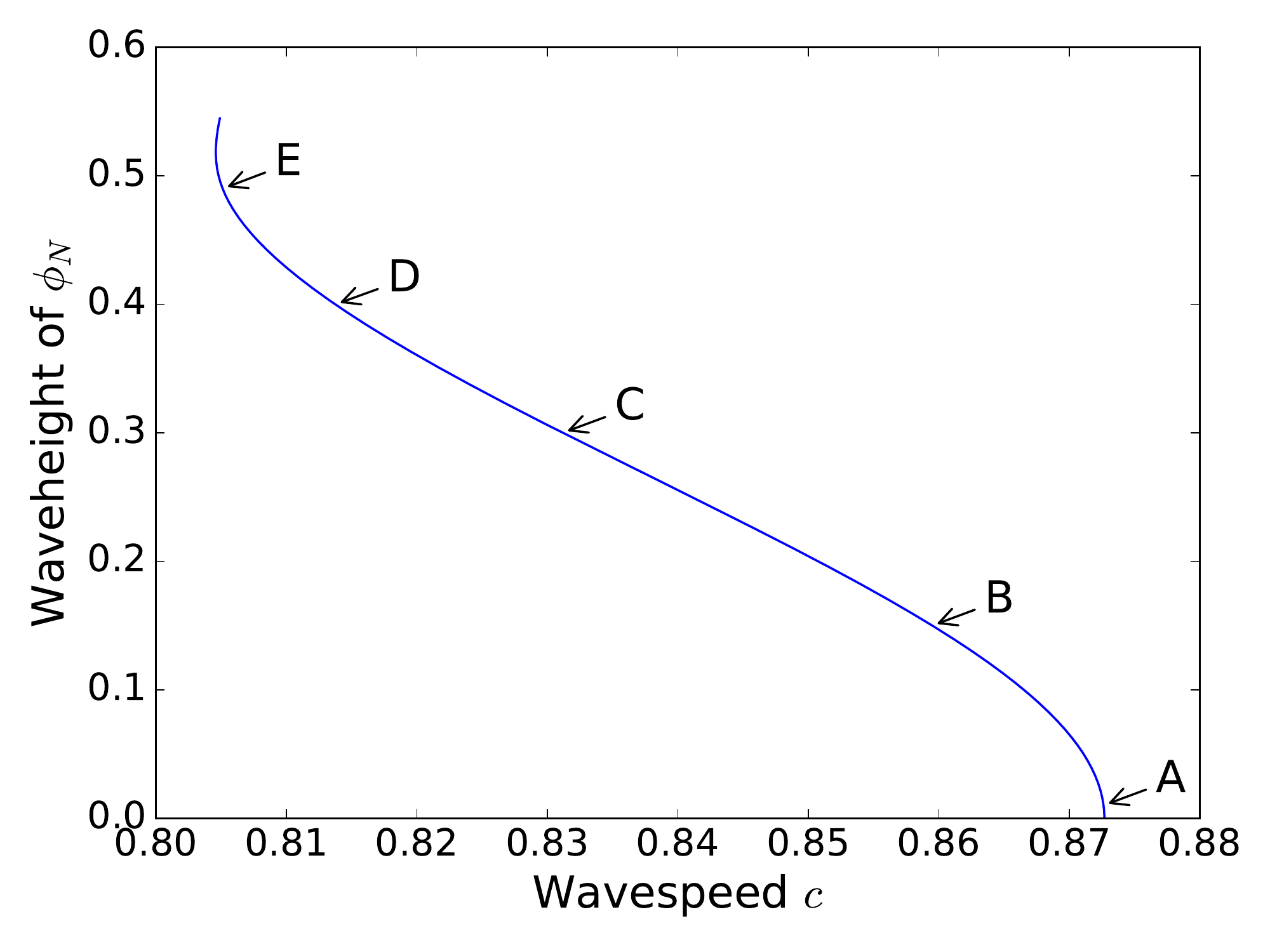}
	\caption{}
	\end{subfigure}%
	\begin{subfigure}[t]{0.45\linewidth}
    \centering
	\includegraphics[scale=0.4]{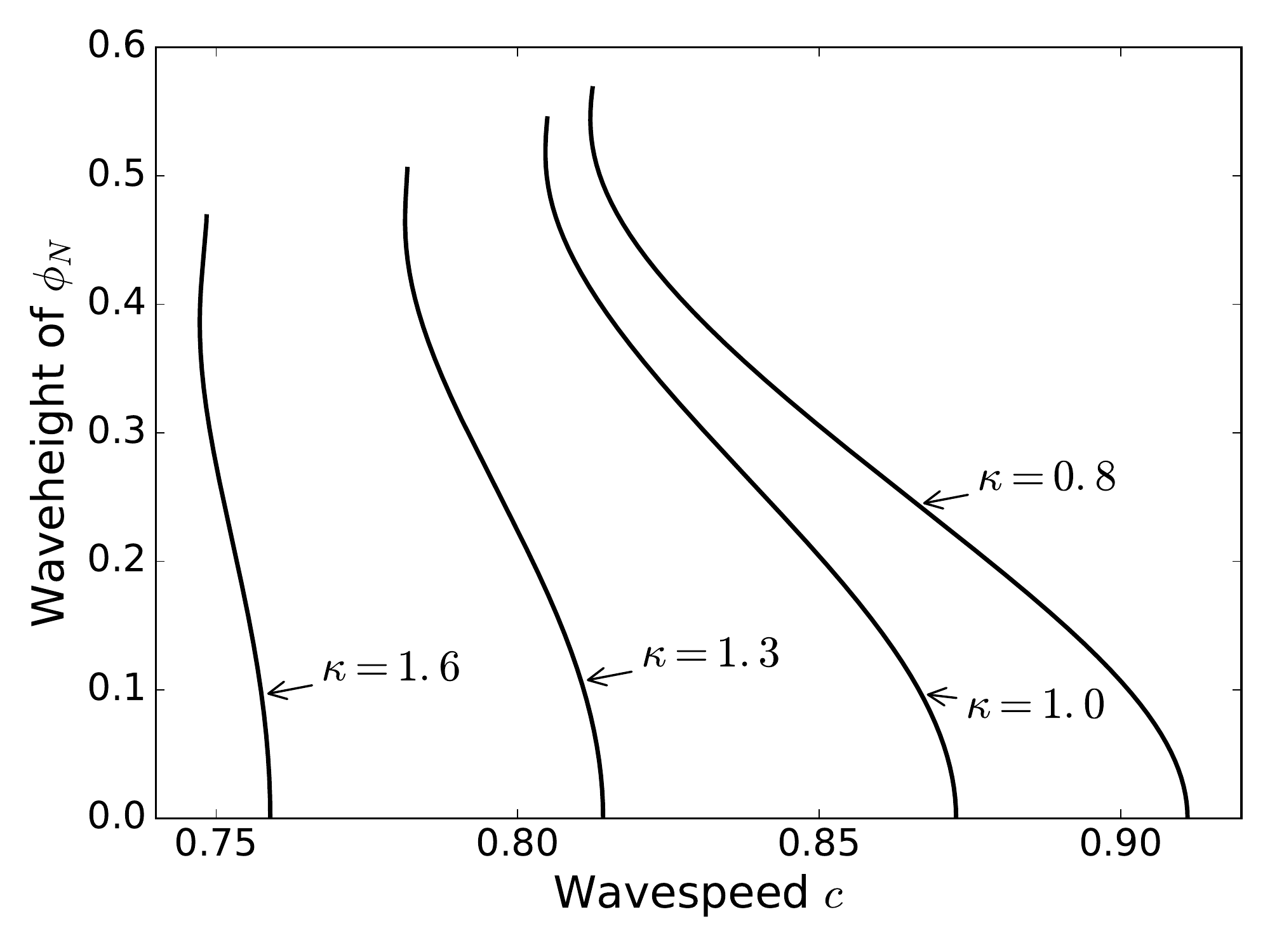}
	\caption{}
	\end{subfigure}
           \caption{(a) A numerical approximation of the global bifurcation branch of $2\pi$-periodic ($\kappa=1$), even traveling solutions of \eqref{E:EJ_profile} is displayed, with specific points A-E labeled for forthcoming computations: Point A at $c \approx 0.8726$, height $\approx 0.01$; Point B at $c \approx 0.8595$, height $\approx 0.15$; Point C at $c \approx 0.8312$, height $\approx 0.30$; Point D at $c \approx 0.8138$, height $\approx 0.40$; Point E at $c \approx 0.8051$, height $\approx 0.49$.  (b) Bifurcation branches of $2\pi/\kappa$-periodic solutions of \eqref{E:EJ_profile} for varying values of $\kappa$.
           Notice all these curves experience a turning point near the top of the branch.}
            \label{F:EJ_global_bif_diagram}
\end{figure}

\subsubsection{Analysis of Waves Bifurcating from Zero}\label{S:bf0}
We begin by studying  the bifurcation of $2\pi/\kappa$-periodic traveling wave solutions of \eqref{E:EJ_profile} bifurcating from the trivial state $\phi=0$.
To this end, using the profile equation \eqref{E:EJ_profile} in combination with the following local bifurcation formulas (see \cite[Proposition 5.1]{EJ16})
\begin{align}
\phi(x;\varepsilon) &:= \varepsilon \cos(\kappa x) + \frac{3c_\kappa \varepsilon^2}{4} \left( \frac{1}{c_\kappa^2-1} + \frac{\cos(2\kappa x)}{c_\kappa^2-c_{2\kappa}^2} \right) + O(\varepsilon^3)  \label{E:EJ_local_bif_formula_phi} \\
c(\varepsilon) &:= c_\kappa + \frac{3\varepsilon^2}{8}\left[-\frac{1}{2c_\kappa} + 3c_\kappa\left( \frac{1}{c_\kappa^2-1} + \frac{1}{2(c_\kappa^2-c_{2\kappa}^2)} \right) \right], 
\quad c_\ell := \sqrt{\frac{\tanh(\ell)}{\ell}},  \label{E:EJ_local_bif_formula_c}
\end{align}
we apply the pseudo-arclength method as discussed in Section \ref{SS:pseudoarclength_method} to compute approximations $\phi_N(x_i)$ of the profile $\phi(x_i)$ at the collocation points $x_i = \dfrac{(2i-1)\pi}{\kappa N}$, continuing while $\phi(0) \approx \phi_N(x_1)  < \gamma$.  A global bifurcation plot of waveheight vs. wavespeed for various wavenumbers $\kappa$ is shown in Figure \ref{F:EJ_global_bif_diagram}(b), where we define
\[
\text{waveheight} := \max_{x \in [0,\pi/\kappa]} \phi(x) - \min_{x \in [0,\pi/\kappa]} \phi(x) = \phi(0)-\phi(\pi/\kappa),
\]
which is well-defined due to the monotonicity properties of even solutions along the global bifurcation branch as demonstrated in \cite{EJ16}.

\begin{figure}[t]
    \centering
    
    \begin{subfigure}[t]{0.5\linewidth}
    \centering
	\includegraphics[scale=0.4]{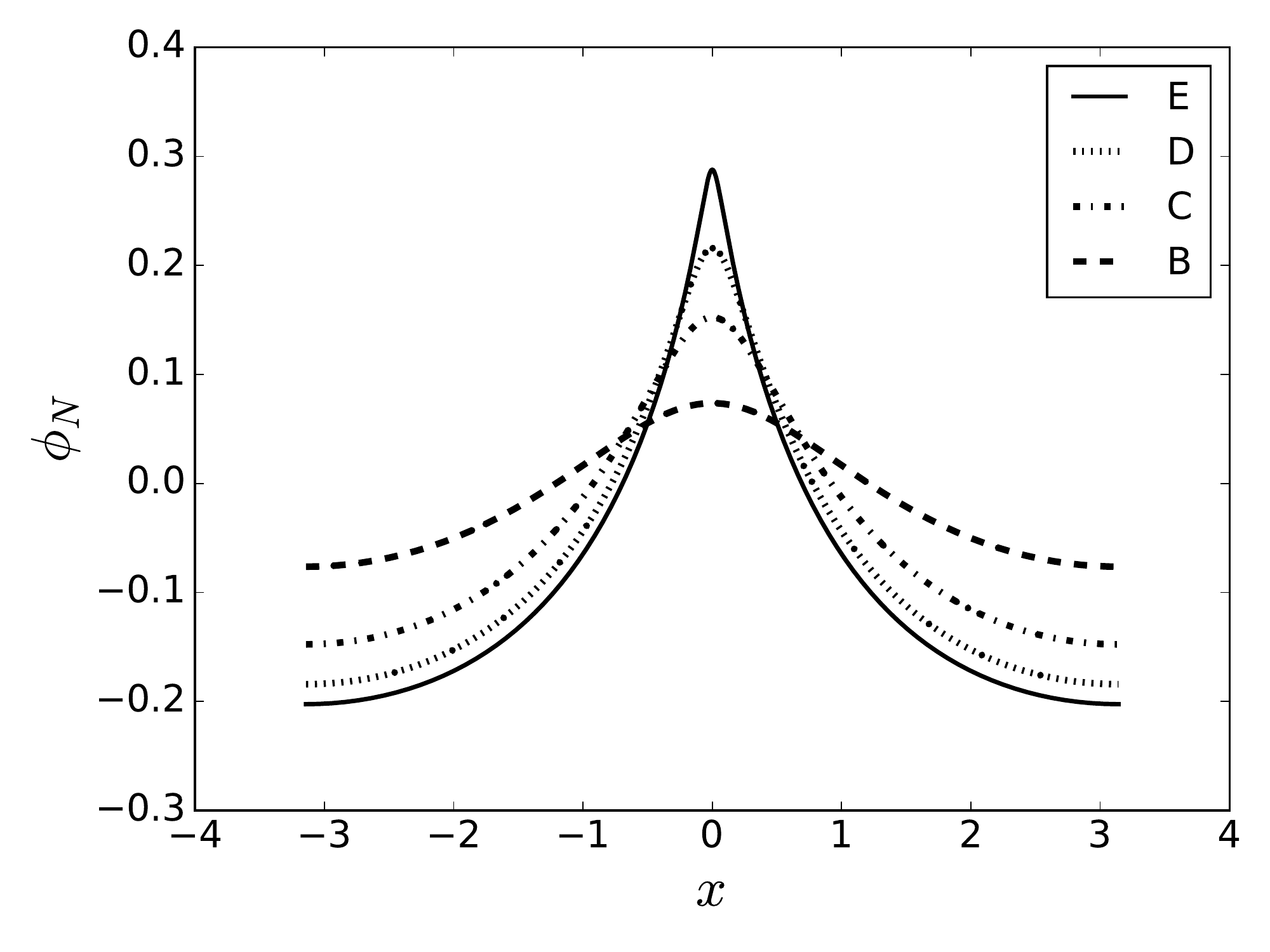}
	\caption{}
	\end{subfigure}%
	\begin{subfigure}[t]{0.5\linewidth}
    \centering
	\includegraphics[scale=0.4]{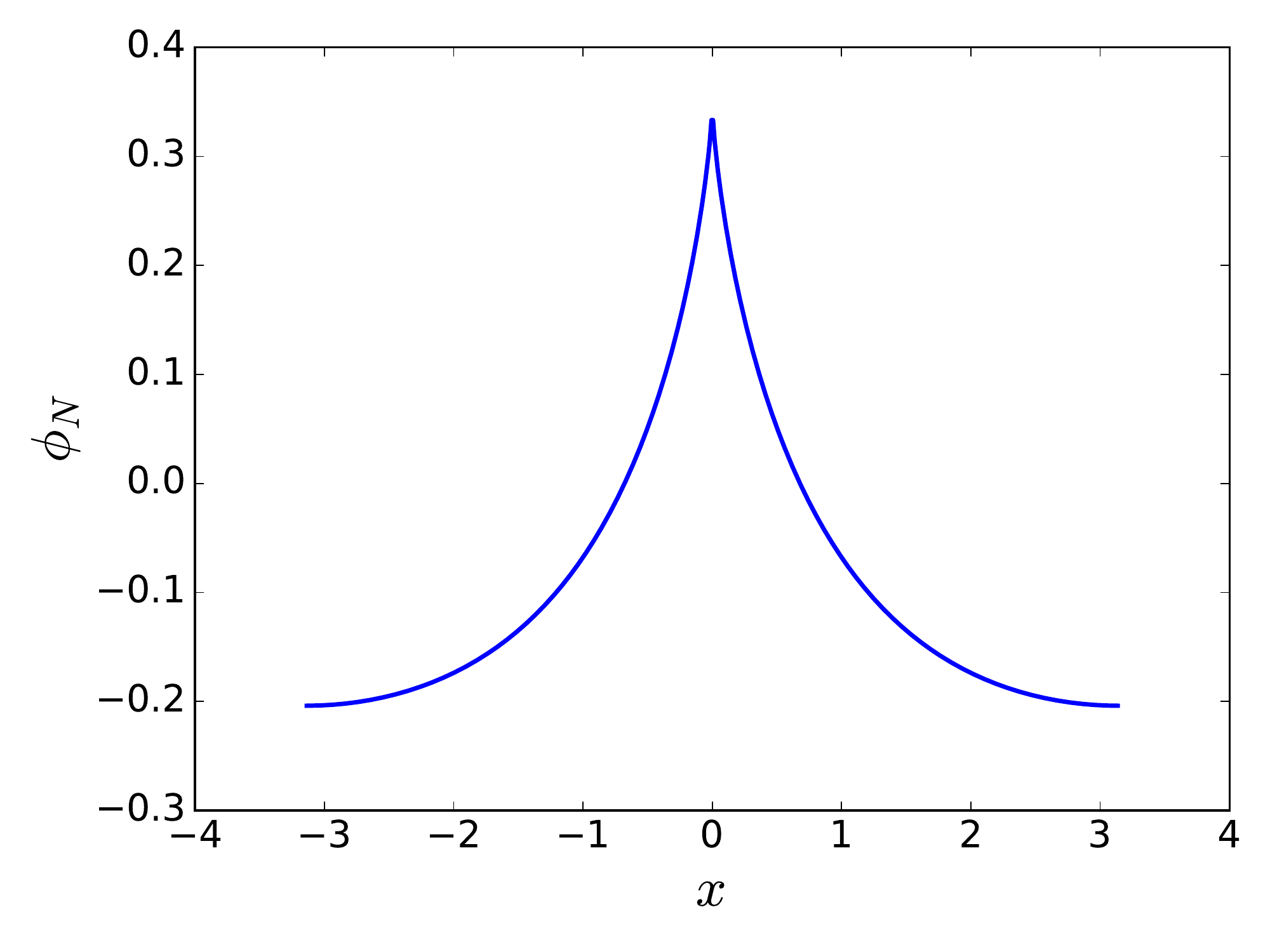}
	\caption{}
	\end{subfigure}
	 \caption{(a) Profiles corresponding to the sampled points on the branch of $2\pi$-periodic solutions of in Figure \ref{F:EJ_global_bif_diagram}(a): Point B at $c \approx 0.8595$, height $\approx 0.15$; Point C at $c \approx 0.8312$, height $\approx 0.30$; Point D at $c \approx 0.8138$, height $\approx 0.40$; Point E at $c \approx 0.8051$, height $\approx 0.49$.  (b) A nearly ``highest" wave high up the bifurcation branch in Figure \ref{F:EJ_global_bif_diagram}(a).}
    \label{F:EJ_profiles}
\end{figure}

To illustrate the our numerical results for this model, we sample a selection of points along the global bifurcation diagram as described in  Figure \ref{F:EJ_global_bif_diagram}(a).
As shown in Figure \ref{F:EJ_profiles}(a), the crest at the maximum of the smooth profiles along the global bifurcation branch becomes sharper as the waveheight increases, appearing
to converge to a non-trivial profile with a singularity at the origin.  The existence and qualitative properties of this highest singular wave has been studied analytically
in \cite{EJ16}, where it was shown the highest wave has a logarithmic cusp of order $\abs{x\ln\abs{x}}$ near the top of the bifurcation curve: see Figure \ref{F:EJ_profiles}(b).
We refer the interested reader to \cite{EJ16} for details.

In addition to the existence of a highest singular wave, we also wish to understand the dynamical stability of the smooth periodic solutions constructed above.  To this end, we wish to determine the spectrum of the 
linearization of \eqref{E:EJ_traveling_wave_model} about such a smooth periodic traveling wave $\phi$.
To linearize the traveling wave system \eqref{E:EJ_traveling_wave_model} about the velocity profile $\phi$ described above and the corresponding height profile $\psi = c\phi - \dfrac{1}{2}\phi^2$, we write
\begin{equation} \label{E:EJ_linearization_expansion}
\begin{array}{rcl}
u(x,t) &=& \phi(x) + \varepsilon v(x,t) + \CalO(\varepsilon^2) \\
\eta(x,t) &=& \psi(x) + \varepsilon w(x,t) + \CalO(\varepsilon^2)
\end{array}
\end{equation}
where here, since we are interested in the stability to \emph{localized} perturbations, we require that $v(\cdot,t),w(\cdot,t)\in L^2(\mathbb{R})$.
for each time $t>0$ for which they are defined.
Substituting these expansions into the first equation of \eqref{E:EJ_traveling_wave_model} yields
\begin{align}
\varepsilon v_t &= (-\psi' - \phi\phi' + c\phi') + \varepsilon \left[ -\phi' v + (c-\phi)v_x - w_x \right] + \CalO(\varepsilon^2) \nonumber \\
\varepsilon w_t &= (-\CalK\phi' - (\psi\phi)' + c\psi') + \varepsilon \left[ -\psi' v - \psi v_x - \CalK v_x - \psi' w + (c-\phi)w \right] + \CalO(\varepsilon^2)
\end{align}
as $\varepsilon \to 0$, where the $\CalO(1)$ terms vanish per the equilibrium system \eqref{E:EJ_traveling_profile_system}.  
Taking $\varepsilon\to 0$ above and decomposing the perturbations as $(v(x,t),w(x,t))=e^{\lambda t}(v(x),w(x))$ yields the spectral problem
\begin{equation}\label{E:EJ_linearize}
\left\{\begin{aligned}
\lambda v&= -\phi' v + (c-\phi)v_x - w_x  \\
\lambda w&= -\psi' v - \psi v_x - \CalK v_x - \psi' w + (c-\phi)w\\
&\qquad=:\mathcal{L}\left[\begin{array}{c}v\\w\end{array}\right]
\end{aligned}\right.
\end{equation}
considered on $L^2(\mathbb{R})\times L^2(\mathbb{R})$.

\begin{figure}[t]
    \centering
    
    \begin{subfigure}[t]{0.33\linewidth}
	\includegraphics[scale=0.26]{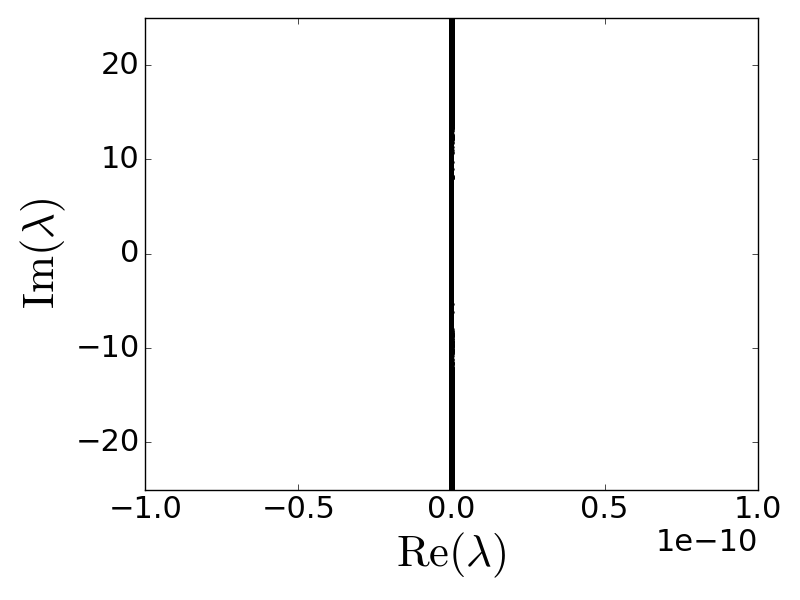}
	\caption{}
	\end{subfigure}%
    \begin{subfigure}[t]{0.33\linewidth}
	\includegraphics[scale=0.26]{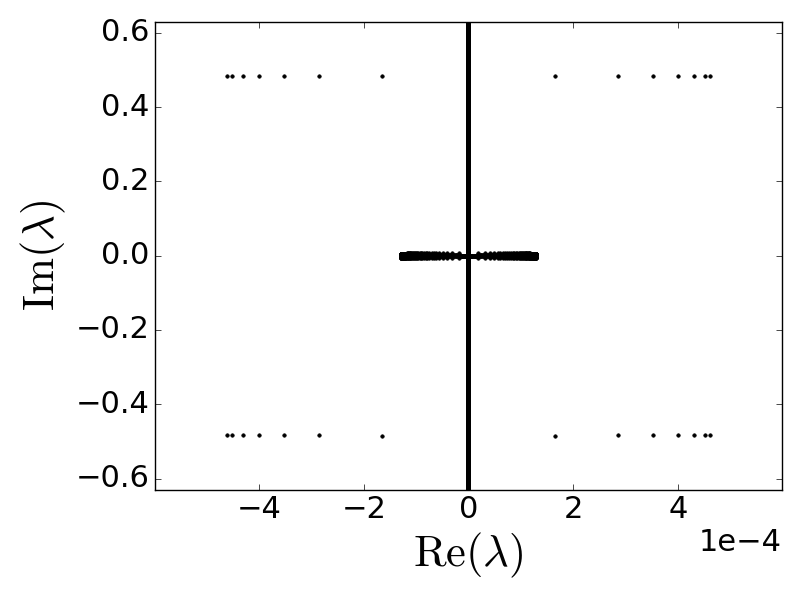}
	\caption{}
	\end{subfigure}%
	\begin{subfigure}[t]{0.33\linewidth}
	\includegraphics[scale=0.26]{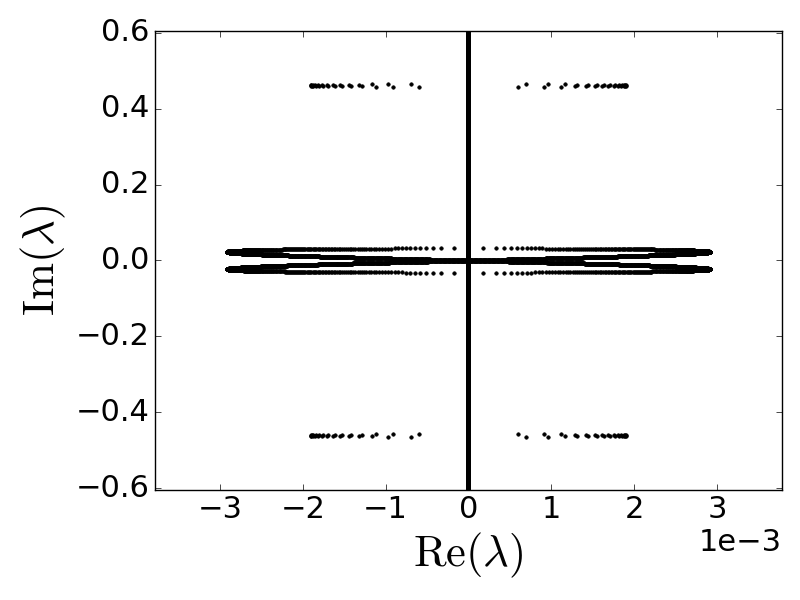}
	\caption{}
	\end{subfigure}%
	
    \begin{subfigure}[t]{0.33\linewidth}
	\includegraphics[scale=0.26]{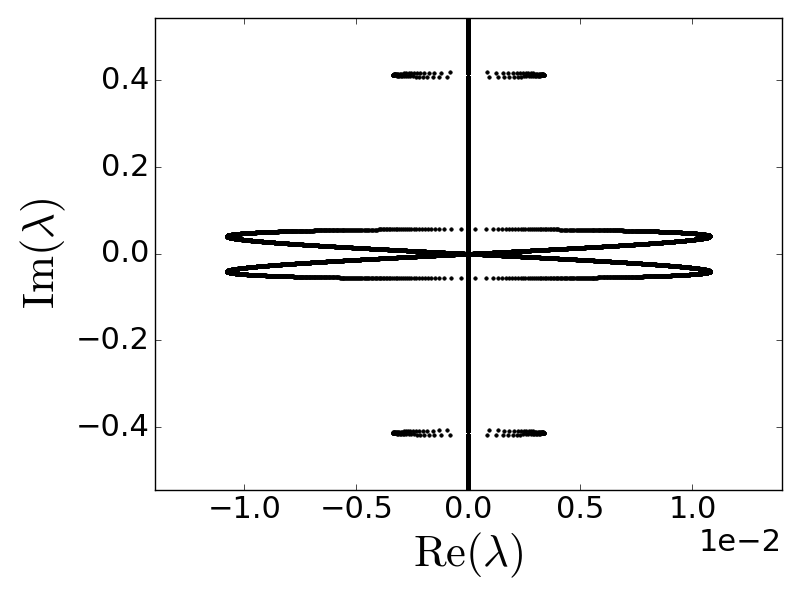}
	\caption{}
	\end{subfigure}%
    \begin{subfigure}[t]{0.33\linewidth}
	\includegraphics[scale=0.26]{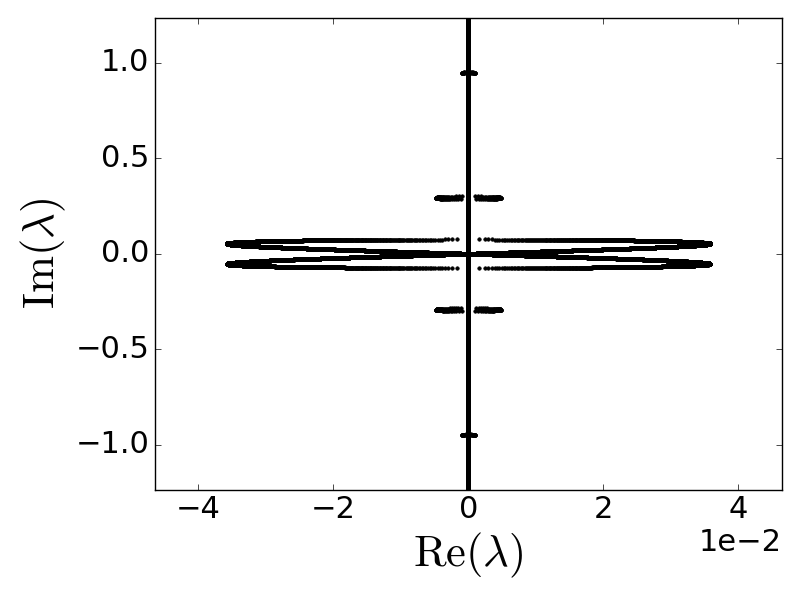}
	\caption{}
	\end{subfigure}%
    \begin{subfigure}[t]{0.33\linewidth}
	\includegraphics[scale=0.26]{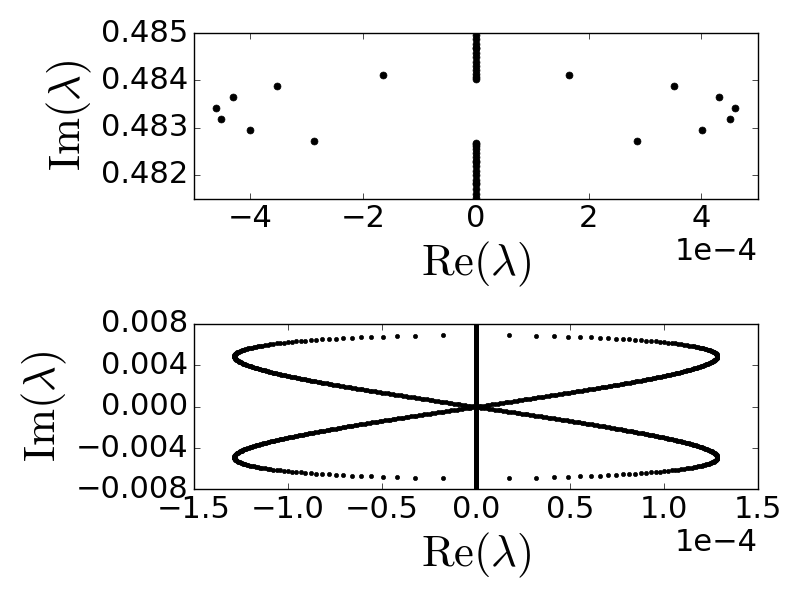}
	\caption{}
	\end{subfigure}

	\caption{(a)-(e) Numerical approximations of the spectrum of $\mathcal{L}$ corresponding to the linearization of \eqref{E:EJ} about the profiles corresponding to the points
A-E, respectively, along the bifurcation branch for $2\pi$-periodic solutions in Figure \ref{F:EJ_global_bif_diagram}(a).  
(f) Zoom-in on the high-frequency instability (top) and the modulational instability (bottom) in the spectral plot (b).}
	\label{F:EJ_spectrum_plots}
\end{figure}

\begin{figure}[t]
    \centering
    
    \begin{subfigure}[t]{0.33\linewidth}
    \centering
	\includegraphics[scale=0.28]{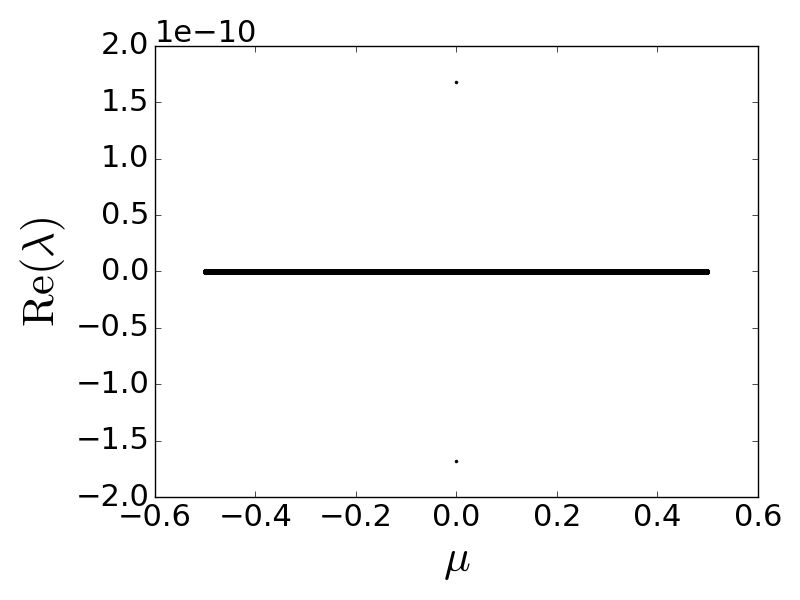}
	\caption{}
	\end{subfigure}%
    \begin{subfigure}[t]{0.33\linewidth}
    \centering
	\includegraphics[scale=0.28]{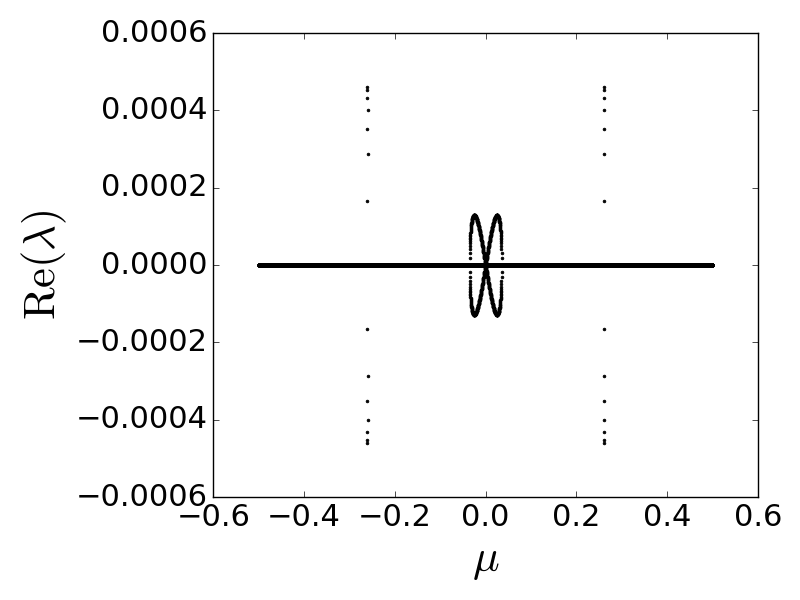}
	\caption{}
	\end{subfigure}%
	\begin{subfigure}[t]{0.33\linewidth}
	\centering
	\includegraphics[scale=0.28]{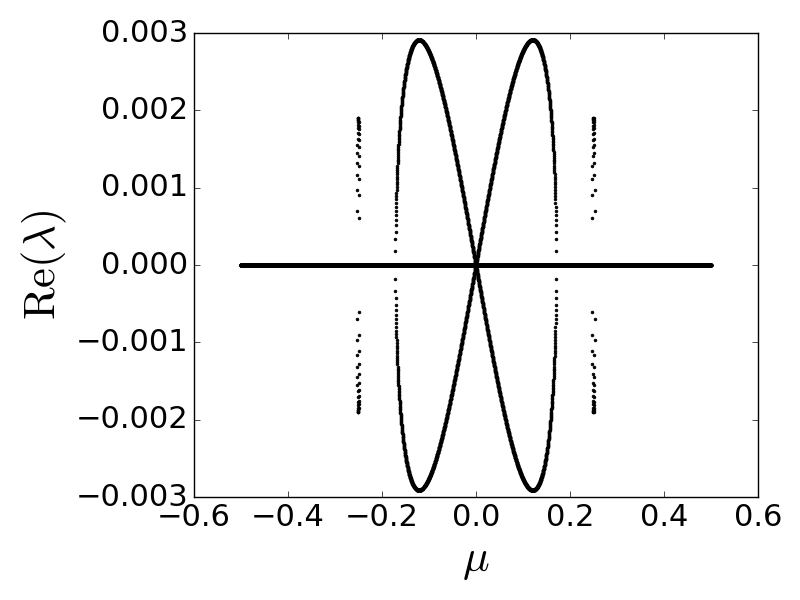}
	\caption{}
	\end{subfigure}%
	
    \begin{subfigure}[t]{0.33\linewidth}
    \centering
	\includegraphics[scale=0.28]{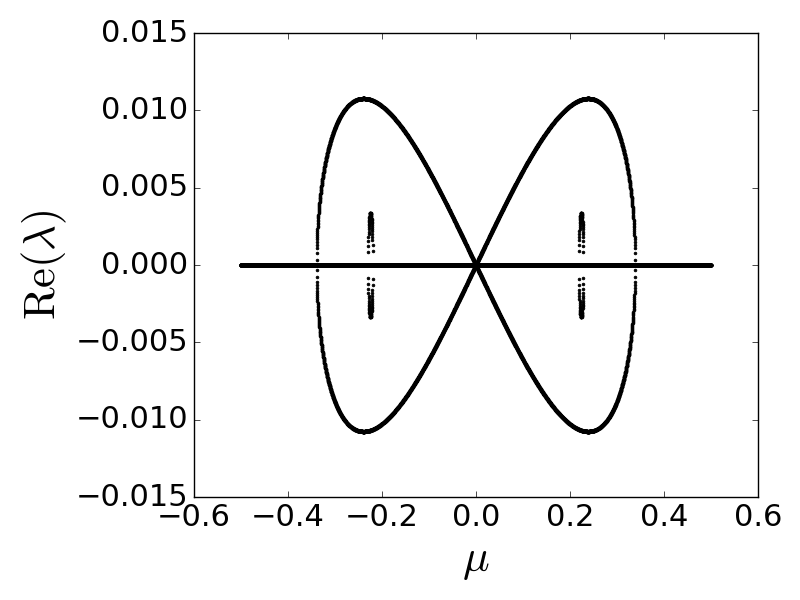}
	\caption{}
	\end{subfigure}%
	\hspace{1em}
    \begin{subfigure}[t]{0.33\linewidth}
    \centering
	\includegraphics[scale=0.28]{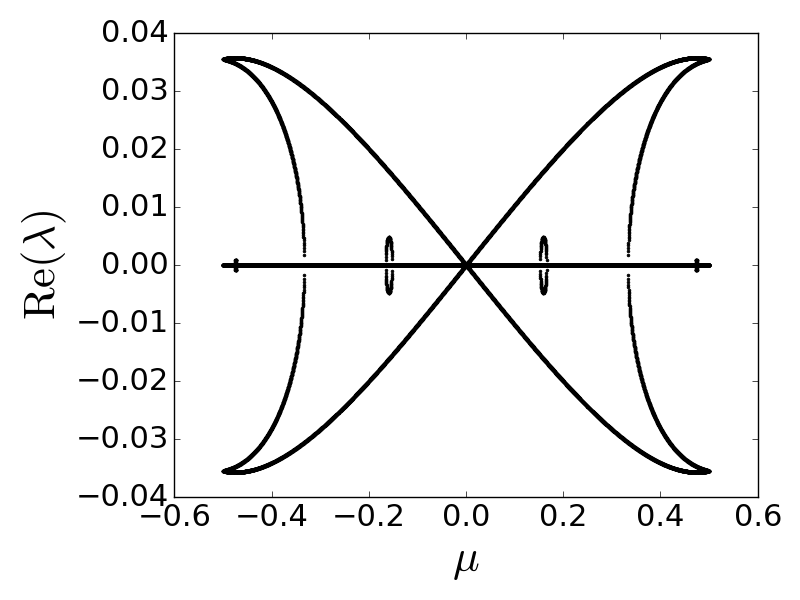}
	\caption{}
	\end{subfigure}

	\caption{(a)-(e) Plots of the growth rates $\Re(\lambda)$ vs. $\mu$ for each of the spectral plots (a)-(e), respectively, in Figure \ref{F:EJ_spectrum_plots}.}
	\label{F:EJ_Re_vs_mu}
\end{figure}

By Floquet theory, we know the spectrum associated to $\mathcal{L}$ is purely essential, containing no isolated eigenvalues of finite multiplicity: see \cite{RS4,BH}.
Indeed, one can show that $\lambda$ belongs to the $L^2(\mathbb{R})$-spectrum of $\mathcal{L}$ if and only if there exists a so-called ``Bloch parameter" $\mu\in[0,1)$ such that
\eqref{E:EJ_linearize} has a bounded solution satisfying
\[
(v,w)(x+2\pi/\kappa) = e^{i\kappa\mu}(v,w)(x)
\] 
for all $x\in\mathbb{R}$.  It follows that the perturbations $(v,w)$ in \eqref{E:EJ_linearize} can be expanded in ``Bloch" form as 
\begin{equation} \label{E:EJ_linearize_floquet}
\begin{array}{rcl}
v(x) =\displaystyle \sum_{l \in \Z} \widehat{V}(l)e^{i\kappa(\mu+l)x},\quad
w(x)=\displaystyle \sum_{l \in \Z} \widehat{W}(l)e^{i\kappa(\mu+l) x}.
\end{array}
\end{equation}
Substituting the expansions \eqref{E:EJ_linearize_floquet} for $v$ and $w$ into \eqref{E:EJ_linearize}(i), we obtain
\begin{align*}
\lambda \sum_{m \in \Z} \widehat{V}(m) e^{i\kappa mx} &= \sum_{m \in \Z} \sum_{l \in \Z} \left[ \left( ic\kappa(\mu+l)\delta_{m,l} - i\kappa(\mu+m)\widehat{\phi}(m-l) \right) \widehat{V}(l) - i\kappa(\mu+l)\widehat{W}(l) \right] e^{i\kappa mx},
\end{align*}
where here $\widehat{\phi}(n)$ denotes the $n$th Fourier coefficient associated to the periodic profile $\phi$.
It follows that each $m \in \Z$ we must have
\begin{equation} \label{E:VHat_eigenvalue_problem}
\lambda \widehat{V}(m) = \sum_{l \in \Z} \left[ \widehat{A}^\mu(m,l)\widehat{V}(l) + \widehat{B}^\mu(m,l)\widehat{W}(l) \right],
\end{equation}
where
\begin{align*}
\widehat{A}^\mu(m,l) &= ic\kappa(\mu+l)\delta_{m,l} - i\kappa(\mu+m)\widehat{\phi}(m-l) \\
\widehat{B}^\mu(m,l) &= -i\kappa(\mu+l)\delta_{m,l}
\end{align*}
and $\delta_{m,n}$ is the Kronecker delta.
Similarly, substituting the expansions \eqref{E:EJ_linearize_floquet} into \eqref{E:EJ_linearize}(ii) yields, for each $m \in \Z$,
\begin{equation} \label{E:WHat_eigenvalue_problem}
\lambda \widehat{W}(m) = \sum_{l \in \Z} \left[ \widehat{C}^\mu(m,l)\widehat{V}(l) + \widehat{D}^\mu(m,l)\widehat{W}(l) \right],
\end{equation}
where
\begin{align*}
\widehat{C}^\mu(m,l) &= -i\kappa(\mu+m)\widehat{\psi}(m-l) - i\kappa(\mu+l)\widehat{\mathcal{K}}(\kappa(\mu+l))\delta_{m,l} \\
\widehat{D}^\mu(m,l) &= ic\kappa(\mu+l)\delta_{m,l} - i\kappa(\mu+m)\widehat{\phi}(m-l).
\end{align*}
Defining the bi-infinite matrices 
\[
\widehat{A}^\mu := [ \widehat{A}^\mu(m,l) ]_{m,l \in \Z}, \quad \widehat{B}^\mu := [ \widehat{B}^\mu(m,l) ]_{m,l \in \Z}, \quad \widehat{C}^\mu := [ \widehat{C}^\mu(m,l) ]_{m,l \in \Z}, \quad \widehat{D}^\mu := [ \widehat{D}^\mu(m,l) ]_{m,l \in \Z}.
\]
entry-wise for row-$m$ and column-$l$ with $m,l \in \Z$, the system \eqref{E:VHat_eigenvalue_problem}-\eqref{E:WHat_eigenvalue_problem} can be written
in block bi-infinite matrix form as
\begin{equation}  \label{E:FFHM_eigenvalue_problem}
\def\arraystretch{1.5}  
\lambda \left[ \begin{array}{c}
\widetilde{V} \\
\hline
\widetilde{W}
\end{array} \right]
= \left[
\begin{array}{c|c}
\widehat{A}^\mu & \widehat{B}^\mu \\
\hline
\widehat{C}^\mu & \widehat{D}^\mu 
\end{array}
\right]
\left[ \begin{array}{c}
\widetilde{V} \\
\hline
\widetilde{W}
\end{array} \right]
=: \widehat{\CalL}^\mu 
\left[ \begin{array}{c}
\widetilde{V} \\
\hline
\widetilde{W}
\end{array}
\right],
\end{equation}
where $\widetilde{V},\widetilde{W}$ denote the bi-infinite arrays
\[
\left\{\begin{aligned}
\widetilde{V} &:= \begin{bmatrix} \ldots & \widetilde{V}(-2) & \widetilde{V}(-1) & \widetilde{V}(0) & \widetilde{V}(1) & \widetilde{V}(2) & \ldots \end{bmatrix}^T \\
\widetilde{W} &:= \begin{bmatrix} \ldots & \widetilde{W}(-2) & \widetilde{W}(-1) & \widetilde{W}(0) & \widetilde{W}(1) & \widetilde{W}(2) & \ldots \end{bmatrix}^T.
\end{aligned}\right.
\]
For each $\mu$, the spectrum of $\widehat{\mathcal{L}}^\mu$ consists of countably many discrete eigenvalues with finite multiplicity and, furthermore, 
by the above considerations we have the spectral decomposition
\[
\sigma_{L^2(\mathbb{R})}(\CalL) = \bigcup_{\mu \in [0,1)} \sigma(\widehat{\CalL}^\mu).
\]

We numerically approximate the spectrum of the bi-infinite matrix $\widehat{\mathcal{L}}^\mu$
by taking a sequence of $\mu$ in a finite discretization of $[0,1)$, truncating each block 
of $\widehat{\CalL}^\mu$ in \eqref{E:FFHM_eigenvalue_problem} to finite dimension, and computing the eigenvalues of the truncated 
matrix using a standard matrix eigenvalue solver.  See Appendix \ref{A:general_FFHM} for further discussion.
See Figure \ref{F:EJ_spectrum_plots} for plots of the spectrum at the sampled points on the bifurcation diagram of $2\pi$-periodic solutions in 
Figure \ref{F:EJ_global_bif_diagram}(a), as well as Figure \ref{F:EJ_Re_vs_mu} for plots of the growth rate $\Re(\lambda)$ vs. the Bloch parameter $\mu$.  
It is known that the full Euler equations exhibit high-frequency instabilities (visually characterized by ``bubbles" of spectrum emanating from the 
imaginary axis from points away from the origin) in small-amplitude periodic traveling waves; however, we discover that while this model 
demonstrates both high-frequency and modulational instabilities for waves of sufficiently high amplitude, a high-frequency instability is 
numerically undetectable for small amplitude waves: see Figure
\ref{F:EJ_spectrum_plots}(a) and Figure \ref{F:EJ_Re_vs_mu}(a) below.
Whether such instabilities are actually nonexistent or simply not detected for small amplitude waves in our numerics is unclear.

\begin{figure}[t]
    \centering
    
    \begin{subfigure}[t]{0.5\linewidth}
    \centering
	\includegraphics[scale=0.4]{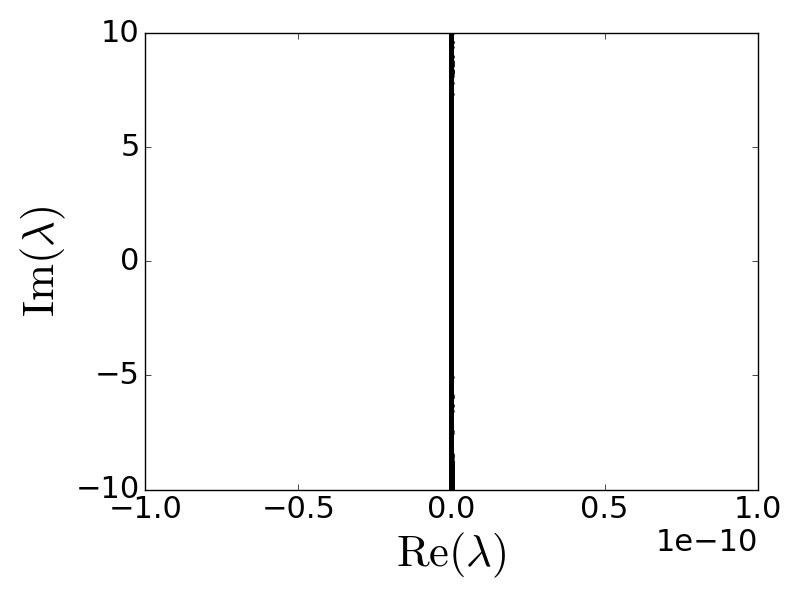}
	\caption{}
	\end{subfigure}%
	\begin{subfigure}[t]{0.5\linewidth}
	\centering
	\includegraphics[scale=0.4]{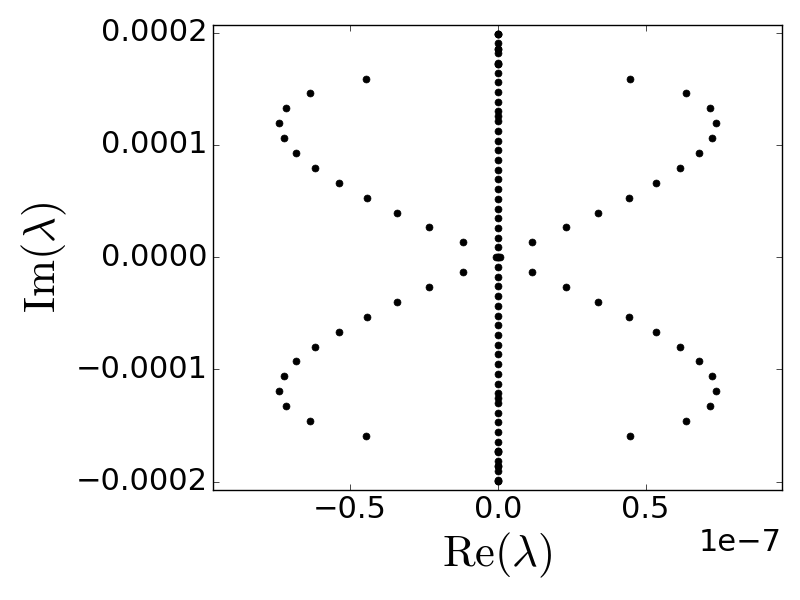}
	\caption{}
	\end{subfigure}
	
    \caption{Evidence of the Benjamin-Feir instability in \eqref{E:EJ} for waveheight $0.01766$.  In(a) we have $\kappa= 1.005$ while in (b) we have $\kappa=1.008$.}
    \label{F:EJ_critical_kappa_comparison}
\end{figure}

Solutions of different wavelength can be generated by varying the value of $\kappa$, as shown in the plot of the bifurcation diagrams for various $\kappa$ in Figure \ref{F:EJ_global_bif_diagram}(b).  In \cite{HP_shallow}, an analysis of asymptotically small amplitude waves shows that there exists a critical value $\kappa_c$ such that all waves are modulationally unstable for all $\kappa > \kappa_c$ where $\kappa_c\approx 1.008$.  
Inspired by this analysis, we performed numerical experiments varying the value of $\kappa$ for small, but fixed, waveheights, finding very good agreement with the results in \cite{HP_shallow}.
Figure \ref{F:EJ_critical_kappa_comparison}(a) demonstrates that waves with waveheight $0.01766$ and corresponding to $\kappa=1.005$ are modulationally stable,
while for $\kappa=1.008$ a modulational instability appears for waves of the same height: see Figure \ref{F:EJ_critical_kappa_comparison}(b).


\begin{figure}[t]
    \centering
    
    \begin{subfigure}[t]{0.5\linewidth}
    \centering
	\includegraphics[scale=0.36]{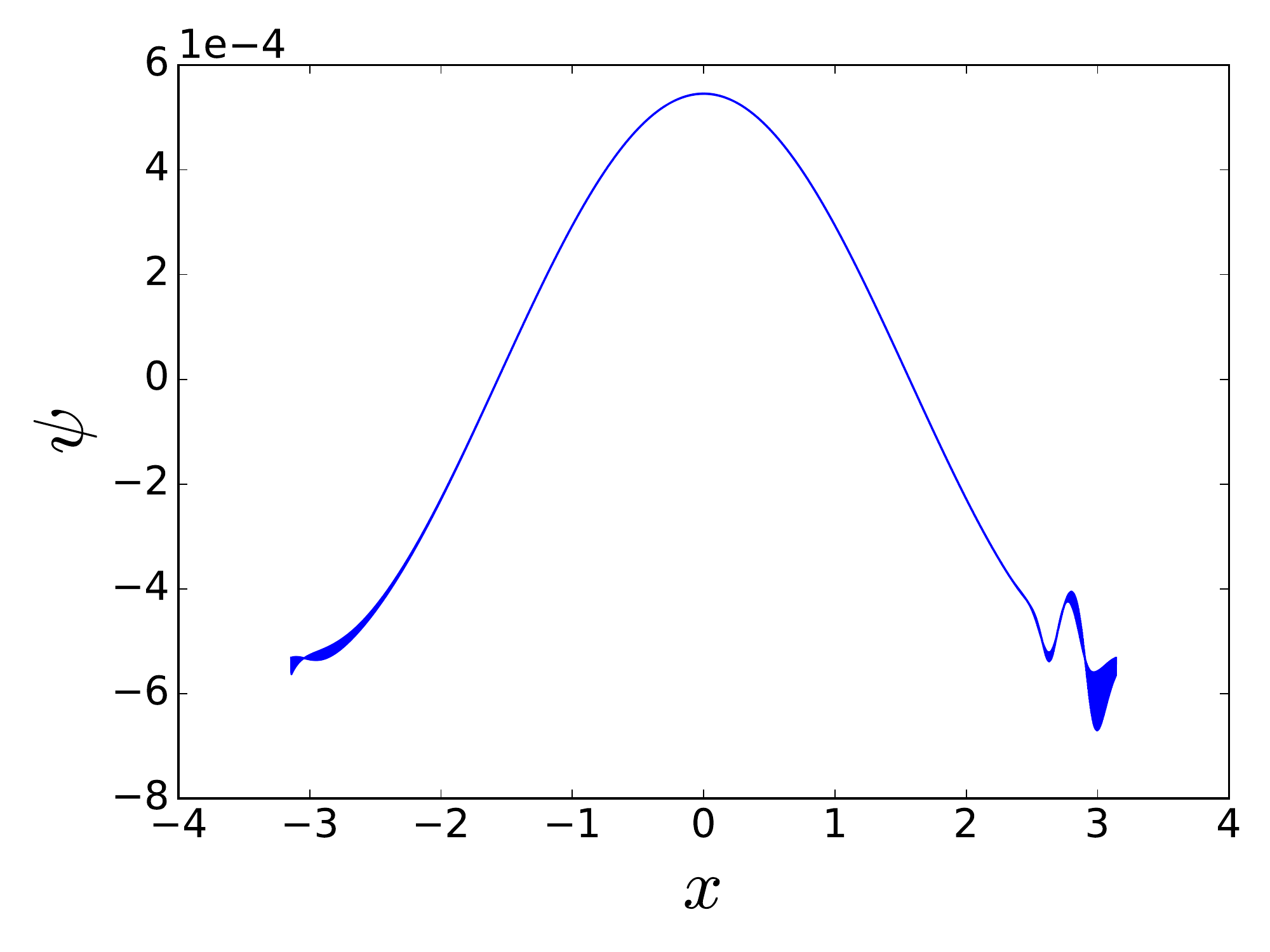}
	\caption{}
	\end{subfigure}%
	\begin{subfigure}[t]{0.5\linewidth}
	\centering
	\includegraphics[scale=0.36]{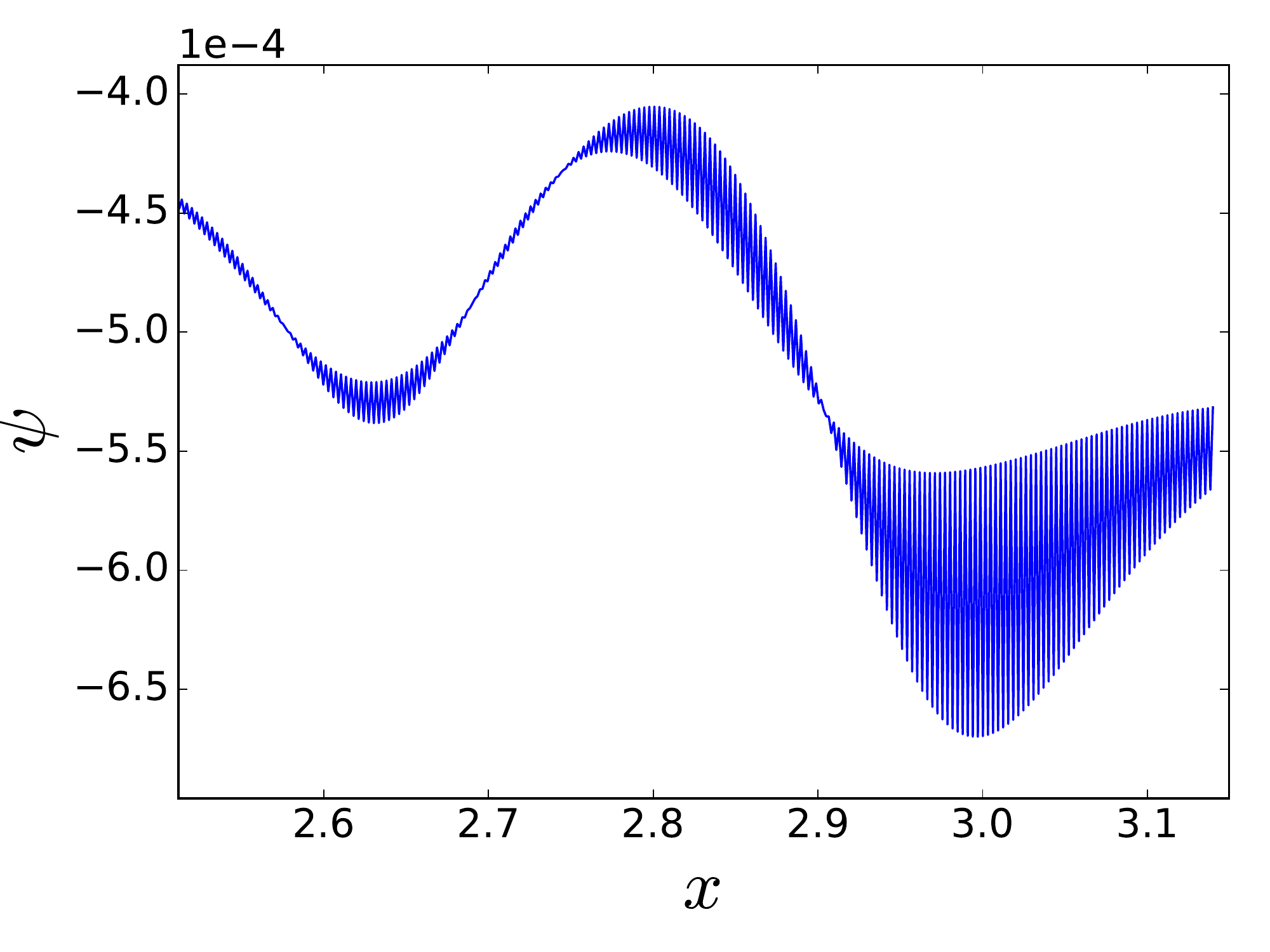}
	\caption{}
	\end{subfigure}
	
    \begin{subfigure}[t]{0.5\linewidth}
    \centering
	\includegraphics[scale=0.36]{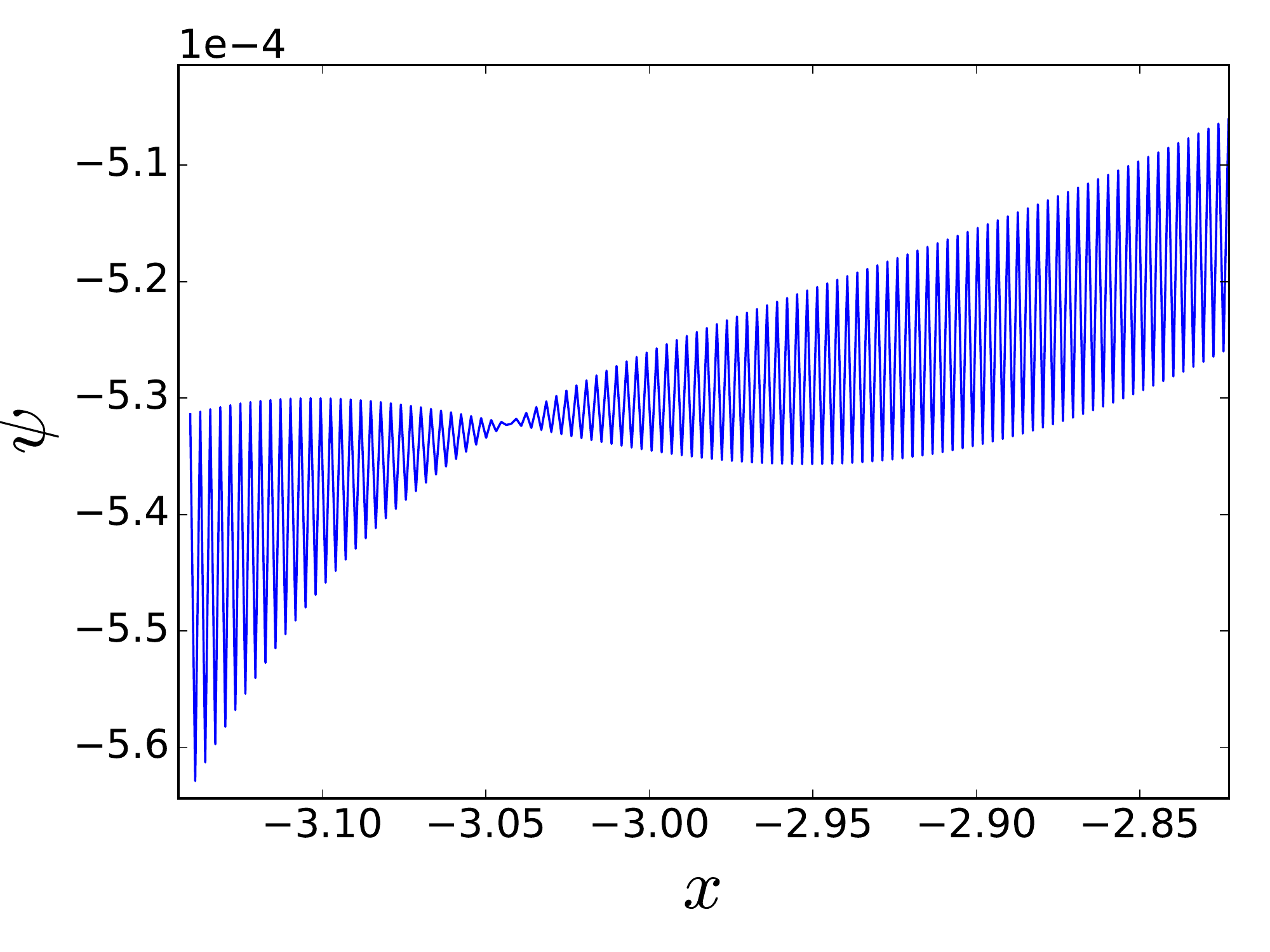}
	\caption{}
	\end{subfigure}%
	\begin{subfigure}[t]{0.5\linewidth}
	\centering
	\includegraphics[scale=0.36]{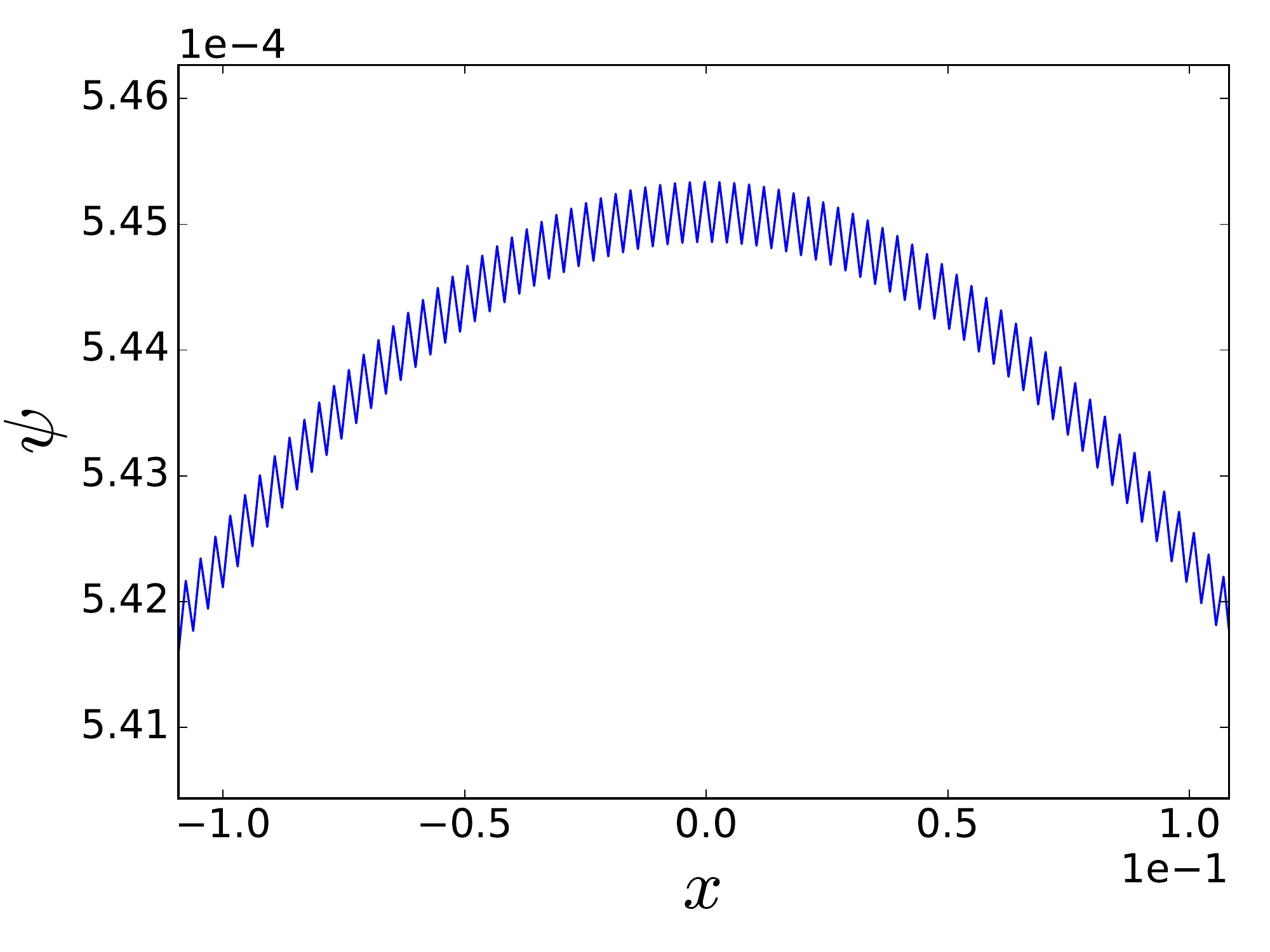}
	\caption{}
	\end{subfigure}	
	
	\caption{Numerical evidence of ill-posedness of \eqref{E:EJ} for initial data $\eta(x,0)=\psi(x)$ having negative infimum.  (a) The result of a short time evolution of a $2\pi$-periodic
wave ($\kappa=1$) with small waveheight $\max\psi - \min\psi \approx 0.00109$, $c \approx 0.872693$.  We zoom in near the rightmost and leftmost of the waves in (b) an (c), respectively.  In (d) we observe that oscillations also are forming
near the positive maxima, although the oscillations there are not as pronounced as near the minima.}
    \label{F:EJ_time_evo}
\end{figure}

\subsubsection{Time Evolution and Ill-Posedness} \label{SSS:ill-posedness}

The numerical findings in the previous section are restricted to waves bifurcating from the zero constant state in the bidirectional Whitham model \eqref{E:EJ}.  This is largely
motivated by the fact that the authors have recently analytically studied the structure of the global bifurcation branch in \cite{EJ16}.  
In an effort to understand the \emph{nonlinear} dynamics about the numerically computed profiles in Figure \ref{F:EJ_global_bif_diagram}(a), we use a sixth-order pseudospectral operator splitting method \cite{Yoshida1990} to evolve \eqref{E:EJ_traveling_wave_model} with the 
numerical solutions as initial data, with the expectation that the initial profile translates as a traveling wave with corresponding wavespeed $c$, 
mapping onto itself after an integral number of temporal periods $T=\frac{2\pi}{\kappa c}$.  Interestingly, however, such time evolution attempts failed, resulting 
in wild oscillations after a short amount of time, e.g. 35.4\% of a temporal period for a $2\pi$-periodic wave of small 
waveheight $\max\psi - \min\psi \approx 0.00109$ and wavespeed $c \approx 0.872693$.  See Figure \ref{F:EJ_time_evo}.  This seems to suggest that either the solutions computed in Figure \ref{F:EJ_global_bif_diagram}(a)
are in fact \emph{not} traveling wave solutions of \eqref{E:EJ}, or that the local evolution in the PDE \eqref{E:EJ} about such waves is ill-posed.

In this direction, we note that it has recently been reported in \cite[Theorem 1]{EPW17}  that \eqref{E:EJ} is locally well-posed for all initial data $\eta(x,0)=\psi(x)$, $u(x,0)=\phi(x)$ with $\inf_{x \in [0,2\pi)} \psi(x) > 0$.  As the height profiles about the waves along the bifurcation curves constructed above all have a \emph{strictly negative infimum}, it follows
that the local well-posedness result in \cite{EPW17} does \emph{not apply} to the waves constructed in the previous section and hence, by extension,
also not to the waves constructed in \cite{EJ16}.  The time evolution results reported in Figure \ref{F:EJ_time_evo} strongly indicate that the positivity assumption
on the height profile $\eta$ in \cite{EPW17} is \emph{sharp} and cannot be removed in general.

In contrast, waves of \eqref{E:EJ} having strictly positive profile that are constructed by bifurcating from a positive constant state exhibit local well-posedness in time.  See Figure \ref{F:EJ_positive_time_evo_15per} and the surrounding discussion.  See Appendix \ref{AA:time_evolution_parameters} for further details about the time evolution.

\subsubsection{Analysis of Waves with Positive Height Profile}\label{SSS:positive}

Given the well-posedness result discussed above and the apparent failure of the waves considered in Section \ref{S:bf0} to have well-posed local dynamics,
it is natural to attempt to construct wavetrain solutions of \eqref{E:EJ} for which the local evolution is well-posed.  With the result of
\cite{EPW17} in mind, we seek to construct solutions of \eqref{E:EJ_traveling_wave_model} that have \emph{strictly positive height profiles} $\eta(x)$ by bifurcating from a 
\emph{positive} constant state. 
As discussed in \cite{EJ16}, the curves of trivial solutions for \eqref{E:EJ_traveling_wave_model} are given by
\[
c \mapsto 0, \quad c \mapsto \Gamma_\pm(c) := \frac{3c \pm \sqrt{8+c^2}}{2},
\]
and the waves studied in Section \ref{SS:analysis_EJ} above bifurcate from the zero-amplitude state, with the height profiles $\eta(x)$ all having strictly negative infima.
Since $\Gamma_-(c)>0$ for all $c>1$, this motivates an attempt to construct waves that bifurcate from the $\Gamma_-$ curve, which (at least for small waveheight)
are guaranteed to fall into the regime of well-posedness: see Figure \ref{F:EJ_positive_global_bif_diagram} below.

\begin{figure}[t]
    \centering
    
    \begin{subfigure}[t]{0.5\linewidth}
    \centering
	\includegraphics[scale=0.3]{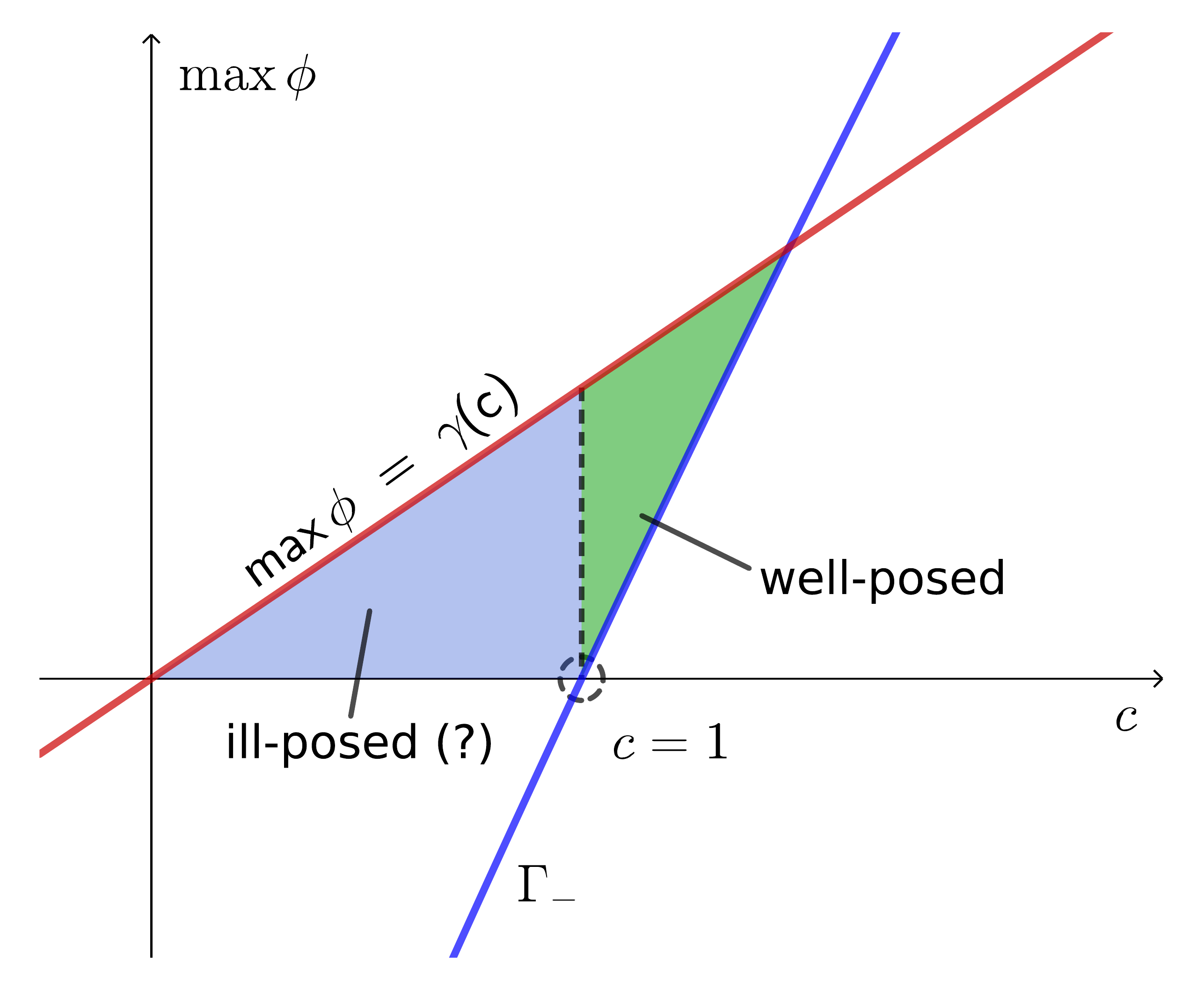}
	\caption{}
	\end{subfigure}%
	\begin{subfigure}[t]{0.5\linewidth}
	\centering
	\includegraphics[scale=0.45]{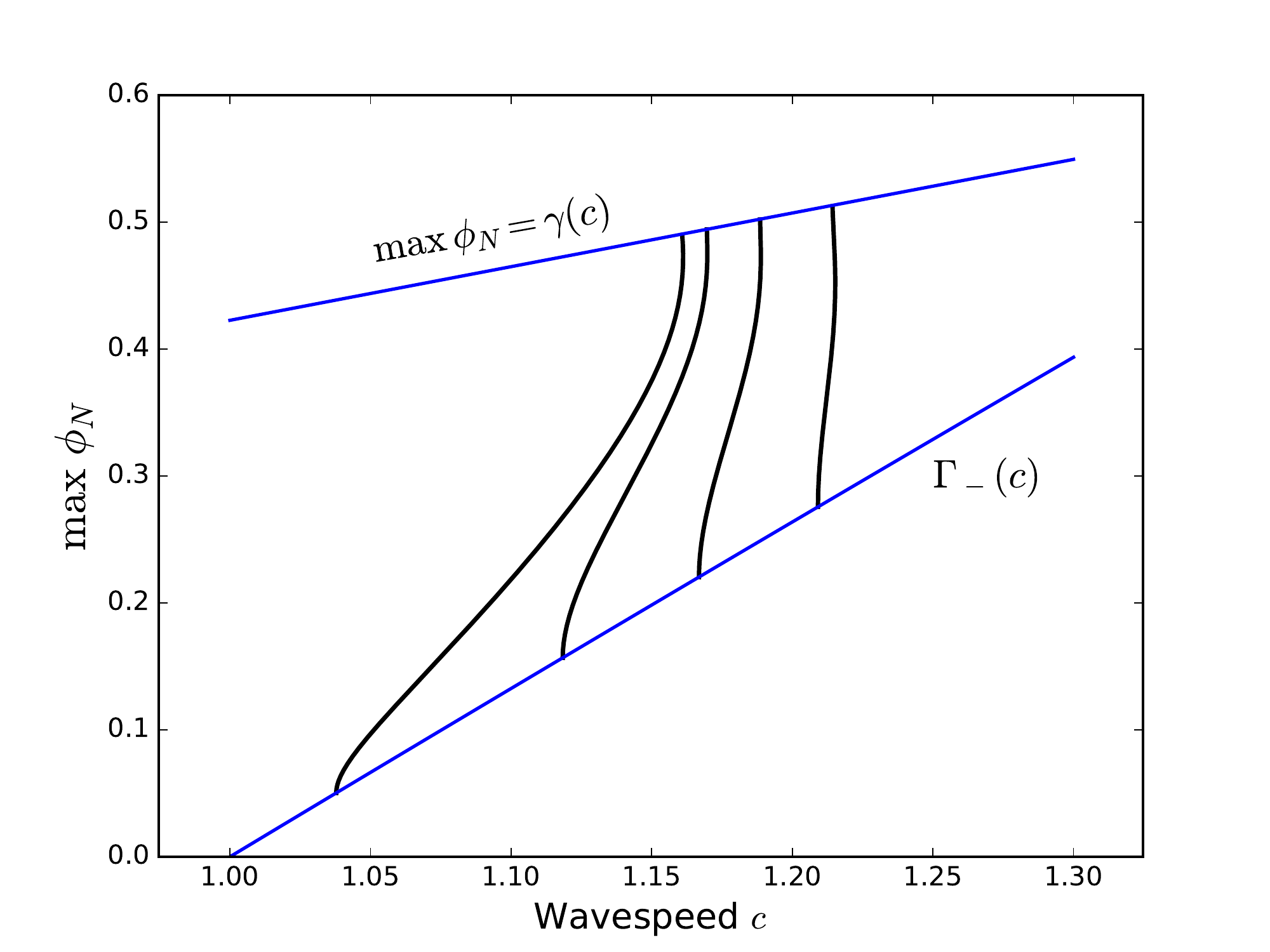}
	\caption{}
	\end{subfigure}
	
	\caption{(a) Bifurcation diagram of \eqref{E:EJ_profile}.  Waves bifurcating from the zero state are likely ill-posed, while waves bifurcating from a 
non-zero constant state (curve $\Gamma_-$) are well-posed.  (b) Numerical approximations of the bifurcation branches of even, one-sided monotone, $2\pi/\kappa$ periodic solutions of \eqref{E:EJ_profile}. From left to right, $\kappa=0.5, 1.0, 1.3, 1.6$.}
    \label{F:EJ_positive_global_bif_diagram}
\end{figure}

\begin{figure}[t]
    \centering
    
    \begin{subfigure}[t]{0.5\linewidth}
    \centering
	\includegraphics[scale=0.4]{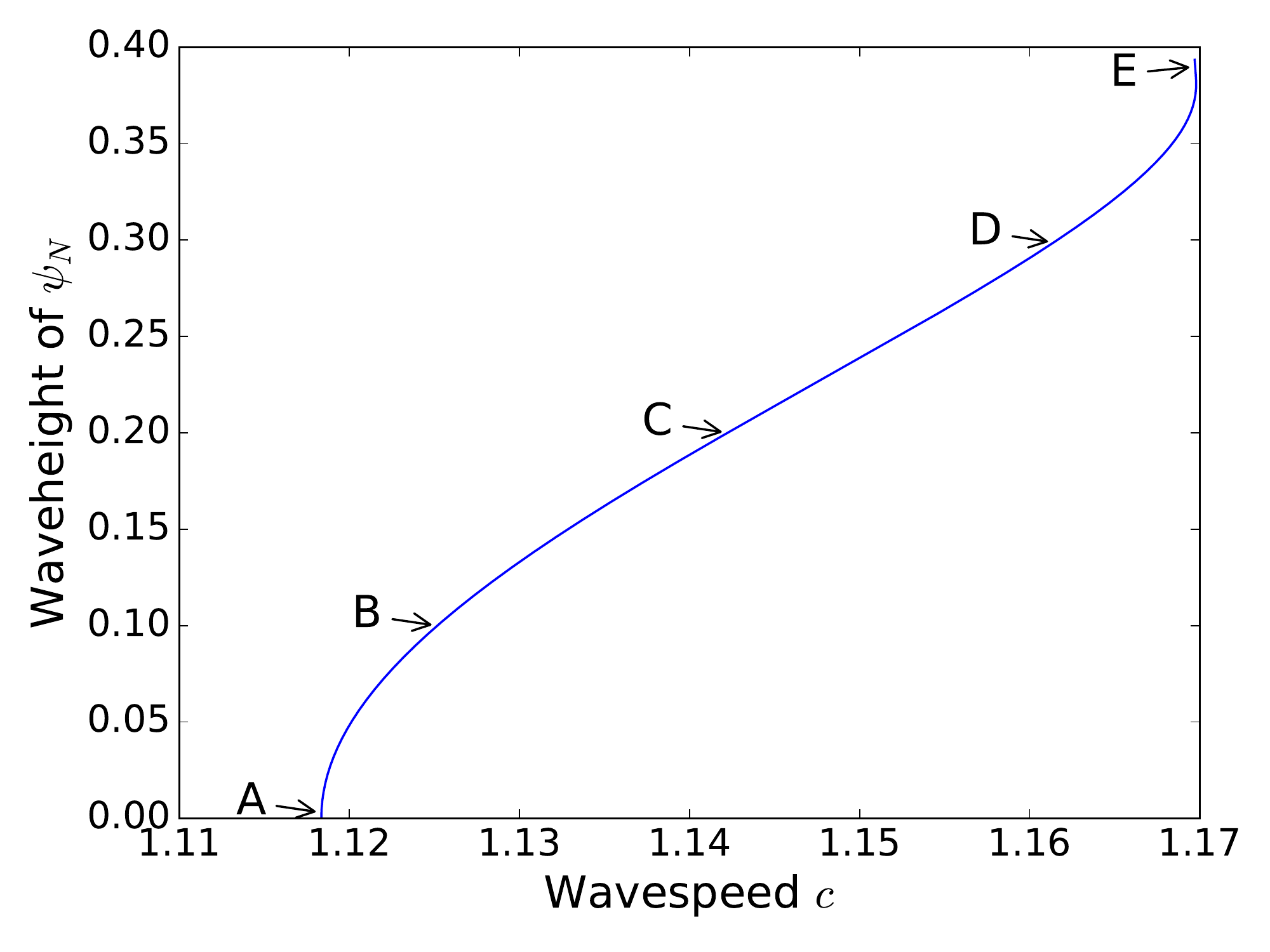}
	\caption{}
	\end{subfigure}%
	\begin{subfigure}[t]{0.5\linewidth}
	\centering
	\includegraphics[scale=0.4]{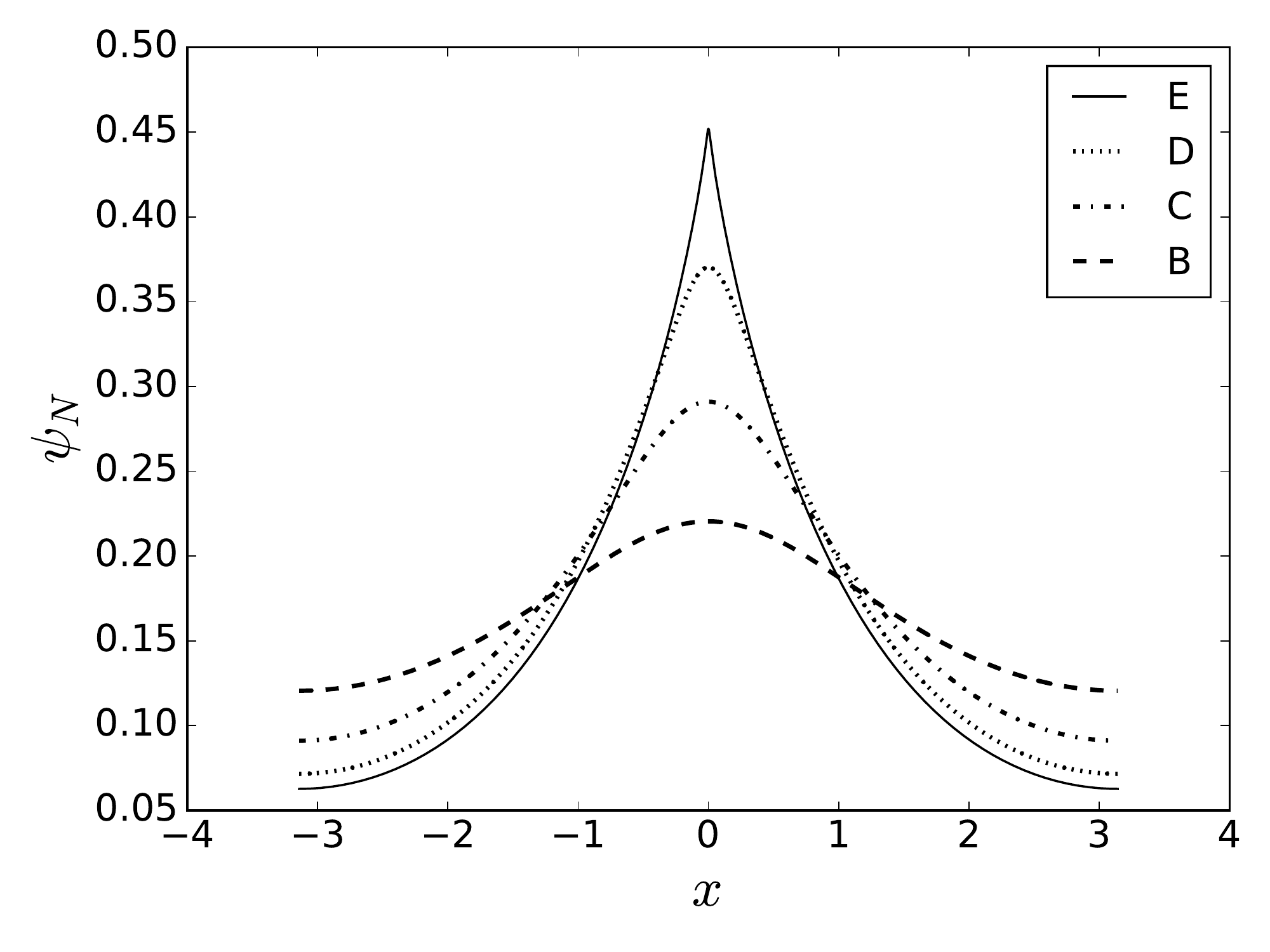}
	\caption{}
	\end{subfigure}
	
	\caption{Positive $2\pi$-periodic equilibrium profiles $\eta(x,t)=\psi(x)$ bifurcating from constant positive state.  In (a), we show the global bifurcation branch with points A-E labeled for forthcoming computations: Point A at $c \approx 1.1184$, height $\approx 0.003$; Point B at $c \approx 1.1252$, height $\approx 0.15$; Point C at $c \approx 0.8312$, height $\approx 0.30$; Point D at $c \approx 0.8138$, height $\approx 0.40$; point E at $c \approx 0.8051$, height $\approx 0.49$.  In (b), we show the height profiles $\psi$
associated to the waves A-E of increasing waveheight along the branch in (a).}
    \label{F:EJ_positive_bif_profiles}
\end{figure}

\begin{figure}[t]
    \centering
    \includegraphics[scale=0.45]{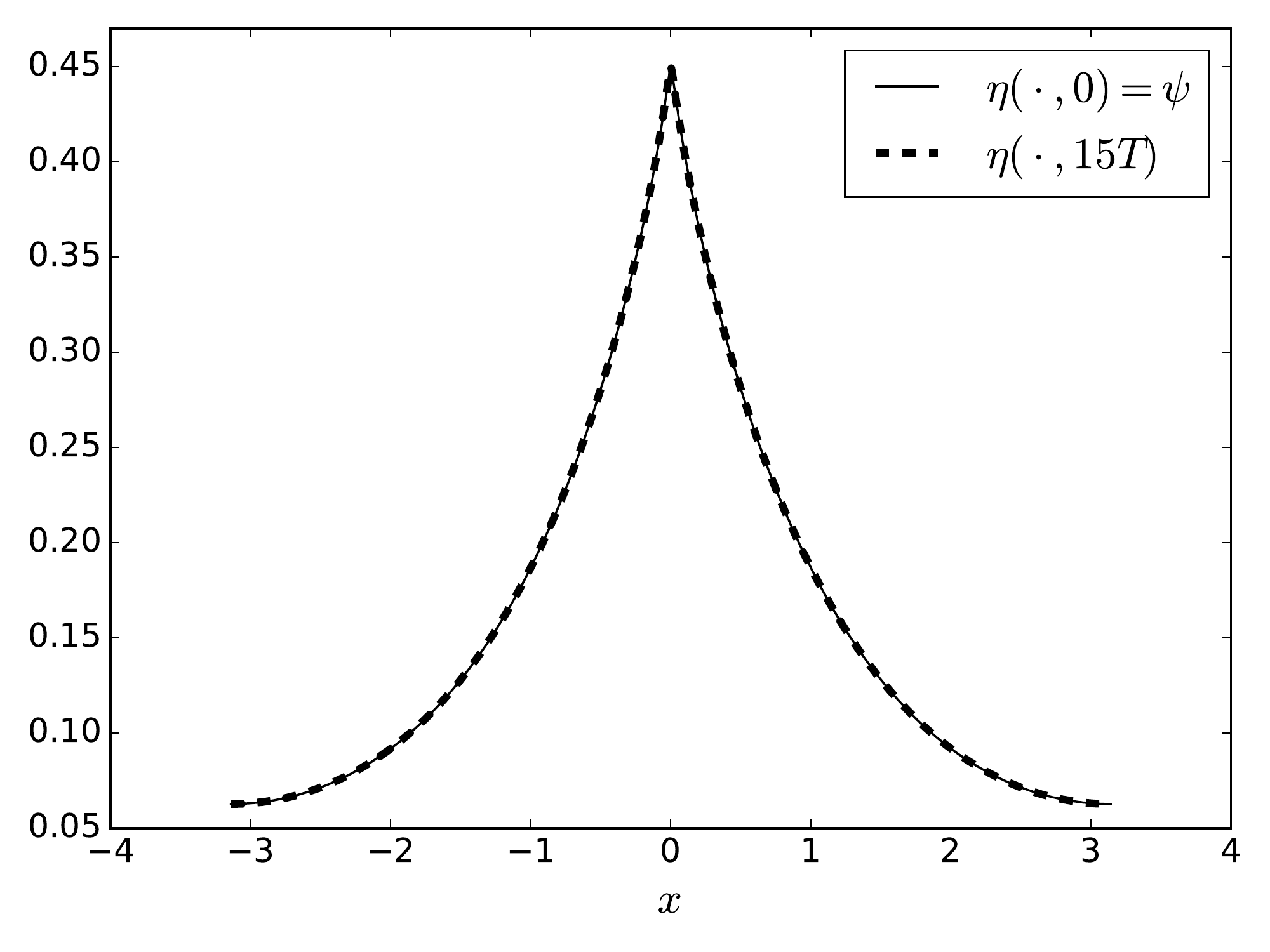}
    \caption{A time evolution of \eqref{E:EJ} using a positive equilibrium profile as initial data appears stationary after $15$ temporal periods, mapping closely onto itself.}
	\label{F:EJ_positive_time_evo_15per}    
\end{figure}

To this end, observe that non-zero constant profiles $\phi$ of \eqref{E:EJ_profile} satisfy
\begin{equation} \label{E:EJ_positive_bifurcation_profile}
Q(\phi;c) := c^2 - 1 - \frac{3}{2}c\phi + \frac{1}{2}\phi^2 = 0.
\end{equation}
Let $(\phi_*,c_*)$ be a solution of \eqref{E:EJ_positive_bifurcation_profile} for $c_* > 1$.  
In order for non-trivial, real-valued, even, periodic solutions to branch from the constant state at $c_*$, the kernel of the Frechet derivative $\delta_\phi Q(\phi_*;c_*)$ must be non-trivial.  If $v \in \ker(\delta_\phi Q(\phi_*;c_*))$ with $v(x) = \sum_{n=0}^\infty \widehat{v}(n)\cos(n \kappa x)$, then
\[
0 = \delta_\phi Q(\phi_*;c_*)v = \left( \CalK - c^2 + 3c\phi_* - \frac{3}{2}\phi_*^2 \right)v.
\]
Taking the Fourier transform yields
\[
0 = \left( \widehat{\CalK}(\kappa n) - c^2 + 3c\phi_* - \frac{3}{2}\phi_*^2 \right) \widehat{v}(n), \quad n \in \N_0,
\]
and we see that if
\begin{equation}  \label{E:EJ_positive_branching_condition}
\widehat{\CalK}(\kappa n_0) - c^2 + 3c\phi_* - \frac{3}{2}\phi_*^2 = 0
\end{equation}
for fixed $n_0 \in \N_0$, then
\[
\Span\{ \cos(\kappa n_0 x) \} = \ker(\delta_\phi Q(\phi_*,c_*)).
\]
Combining \label{E:EJ_positive_profile} and \eqref{E:EJ_positive_branching_condition} yields a system of necessary conditions for a non-trivial branch of $2\pi/(\kappa n_0)$-periodic
traveling wave solutions of \eqref{E:EJ} to appear:
\begin{equation} \label{E:EJ_positive_branching_system}
\left\{
\begin{array}{rcl}
\displaystyle c_*^2 - 1 - \frac{3}{2}c_*\phi_* + \frac{1}{2}\phi_*^2 &=& 0 \vspace{0.5em} \\
\displaystyle \widehat{\mathcal{K}}(\kappa n_0) - c_*^2 + 3c_*\phi_* - \frac{3}{2}\phi_*^2 &=& 0.
\end{array}
\right.
\end{equation}
For positive wavespeed, this system has a unique solution $(c_*(n_0), \phi_*(n_0))$ for each $n_0 \in \N_0$.  In the case of $n_0=1$, corresponding to $2\pi$-periodic
solutions of \eqref{E:EJ_profile}, we find
\[
c_* \approx 1.11834, \quad \phi_* \approx 0.15677.
\]
By standard Lyapunov-Schmidt arguments, it can be shown that a local bifurcation of $2\pi$-periodic wavetrains of asymptotically small waveheight occurs at $(\phi_*,c_*)$.
In particular, one finds the following local bifurcation formulas:
\begin{align*}
\phi(x;a) &:= \phi_* + a\cos(\kappa n_0 x) + O(a^2) \\
c(a) &:= c_* + \frac{3-4q}{24\phi_* - 16c_*} a^2 + O(a^3), \quad q := \frac{3}{4}\left[ \frac{\phi_* - c_*}{\widehat{\CalK}(2\kappa n_0) - c_*^2 + 3c_*\phi_* - \frac{3}{2}\phi_*^2} \right].
\end{align*}
Using the methods of Section \ref{S:bifurcation_methods} to approximate the wave profiles, we obtain a branch of solutions, depicted in Figure \ref{F:EJ_positive_bif_profiles} as waveheight vs. wavespeed, with all the height profiles $\psi$ possessing a strictly positive infimum and an increasingly sharp crest for large waveheights.   In particular, we note that the time evolution about these waves indeed appears to be \emph{well-posed}. In fact, a time evolution of a $2\pi$-periodic solution of \eqref{E:EJ} with wavespeed $c \approx 1.1698$ and $\max\psi-\min\psi \approx 0.387$ over 15 temporal periods $15T = 15 \cdot 2\pi/c$ closely mapped the initial data onto itself, with residual
\[
\|\eta(\cdot,15T) - \psi\|_{L^2} \approx 3.5 \times 10^{-6}.
\]
See Figure \ref{F:EJ_positive_time_evo_15per}.
This demonstrates that these computed profiles with positive minimum indeed form traveling wave solutions of \eqref{E:EJ} that exhibit locally well-posed time evolution.

\begin{figure}[t]
    \centering
    
    \begin{subfigure}[t]{0.5\linewidth}
    \centering
	\includegraphics[scale=0.3]{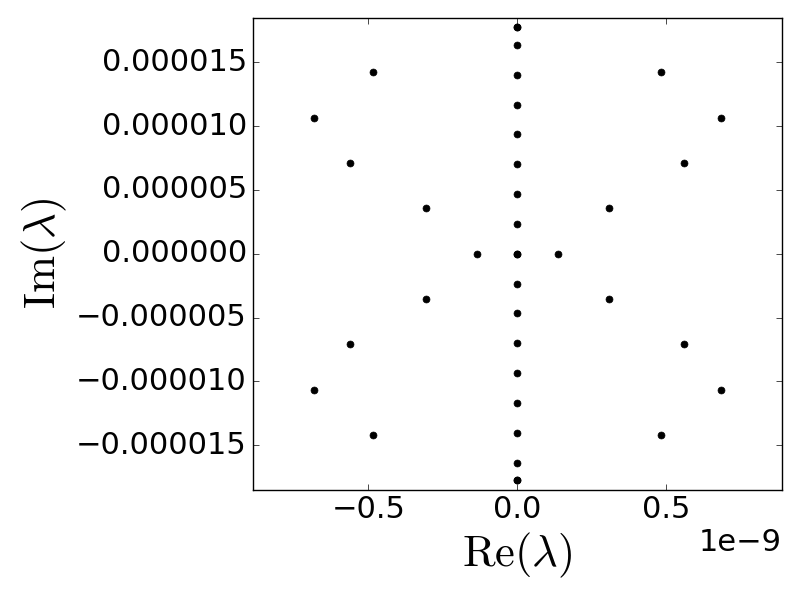}
	\caption{}
	\end{subfigure}%
	\begin{subfigure}[t]{0.5\linewidth}
	\centering
	\includegraphics[scale=0.3]{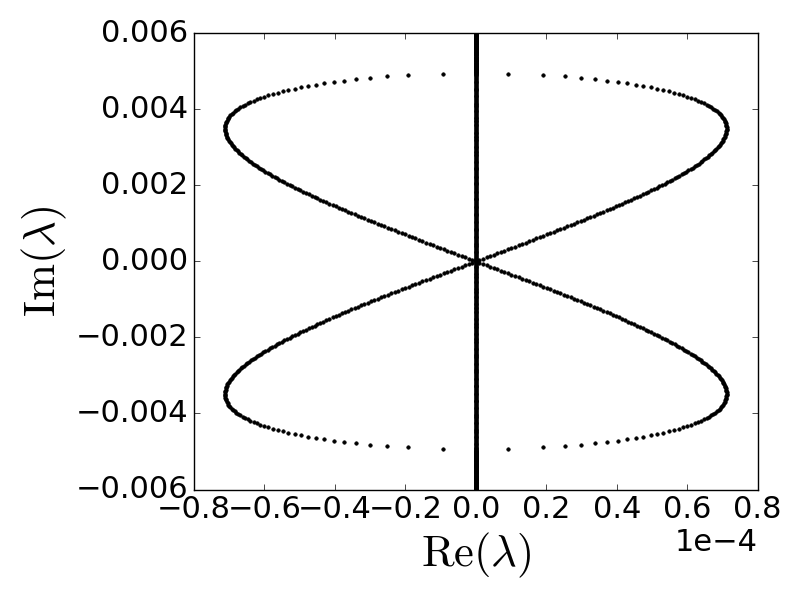}
	\caption{}
	\end{subfigure}
	
    \begin{subfigure}[t]{0.5\linewidth}
    \centering
	\includegraphics[scale=0.3]{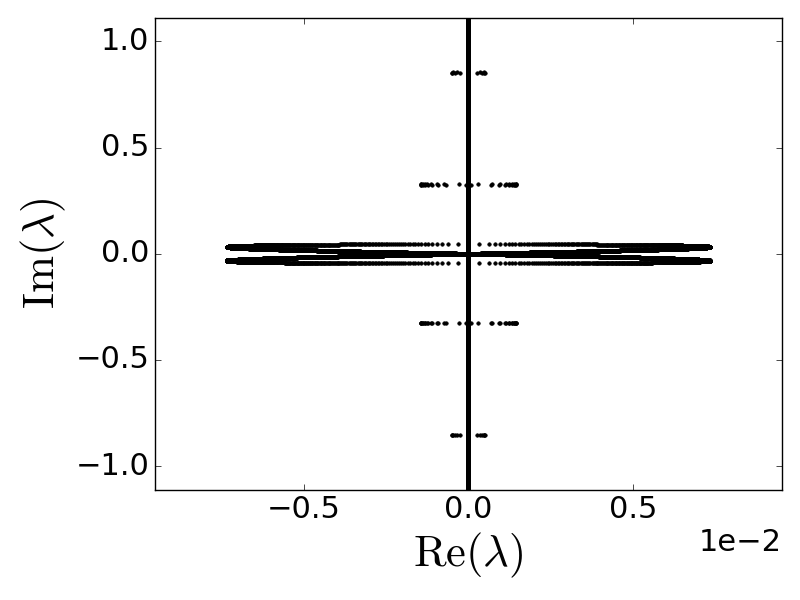}
	\caption{}
	\end{subfigure}%
	\begin{subfigure}[t]{0.5\linewidth}
	\centering
	\includegraphics[scale=0.3]{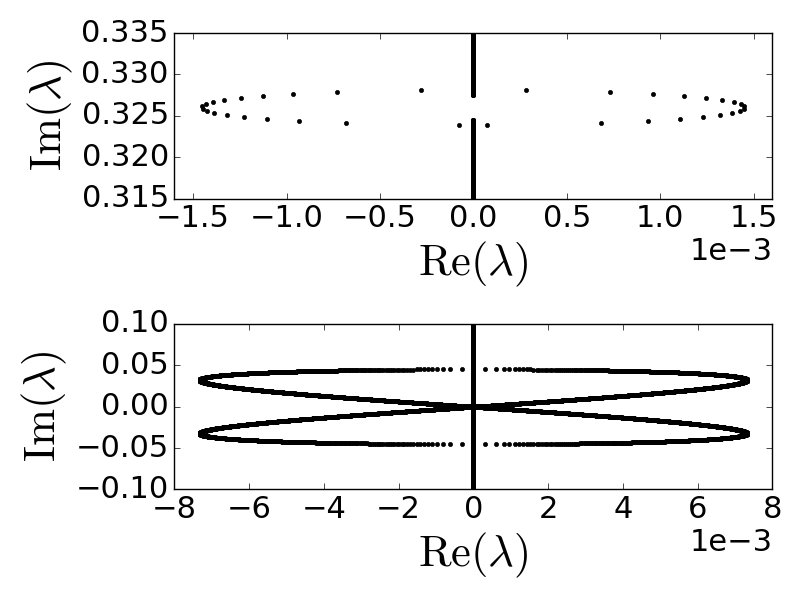}
	\caption{}
	\end{subfigure}	
	
    \caption{Spectral plots for $2\pi$-periodic positive waves for \eqref{E:EJ} at points  A (a), B (b), and D (c) along the bifurcation branch of $2\pi$-periodic solutions in Figure \ref{F:EJ_positive_bif_profiles}.  (d) Zoom-in on the high-frequency instability (top) and the modulational instability (bottom) in the spectral plot (c).}
    \label{F:EJ_positive_spectrum_plots}
\end{figure}

Concerning the stability of these $2\pi$-periodic positive height waves, we find that they behave similarly to ill-posed waves discussed in the previous section, with waves of sufficiently large
height demonstrating modulational and high-frequency instabilities: see Figure \ref{F:EJ_positive_spectrum_plots}.  In fact,  Figure \ref{F:EJ_positive_spectrum_plots}(a) suggests
that  $2\pi$-periodic waves of small amplitude exhibit modulational instability.

\subsection{Analysis of System \eqref{E:HPShallow_model}}  \label{SS:analysis_HP}

We now turn our attention to the bidirectional Whitham model \eqref{E:HPShallow_model}:
\[
\left\{
\begin{array}{rcl}
u_t  &=& -\CalK \eta_x - uu_x \\
\eta_t &=& -u_x - (\eta u)_x.
\end{array}
\right.
\]
In particular, we are interested in performing an analogous study for \eqref{E:HPShallow_model} that was performed for the model \eqref{E:EJ} in Section \ref{SS:analysis_EJ} above.  
Below, we compute global bifurcation diagrams, including large amplitude solutions, and analyze the spectrum of the linearization of \eqref{E:HPShallow_model} about these solutions.  

\begin{figure}[t]
    \centering
    
    \begin{subfigure}[t]{0.5\linewidth}
    \centering
	\includegraphics[scale=0.35]{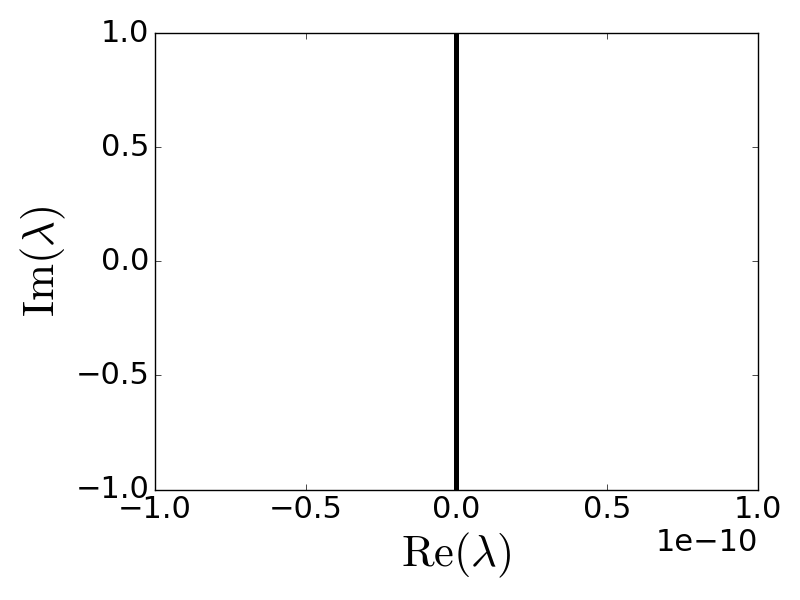}
	\caption{}
	\end{subfigure}%
	\begin{subfigure}[t]{0.5\linewidth}
	\centering
	\includegraphics[scale=0.35]{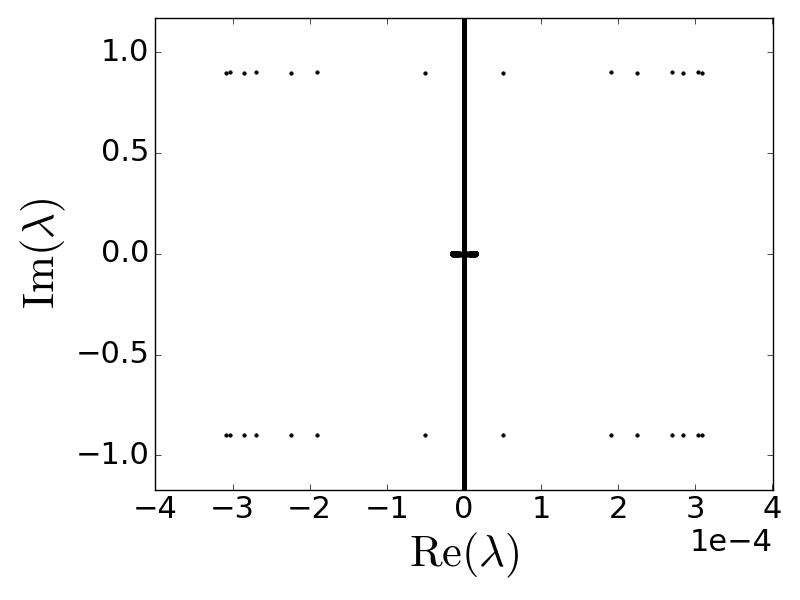}
	\caption{}
	\end{subfigure}
	
    \caption{Spectral plots for small $2\pi/\kappa$-periodic solutions of \eqref{E:HPShallow_profile} with sub-critical $\kappa=1.609 < \kappa_c$).  
   In (a), the waveheight of the underlying wave $\psi$ is $0.03$, while in (b) it is $0.15$.}
    \label{F:HPShallow_subcrit_spectrum_plots}
\end{figure}

As mentioned in the introduction, this model was proposed and analyzed in \cite{HP_shallow}, which examines the local bifurcation and spectral stability
of asymptotically small amplitude periodic traveling waves, finding that this model exhibits the Benjamin-Feir instability.  Precisely, they prove the existence of a 
critical wave number $\kappa_c \approx 1.610$ such that asymptotically small-amplitude wavetrains of period $2\pi/\kappa$ are modulationally when $\kappa > \kappa_c$ (``super-critical"), 
while they are spectrally stable for $0 < \kappa < \kappa_c$ (``subcritical").  Below, we numerically confirm this result for small amplitude waves, and also demonstrate that large amplitude waves are spectrally unstable in both the sub-critical and super-critical regimes.  See Figure \ref{F:HPShallow_supcrit_spectrum_plots} for spectrum plots in the super-critical case and Figure \ref{F:HPShallow_subcrit_spectrum_plots} for spectrum plots the sub-critical case.  Furthermore, from our experiments we make a conjecture regarding the nonexistence of a singular wave of greatest height.

\begin{figure}[t]
    \centering
    
    \begin{subfigure}[t]{0.5\linewidth}
    \centering
	\includegraphics[scale=0.3]{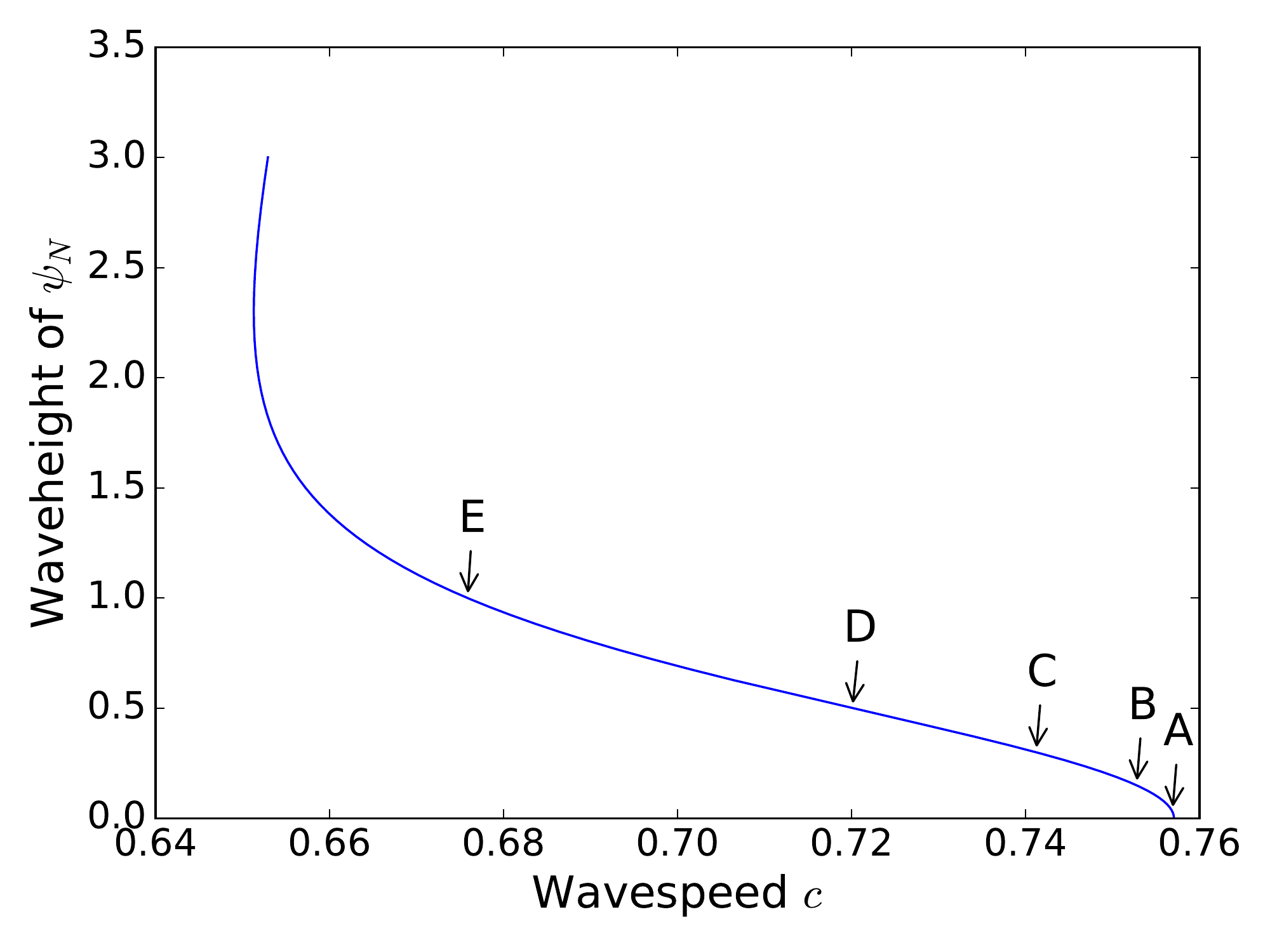}
	\caption{}
	\end{subfigure}%
	\begin{subfigure}[t]{0.5\linewidth}
	\centering
	\includegraphics[scale=0.3]{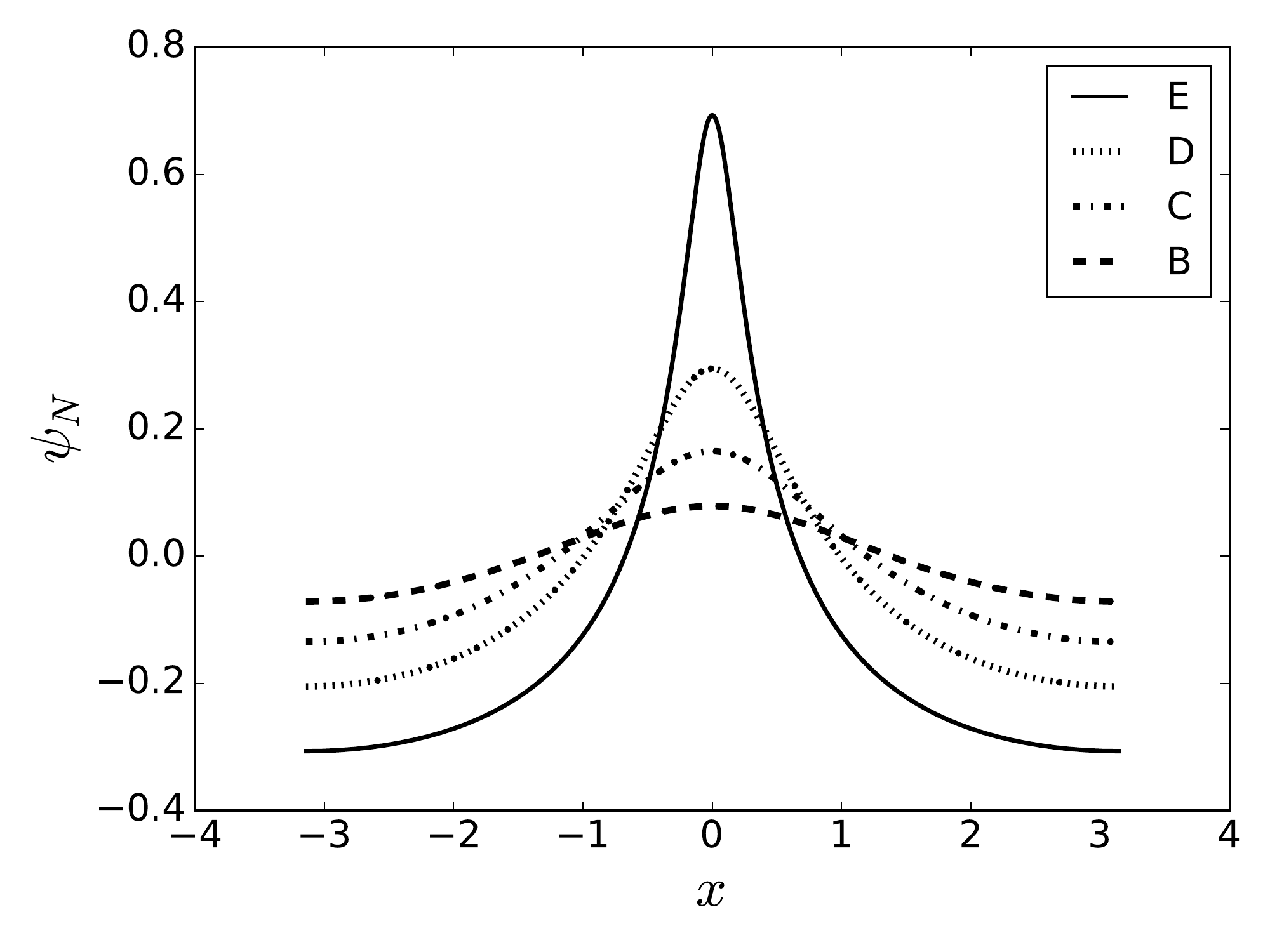}
	\caption{}
	\end{subfigure}
	
    \begin{subfigure}[t]{0.5\linewidth}
    \centering
	\includegraphics[scale=0.3]{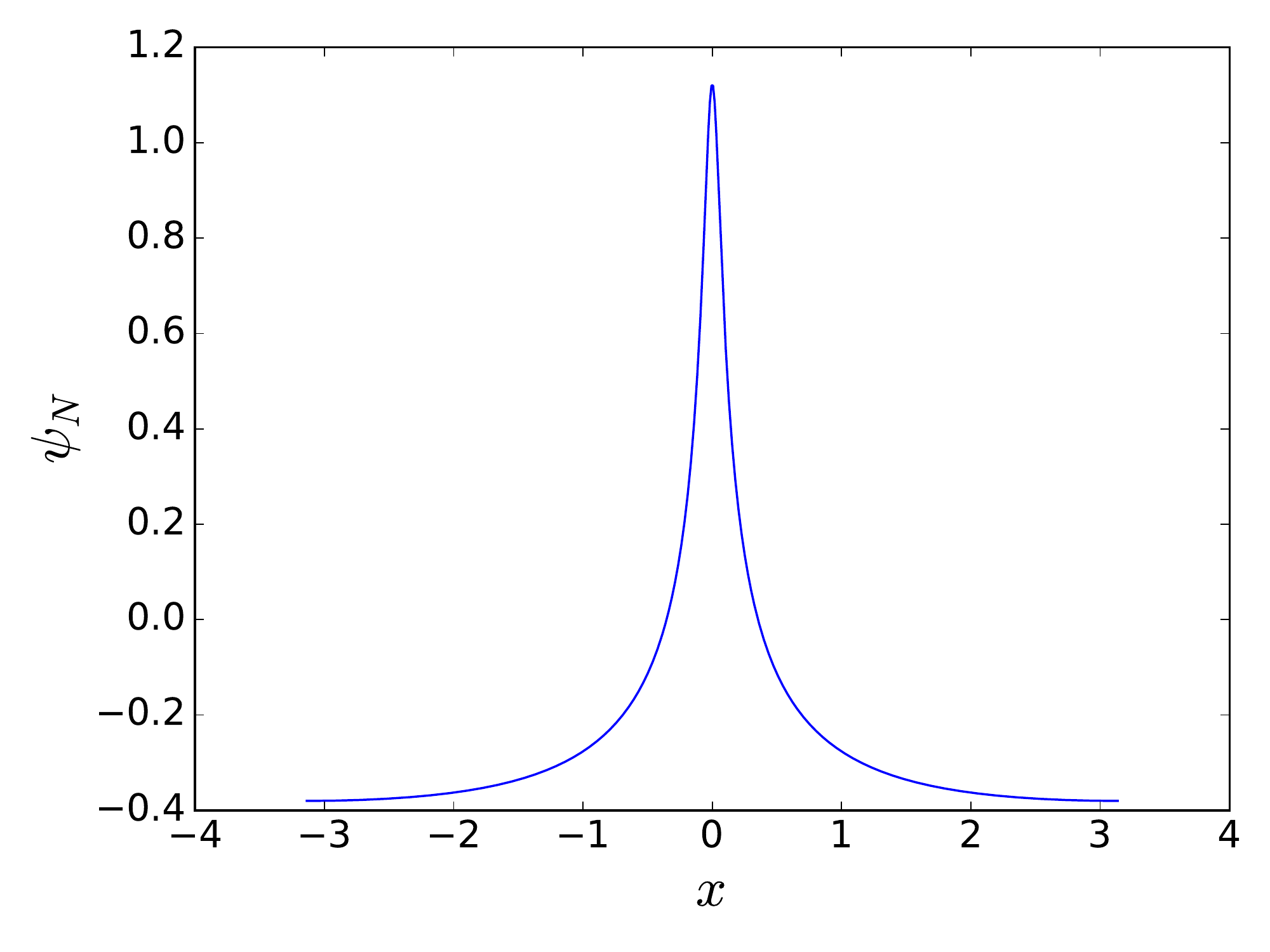}
	\caption{}
	\end{subfigure}%
	\begin{subfigure}[t]{0.5\linewidth}
	\centering
	\includegraphics[scale=0.3]{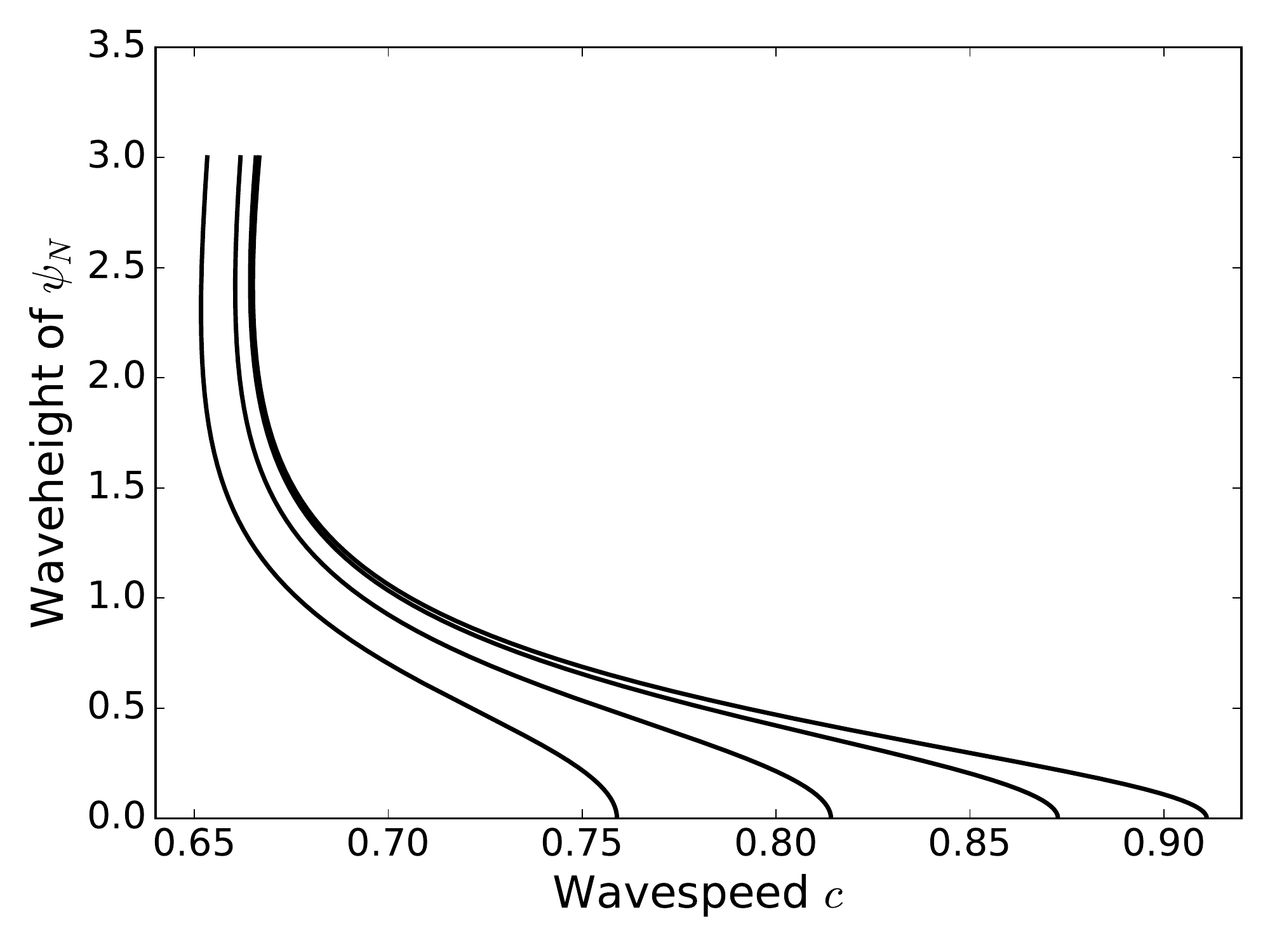}
	\caption{}
	\end{subfigure}	
	
	\caption{Bifurcation branches of \eqref{E:HPShallow_profile} for super-critical $\kappa=1.611 > \kappa_c$.  (a) Locations on the bifurcation branch of periodic solutions corresponding to super-critical $\kappa=1.611$ of \eqref{E:HPShallow_profile} that are sampled for profile and spectral plots: Point A at $c \approx 0.7569$, height $\approx 0.03$; Point B at $c \approx 7528$, height $\approx 0.15$; Point C at $c \approx 0.7412$, height $\approx 0.30$; Point D at $c \approx 0.7201$, height $\approx 0.50$; Point E at $c \approx 0.6759$, height $\approx 1.00$.  (b) Profiles corresponding to the points along the bifurcation branch labeled in (a).
(c) A super-critical wave ($\kappa=1.611$) of large waveheight; $c \approx 0.6745$, $\text{waveheight} \approx 1.500$.  The crest of the wave is still smooth despite the large waveheight.
(d)   Numerical approximations of the $2\pi/\kappa$ periodic solutions of \eqref{E:HPShallow_profile} with wavespeed $c$.  From left to right, $\kappa=1.6, 1.3, 1, 0.8$.}
    \label{F:HPShallow_supcrit_bifurcation_and_profiles}
\end{figure}

\begin{figure}[t]
    \centering
    
    \begin{subfigure}[t]{0.33\linewidth}
	\includegraphics[scale=0.26]{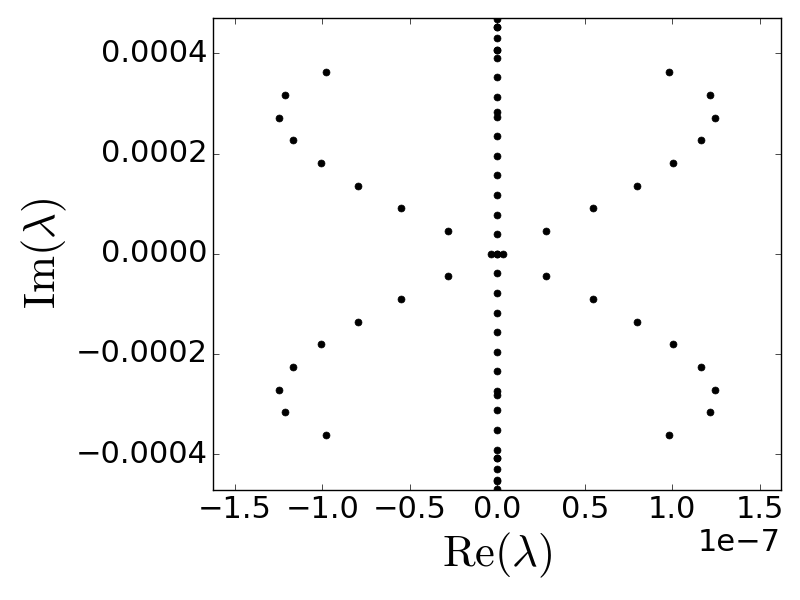}
	\caption{}
	\end{subfigure}%
    \begin{subfigure}[t]{0.33\linewidth}
	\includegraphics[scale=0.26]{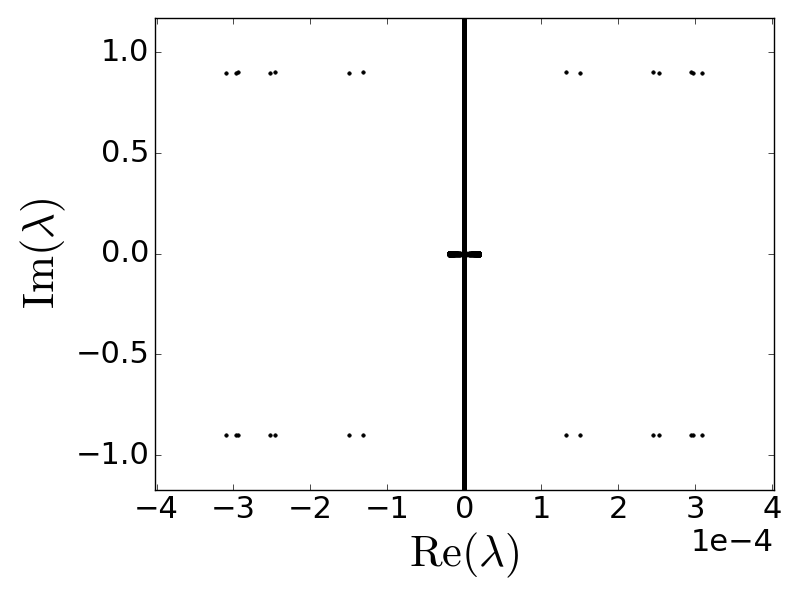}
	\caption{}
	\end{subfigure}%
	\begin{subfigure}[t]{0.33\linewidth}
	\includegraphics[scale=0.26]{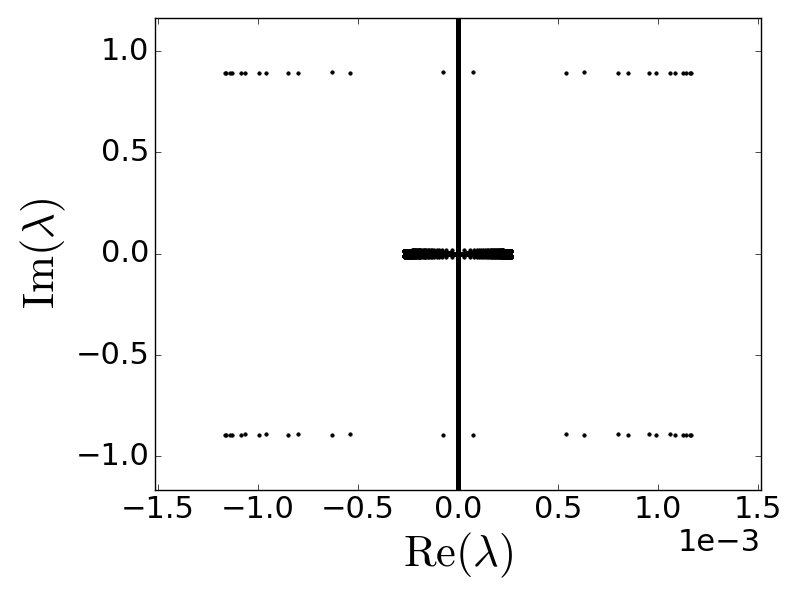}
	\caption{}
	\end{subfigure}%
	
    \begin{subfigure}[t]{0.33\linewidth}
	\includegraphics[scale=0.26]{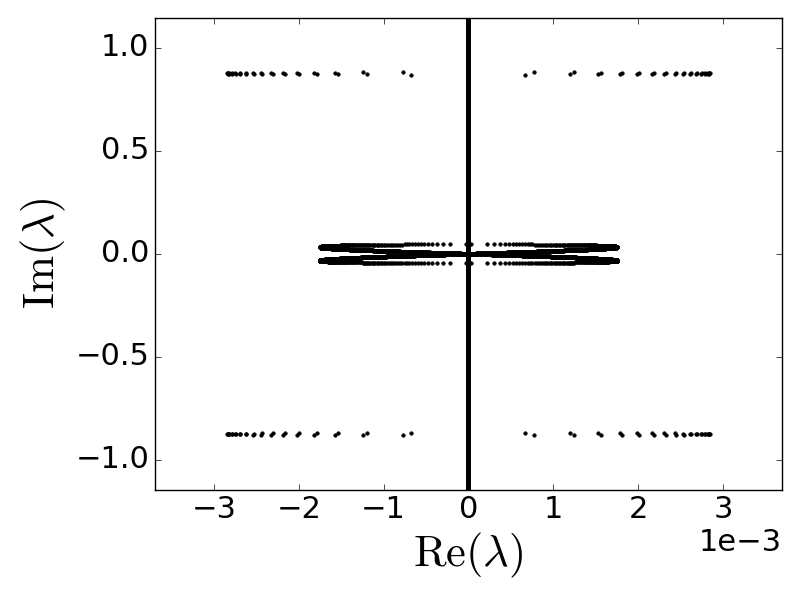}
	\caption{}
	\end{subfigure}%
    \begin{subfigure}[t]{0.33\linewidth}
	\includegraphics[scale=0.26]{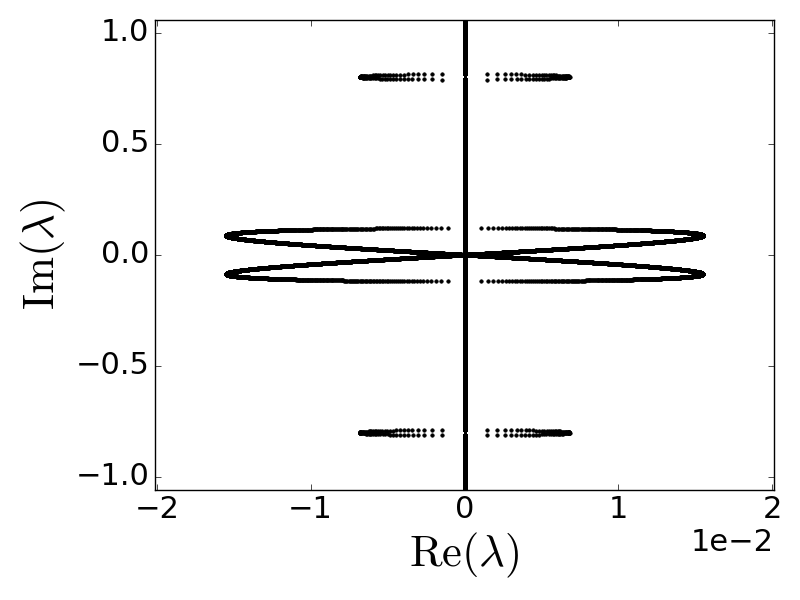}
	\caption{}
	\end{subfigure}%
    \begin{subfigure}[t]{0.33\linewidth}
	\includegraphics[scale=0.26]{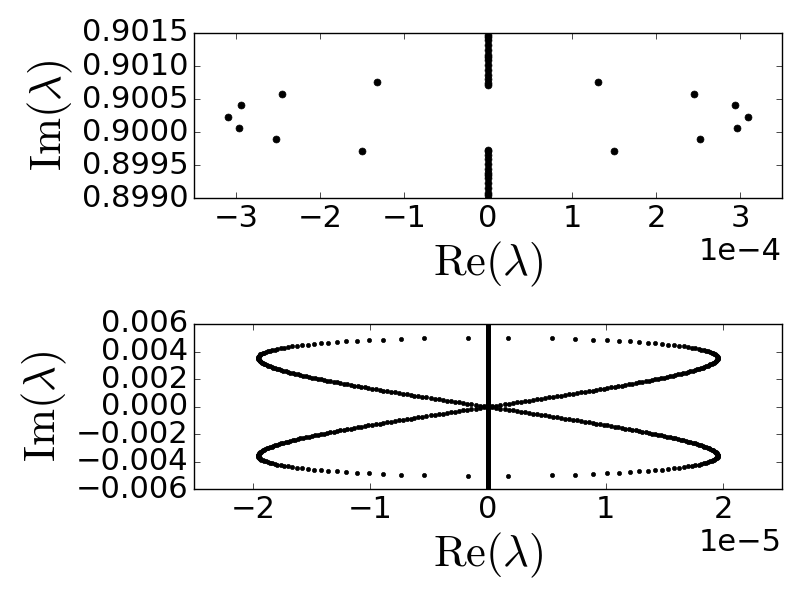}
	\caption{}
	\end{subfigure}

	\caption{(a)-(e) Spectral plots for $2\pi/\kappa$-periodic waves A-E selected along the bifurcation branch in Figure \ref{F:HPShallow_supcrit_bifurcation_and_profiles}(a) with super-critical 
$\kappa=1.611$.
(f) Zoom-in on the high-frequency instability (top) and modulational instability (bottom) for the spectrum in (b). }
    \label{F:HPShallow_supcrit_spectrum_plots}
\end{figure}

To generate profiles, we follow the program of the numerical methods described in Section \ref{S:bifurcation_methods}.  First, we transform \eqref{E:HPShallow_model} to traveling wave coordinates:
\[
\left\{
\begin{array}{rcl}
u_t &=& cu_x - \CalK \eta_x - uu_x =: F(u,u_x,\eta,\eta_x;c) \\
\eta_t &=& c\eta_x - u_x - (\eta u)_x =: G(u,u_x,\eta,\eta_x;c).
\end{array}
\right.
\]
An equilibrium solution $(u,\eta) (x,t)= (\phi,\psi)(x)$ of the above traveling wave system satisfies
\[
\left\{
\begin{array}{rcl}
-c\phi' + \CalK \psi' + \phi\phi' &=& 0 \\
-c\psi' + \phi' + (\psi\phi)' &=& 0
\end{array}
\right.
\]
which, upon integrating and setting integration constants to zero, yields 
\[
\left\{
\begin{array}{rcl}
-c\phi + \CalK \psi + \frac{1}{2}\phi^2 &=& 0 \\
-c\psi + \phi + \psi\phi &=& 0.
\end{array}
\right.
\]
Then resolving $\phi$ in terms of $\psi$ yields the scalar profile equation
\begin{equation} \label{E:HPShallow_profile}
\CalK \psi = \frac{c^2 \psi}{1+\psi} - \frac{c^2 \psi^2}{2(1+\psi)^2} =: g(c,\psi).
\end{equation}
Using the methods in Section \ref{SS:cosine_collocation_method}, we numerically compute the bifurcation branch and periodic profiles with various waveheights with super-critical $\kappa=1.611$: see  Figure \ref{F:HPShallow_supcrit_bifurcation_and_profiles}(a)-(c).  Further, in Figure \ref{F:HPShallow_supcrit_bifurcation_and_profiles}(d) we compute  bifurcation branches of $2\pi/\kappa$-periodic profiles for various $\kappa$.

In contrast to the model \eqref{E:EJ} analyzed in Section \ref{SS:analysis_EJ} above, we conjecture that the bifurcation branches associated with \eqref{E:HPShallow_model} will \emph{not possess peaked/cusped waves of maximum height}.  In fact, we believe the profiles
along the global bifurcation branches will be smooth for arbitrarily large waveheight.  To understand this, observe that differentiating the profile equation \eqref{E:HPShallow_profile} and rearranging yields
\[
\psi' = \frac{c^2}{(1+\psi)^3} \CalK\psi'.
\]
Since $\CalK \psi'$ has the same regularity as $\psi$, we have that $\psi'$ is smooth so long as the wave speed $c$ remains bounded away from zero and the
profiles $\psi$ remain bounded away from $\psi=-1$ along the bifurcation branch.  From  Figure \ref{F:HPShallow_supcrit_bifurcation_and_profiles}, it seems plausible
that \emph{both} of these conditions hold uniformly along the bifurcation branch, indicating that the bifurcation branch does not terminate in a highest wave.
See Figure \ref{F:HPShallow_supcrit_bifurcation_and_profiles}(d) for a periodic profile of large waveheight having a visibly smooth crest.
We leave the analytical verification of this conjecture as an interesting open problem.

\begin{remark}
Since this model is conjectured to not possess peaked/cusped waves of maximum height, the waveheight vs. wavespeed plots in Figures \ref{F:HPShallow_supcrit_bifurcation_and_profiles}(a), (d) were stopped at height $3$. As seen in Figure \ref{F:HPShallow_supcrit_bifurcation_and_profiles}(a), a turning point occurs near waveheight $3$, after which the wavespeed increases as the waveheight increases.  Numerically, a second turning point is also observed for waves of larger height, but we believe that this second turning point is only due to truncation, at its location varies considerably as larger numbers of Fourier modes $N$ are used.  See Figure \ref{F:HPShallow_secondary_turning_point}.

\begin{figure}[t]
    \centering
    \includegraphics[scale=0.5]{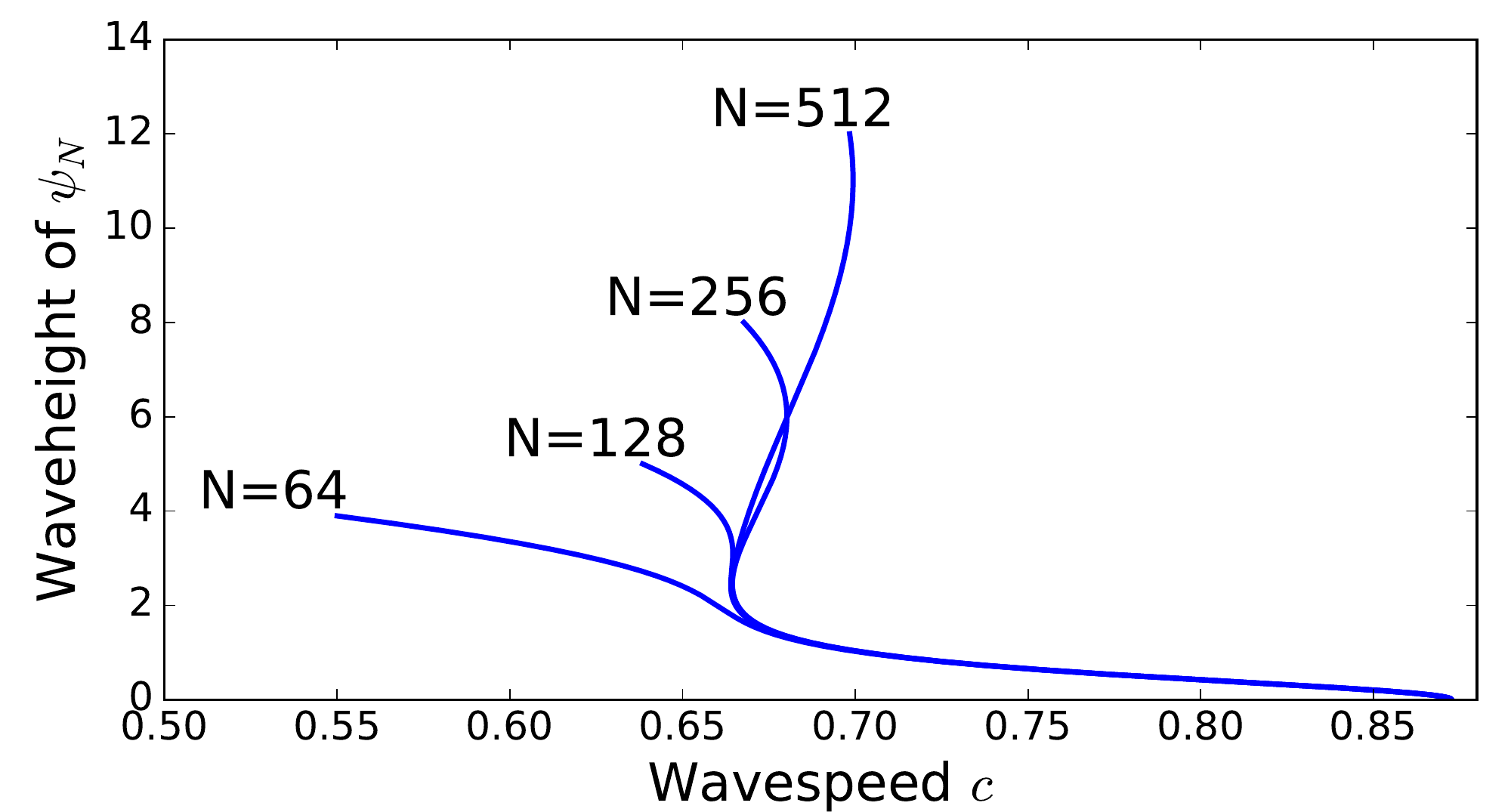}
    \caption{Secondary turning points in the bifurcation diagrams for system \eqref{E:HPShallow_model} are present at large waveheight, but this is believed to be purely due to truncation since its location increases as the number of Fourier modes $N \in \{64,128,256,512\}$ increases.}
	\label{F:HPShallow_secondary_turning_point}    
\end{figure}
\end{remark}

\begin{figure}[t]
    \centering
    
    \begin{subfigure}[t]{0.33\linewidth}
    \centering
	\includegraphics[scale=0.28]{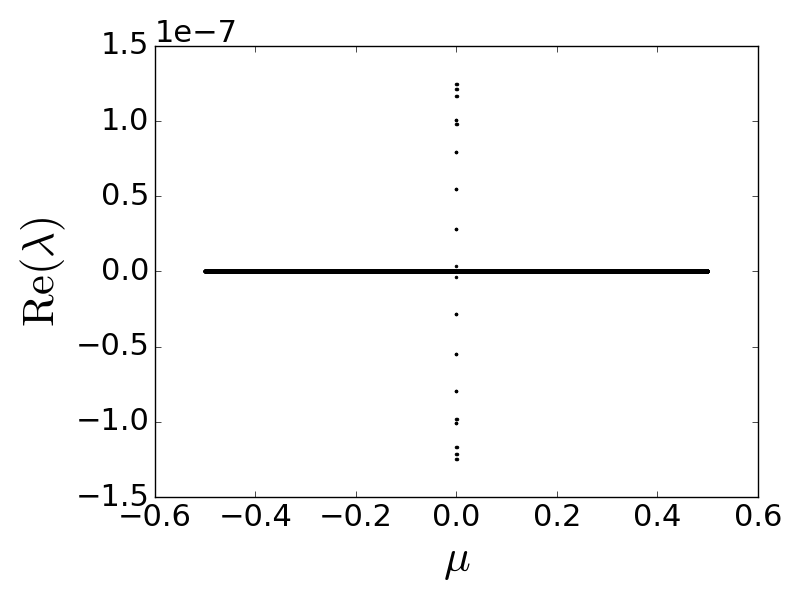}
	\caption{}
	\end{subfigure}%
    \begin{subfigure}[t]{0.33\linewidth}
    \centering
	\includegraphics[scale=0.28]{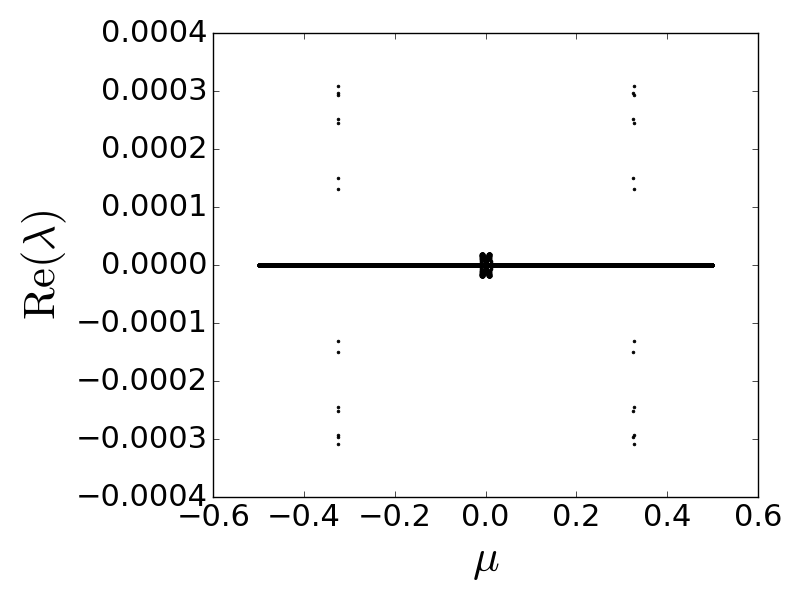}
	\caption{}
	\end{subfigure}%
	\begin{subfigure}[t]{0.33\linewidth}
	\centering
	\includegraphics[scale=0.28]{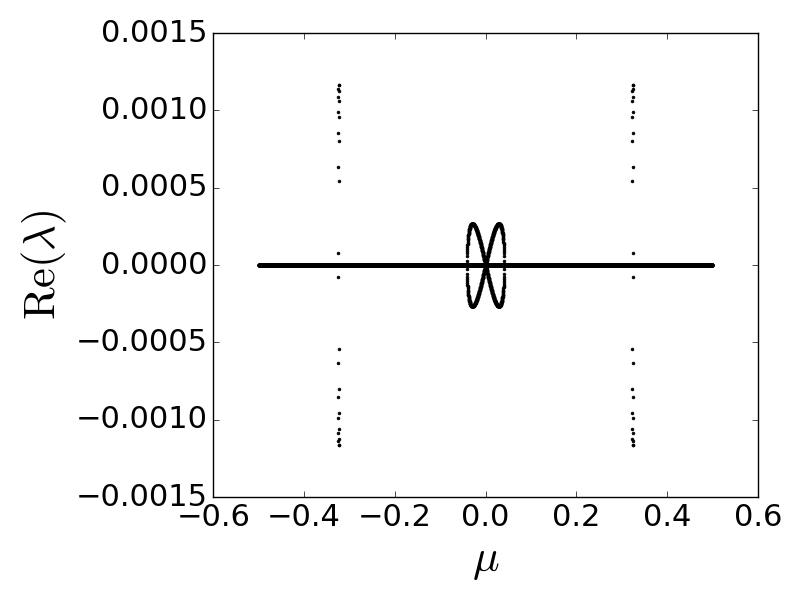}
	\caption{}
	\end{subfigure}%
	
    \begin{subfigure}[t]{0.33\linewidth}
    \centering
	\includegraphics[scale=0.28]{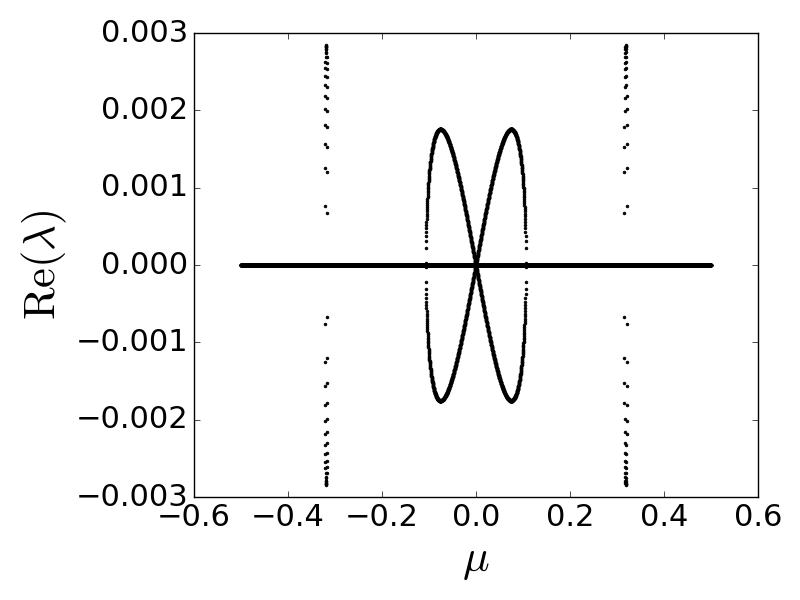}
	\caption{}
	\end{subfigure}%
	\hspace{1em}
    \begin{subfigure}[t]{0.33\linewidth}
    \centering
	\includegraphics[scale=0.28]{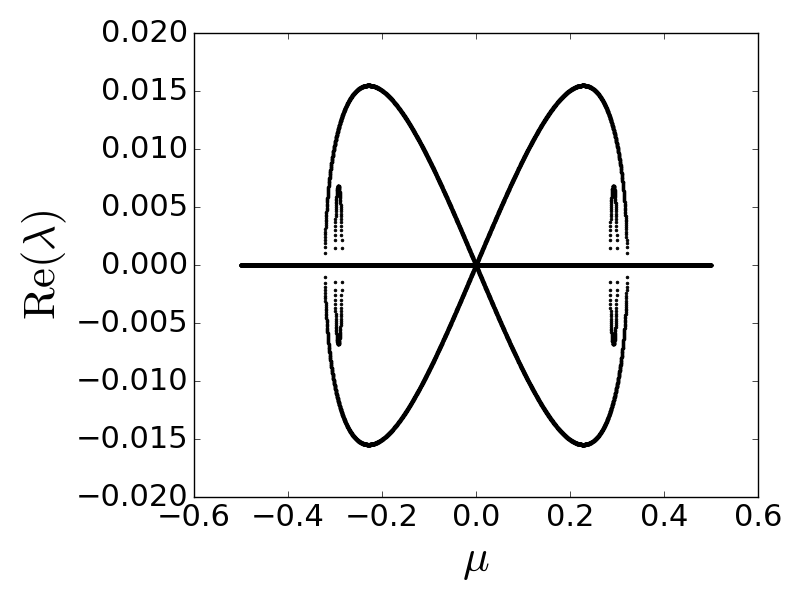}
	\caption{}
	\end{subfigure}

     \caption{(a)-(e) $\Re(\lambda)$ vs. $\mu$ for each of the spectral plots (a)-(e), respectively, for the super-critical $\kappa$ in Figure \ref{F:HPShallow_supcrit_spectrum_plots}.}
    \label{F:HPShallow_supcrit_Res}
\end{figure}

To study the dynamical stability of the periodic wavetrains computed in Figure \ref{F:HPShallow_supcrit_bifurcation_and_profiles}, we use the general methods of Appendix \ref{A:general_FFHM}.
In particular, we build  truncated bi-infinite matrices whose eigenvalues approximate the spectrum of the linearization:
\begin{align*}
\widehat{A}^\mu(m,l) &= ic\kappa(\mu+l)\delta_{m,l} - i\kappa(\mu+m)\widehat{\phi}(m-l) \\
\widehat{B}^\mu(m,l) &= -i\kappa(\mu+l)\widehat{\CalK}(\kappa(\mu+l)) \delta_{m,l} \\
\widehat{C}^\mu(m,l) &= -i\kappa(\mu+l)\delta_{m,l} - i\kappa(\mu+m)\widehat{\psi}(m-l) \\
\widehat{D}^\mu(m,l) &= ic\kappa(\mu+l)\delta_{m,l} - i\kappa(\mu+m)\widehat{\phi}(m-l).
\end{align*}
Plots of the spectrum at the sampled points along the bifurcation branch in Figure \ref{F:HPShallow_supcrit_bifurcation_and_profiles}(a) are shown in Figure \ref{F:HPShallow_supcrit_spectrum_plots}.  In particular, Figure \ref{F:HPShallow_supcrit_spectrum_plots}(a) demonstrates a modulational instability for a small waveheight, while the other plots for larger waveheights show both a modulational instability and a high-frequency instability.  Moreover, plots of the growth rate $\Re(\lambda)$ vs. $\mu$ for the these spectra are shown in Figure \ref{F:HPShallow_supcrit_Res}.  
Moreover, the spectral stability of small-amplitude sub-critical solutions is shown in Figure \ref{F:HPShallow_subcrit_spectrum_plots}(a), as predicted by analytical theory in \cite{HP_shallow}.  However, even for $\kappa$ in the sub-critical regime, waves of sufficiently large waveheight develop both modulational and high-frequency instabilities (see Figure \ref{F:HPShallow_subcrit_spectrum_plots}(b)) while small amplitude waves are spectrally stable.  Moreover, waves are stable with respect to co-periodic perturbations.

\subsection{Boussinesq-Whitham}  \label{SS:BW}

Finally, we now turn our attention to the scalar \emph{Boussinesq-Whitham} model \eqref{E:BoussinesqWhitham_equation} proposed in \cite[Section 5]{DT1}:
\[
u_{tt} = \partial_x^2 (u^2 + \CalK u).
\]
Since \cite{DT1} does not provide details for the local bifurcation theory of this model, we sketch the Lyapunov-Schmidt reduction here.  Writing the second order equation \eqref{E:BoussinesqWhitham_equation} as a first order system, we have
\[
\left\{
\begin{array}{rcl}
u_t &=& \eta \\
\eta_t &=& \partial_x^2 (u^2 + \CalK u),
\end{array}
\right.
\]
and changing to traveling wave coordinates yields
\[
\left\{
\begin{array}{rcl}
u_t &=& \eta + cu_x =: F(u,u_x,\eta,\eta_x;c) \\
\eta_t &=& c\eta_x + \partial_x^2 (u^2 + \CalK u) := G(u,u_x,\eta,\eta_x;c).
\end{array}
\right.
\]
Equilibrium solutions $(u,\eta)=(\phi,\psi)$ of this system satisfy
\begin{equation} \label{E:BoussinesqWhitham_equilibrium_system}
\left\{
\begin{array}{rcl}
-c\phi' &=& \psi \\
-c\psi' &=& (\phi^2 + \CalK \phi)'',
\end{array}
\right.
\end{equation}
in which we can resolve $\psi$ in terms of $\phi$ by integrating the second equation, setting the resulting constant of integration to zero, yielding the relation
\[
-c\psi = (\phi^2 + \CalK \phi)'.
\]
Using this relationship in the first equation of \eqref{E:BoussinesqWhitham_equilibrium_system} and integrating, again setting the constant of integration to zero, yields the scalar profile
equation
\begin{equation}  \label{E:BoussinesqWhitham_profile}
\CalK \phi = c^2 \phi - \phi^2 =: g(c,\phi).
\end{equation}
for the velocity profile $\phi$.

Concerning the smoothness of solutions of \eqref{E:BoussinesqWhitham_profile}, note that differentiating \eqref{E:BoussinesqWhitham_profile} implies that
\begin{equation} \label{E:BoussinesqWhitham_smoothness}
\CalK \phi' = (c^2 - 2\phi)\phi' \implies \phi' = \frac{\CalK \phi'}{c^2-2\phi}.
\end{equation}
As before, since $\CalK \phi'$ has the same regularity as $\phi$, we have by \eqref{E:BoussinesqWhitham_smoothness} that $\phi'$ is as smooth as $\phi$ so long as 
$c^2-2\phi~\neq~0$.  Precisely, we see that smoothness should be expected to break down if the bifurcation branch of periodic solutions of \eqref{E:BoussinesqWhitham_profile}
intersects the curve $\phi=c^2/2$ non-trivially.  In our numerical bifurcation calculations, we thus stop the continuation algorithm when the maximum value of the approximated wave exceeds $c^2/2$.

In order to use \eqref{E:BoussinesqWhitham_profile} to perform a numerical continuation, we first obtain local bifurcation curves from zero amplitude in order to accurately seed the continuation algorithm.  Recasting \eqref{E:BoussinesqWhitham_profile} as a solution of
\[
Q(\phi;c) := \CalK \phi + \phi^2 - c^2\phi = 0,
\]
the variational derivative of $Q$ with respect to $\phi$ evaluated at $\phi=0$ is
\begin{equation} \label{E:BW_variational_deriv}
\delta_\phi Q(0;c) = \CalK - c^2 \mathrm{Id}.
\end{equation}
By the implicit function theorem, no bifurcation can occur unless $\ker(\CalK - c^2 \mathrm{Id})$ is nontrivial.  Restricting to even, real-valued, $2\pi/\kappa$-periodic functions, we find that for wavespeed $c_\kappa := \sqrt{\widehat{\CalK}(\kappa)} = \sqrt{\tanh(\kappa)/\kappa}$, the kernel of \eqref{E:BW_variational_deriv} is indeed non-trivial, with
\[
\ker(\CalK - c_\kappa^2 \mathrm{Id}) = \Span\{\cos(\kappa x)\}.
\]
Using a standard Lyapunov-Schmidt reduction argument, one can now prove the existence of a local bifurcation curve $\left\{(c(\varepsilon),\phi(x;\varepsilon))\right\}_{|\varepsilon|\ll 1}$ 
of non-trivial $2\pi/\kappa$-periodic traveling wave solutions $\phi(x;\varepsilon)$ of \eqref{E:BoussinesqWhitham_profile} with wave speed $c(\varepsilon)$.
In fact, we obtain the following asymptotic formulas describing the local bifurcation curve for $|\varepsilon|\ll 1$:
\begin{align*}
c(\varepsilon) &= c_\kappa + \frac{\varepsilon^2}{4c_\kappa} \left[ \frac{1}{c_\kappa^2-1} + \frac{1}{c_\kappa^2 - c_{2\kappa}^2} \right] + \CalO(\varepsilon^3) \\
\phi(x;\varepsilon) &= \varepsilon\cos(\kappa x) + \frac{\varepsilon^2}{2} \left[ \frac{1}{c_\kappa^2-1} + \frac{1}{c_\kappa^2 - c_{2\kappa}^2} \right] \cos(2\kappa x) + \CalO(\varepsilon^3),
\end{align*}
where $c_\ell = \sqrt{\widehat{\CalK}(\ell )} = \sqrt{\tanh(\ell)/\ell}$.

\begin{figure}[t]
    \centering
    
    \begin{subfigure}[t]{0.5\linewidth}
    \centering
	\includegraphics[scale=0.30]{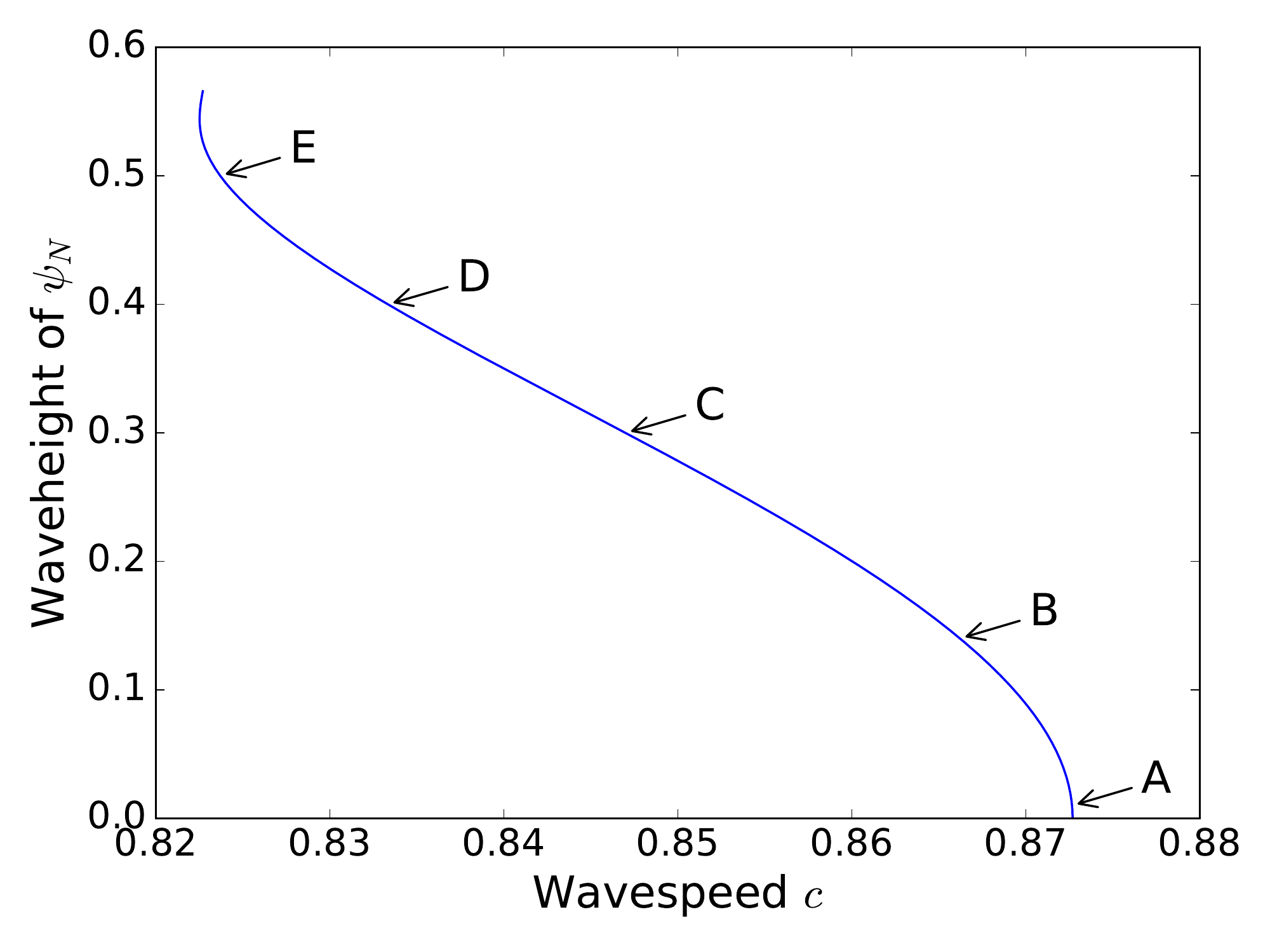}
	\caption{}
	\end{subfigure}%
	\begin{subfigure}[t]{0.5\linewidth}
	\centering
	\includegraphics[scale=0.30]{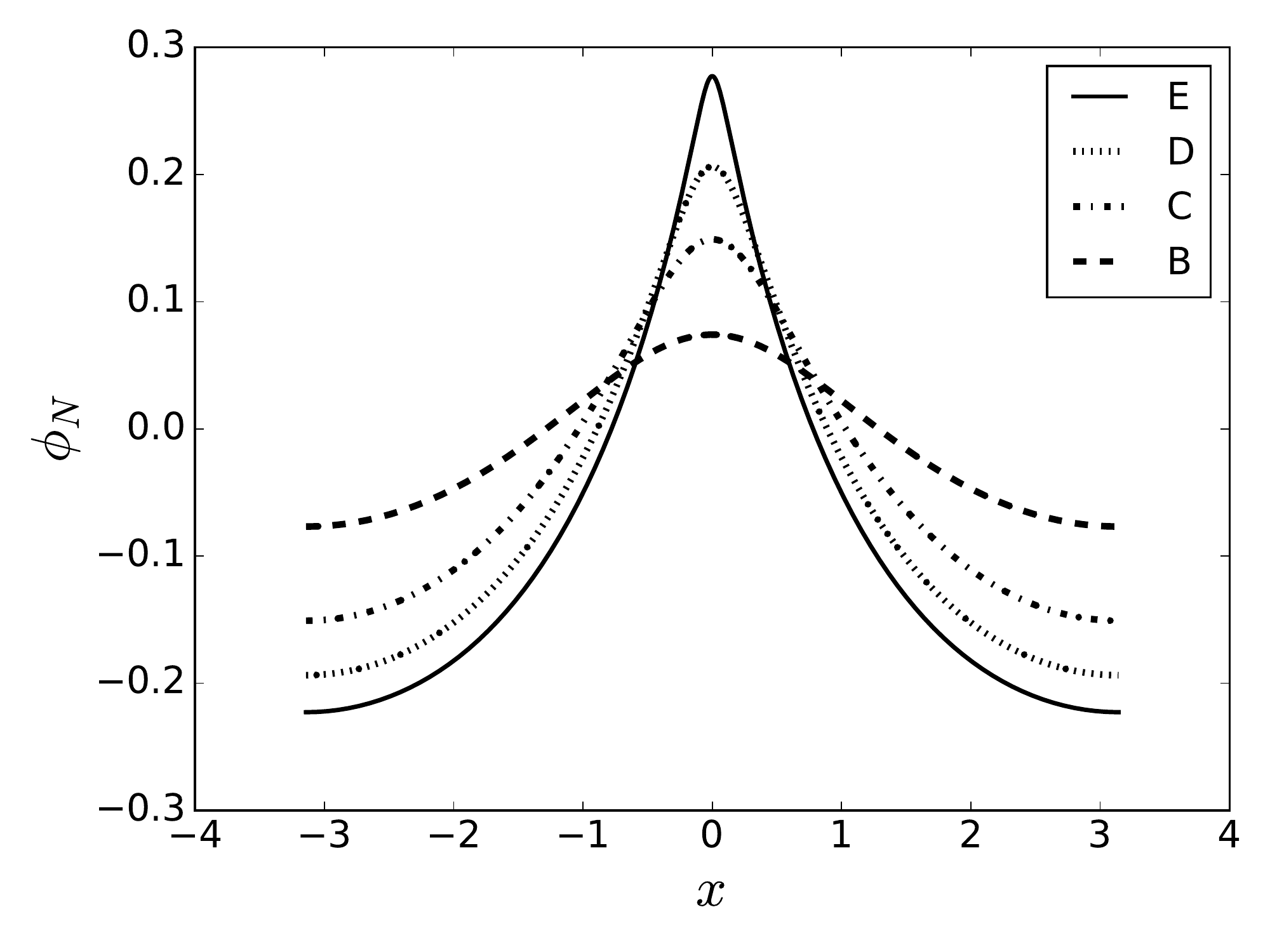}
	\caption{}
	\end{subfigure}
	
    \begin{subfigure}[t]{0.5\linewidth}
    \centering
	\includegraphics[scale=0.30]{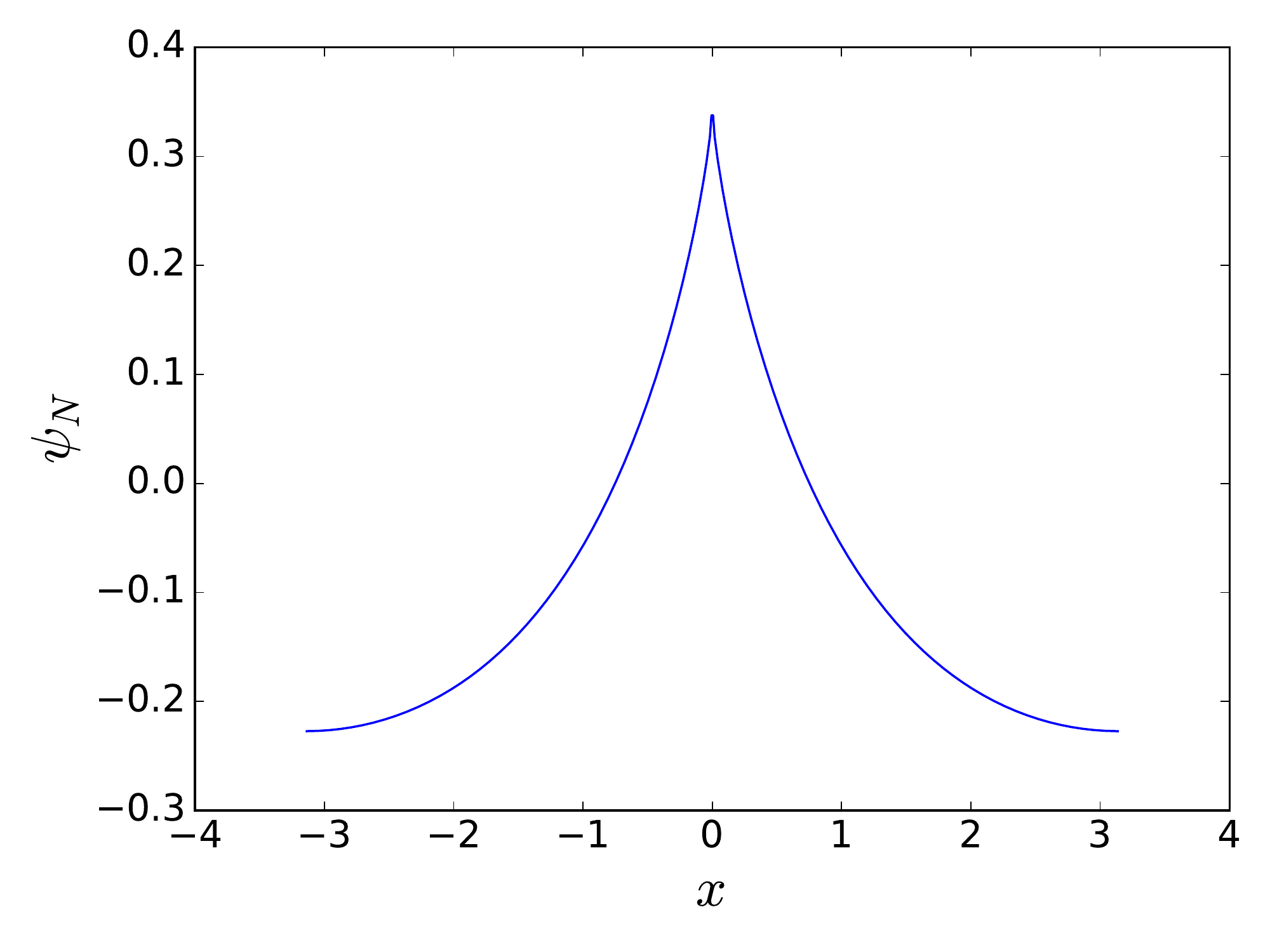}
	\caption{}
	\end{subfigure}%
	\begin{subfigure}[t]{0.5\linewidth}
	\centering
	\includegraphics[scale=0.30]{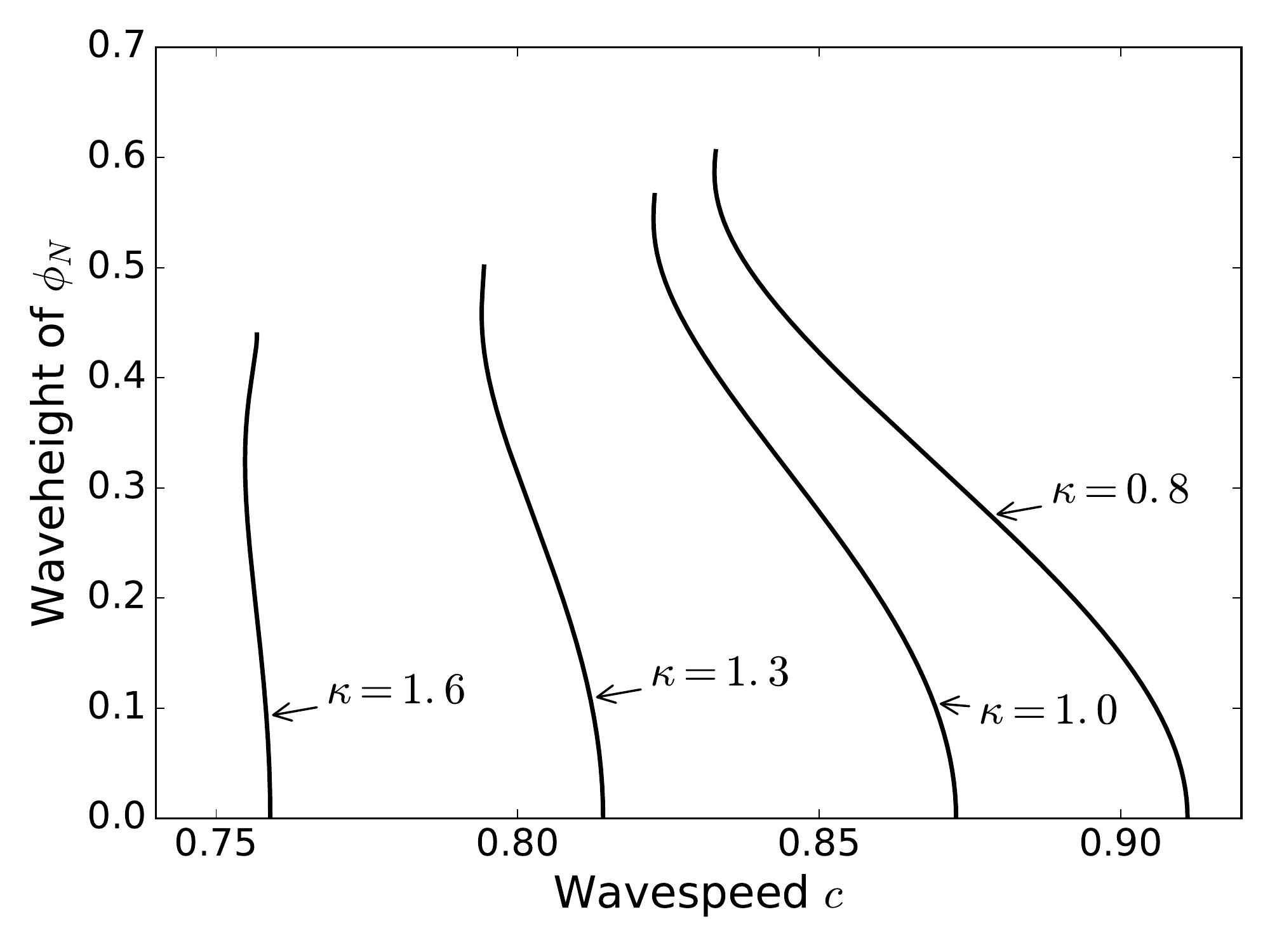}
	\caption{}
	\end{subfigure}	
	
    \caption{(a) Numerical approximation of the bifurcation branch of even, $2\pi$-periodic solutions of \eqref{E:BoussinesqWhitham_profile}, with  with specific points A-E labeled for forthcoming computations: Point A at $c \approx 0.8727$, height $\approx 0.01$; Point B at $c \approx 0.8662$, height $\approx 0.14$; Point C at $c \approx 0.8470$, height $\approx 0.30$; Point D at $c \approx 0.8333$, height $\approx 0.4$; Point E at $c \approx 0.8237$, height $\approx 0.50$.   (b) the profiles associated to the sampled points A-E in (a).  (c) A (nearly peaked) $2\pi$-periodic profile near the top of the bifurcation branch of $2\pi$-periodic solutions, corresponding to $c \approx 0.8227$, $\text{waveheight} \approx 0.5650$.  (d) The bifurcation branches of $2\pi/\kappa$-periodic solutions of \eqref{E:BoussinesqWhitham_profile} for varying $\kappa$.}
    \label{F:BoussinesqWhitham_bif}
\end{figure}

Using the above bifurcation formulas, along with the methods of Sections \ref{SS:cosine_collocation_method} and \ref{SS:pseudoarclength_method}, 
we can generate numerical bifurcation diagrams of $2\pi/\kappa$-periodic even solutions for this model: see Figures \ref{F:BoussinesqWhitham_bif}(a), (d).
Numerical approximations of the profiles at the sampled locations on the bifurcation branch of $2\pi$-periodic solutions are displayed in Figure \ref{F:BoussinesqWhitham_bif}(a)
are presented in Figure \ref{F:BoussinesqWhitham_bif}(b).  
Near the top of the bifurcation branch of $2\pi$-periodic solutions, the profiles begin to display a cusp singularity, similar to what was observed in the model \eqref{E:EJ} in Section \ref{SS:analysis_EJ} above.  See Figure \ref{F:BoussinesqWhitham_bif}(c).

\begin{figure}[t]
    \centering
    
    \begin{subfigure}[t]{0.33\linewidth}
	\includegraphics[scale=0.26]{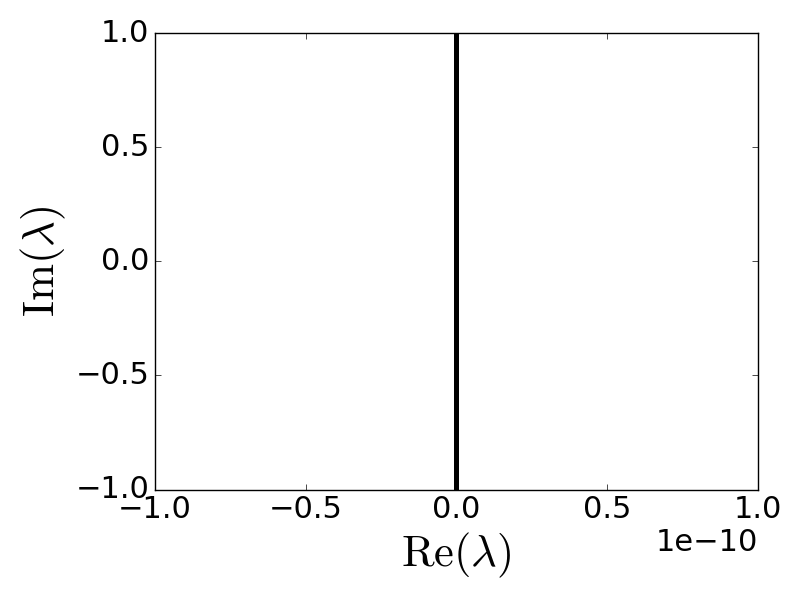} 
	\caption{}
	\end{subfigure}%
    \begin{subfigure}[t]{0.33\linewidth}
	\includegraphics[scale=0.26]{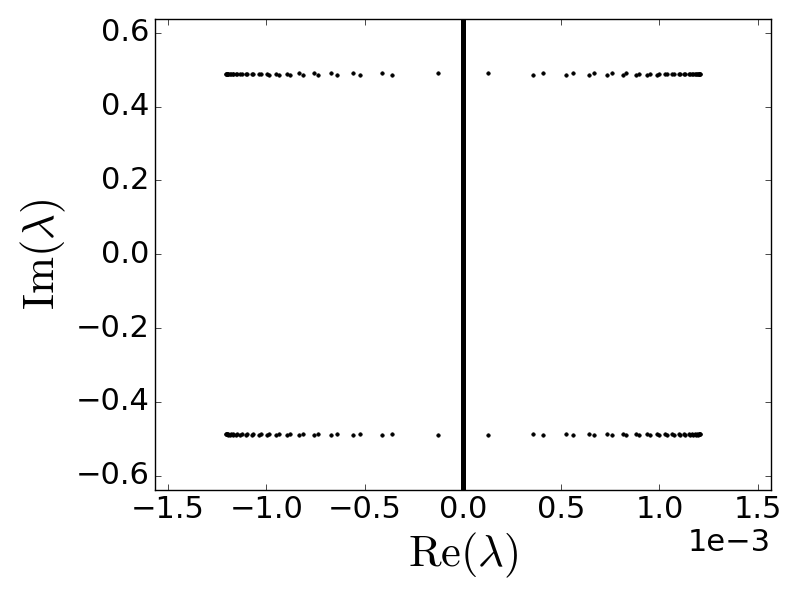}
	\caption{}
	\end{subfigure}%
	\begin{subfigure}[t]{0.33\linewidth}
	\includegraphics[scale=0.26]{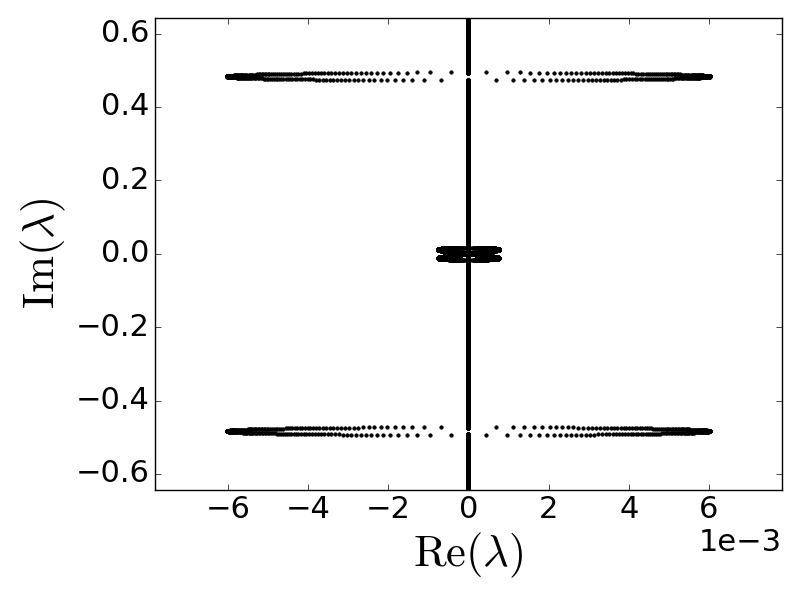}
	\caption{}
	\end{subfigure}%
	
    \begin{subfigure}[t]{0.33\linewidth}
	\includegraphics[scale=0.26]{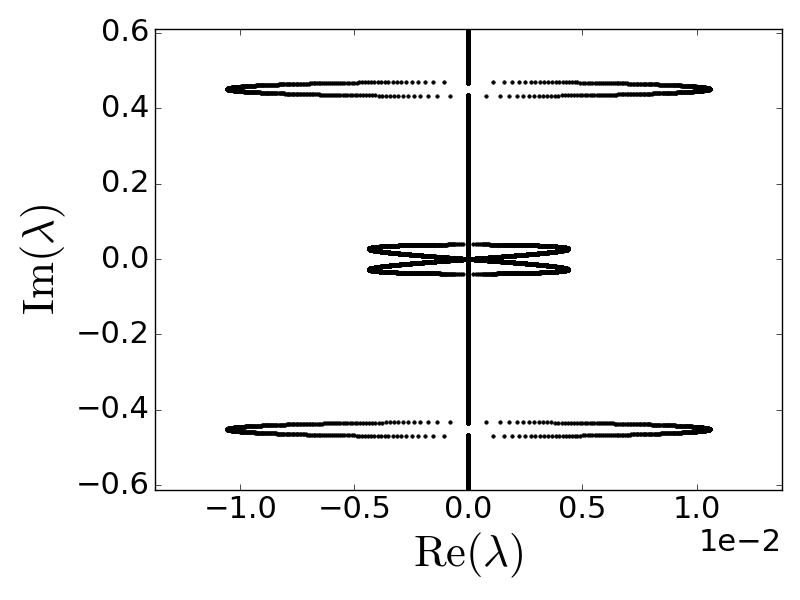}
	\caption{}
	\end{subfigure}%
    \begin{subfigure}[t]{0.33\linewidth}
	\includegraphics[scale=0.26]{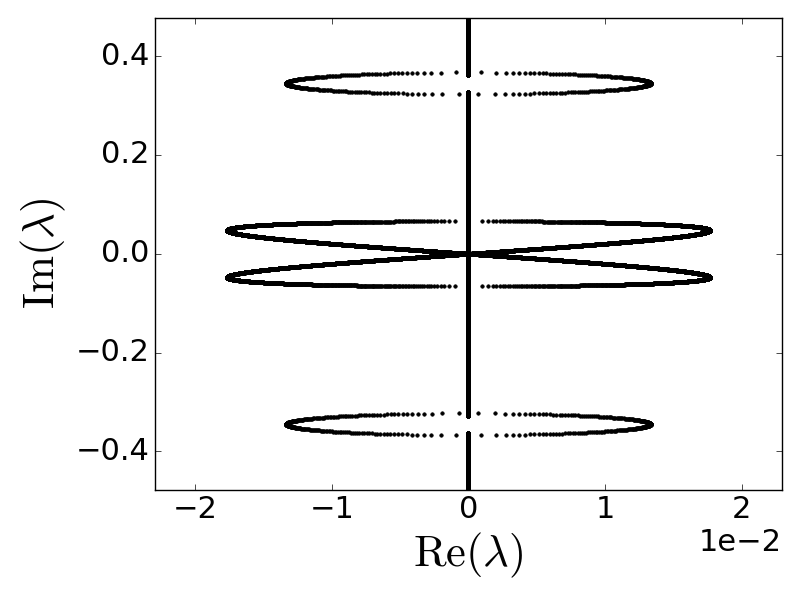}
	\caption{}
	\end{subfigure}%
    \begin{subfigure}[t]{0.33\linewidth}
	\includegraphics[scale=0.26]{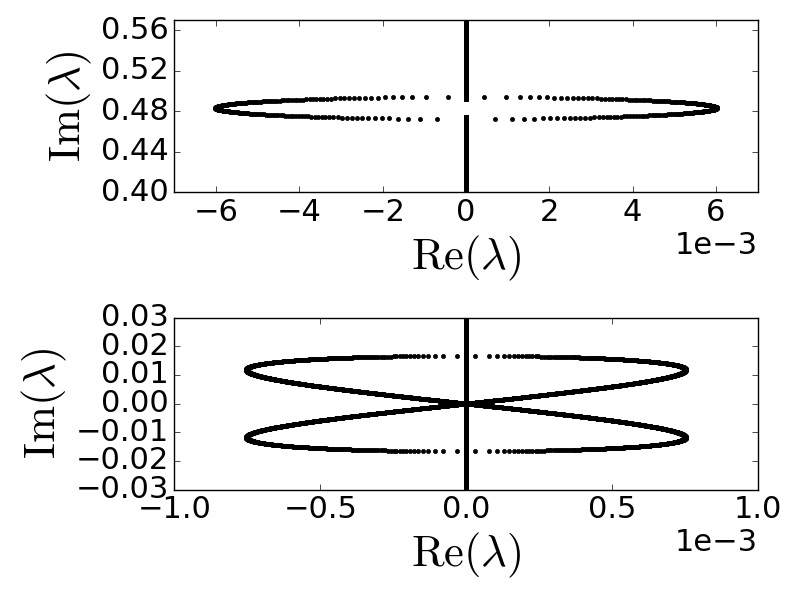}
	\caption{}
	\end{subfigure}

 	\caption{(a)-(e) Spectral plots for Boussinesq-Whitham at the points A-E, respectively, sampled along the bifurcation branch for $2\pi$-periodic solutions in Figure \ref{F:BoussinesqWhitham_bif}(a).
 (f) A zoom-in on the high-frequency instability (top) and modulational instability (bottom) for the spectral plot (c).}
    \label{F:BoussinesqWhitham_spectrum_plots}
\end{figure}

\begin{figure}[t]
    \centering
    
    \begin{subfigure}[t]{0.33\linewidth}
    \centering
	\includegraphics[scale=0.28]{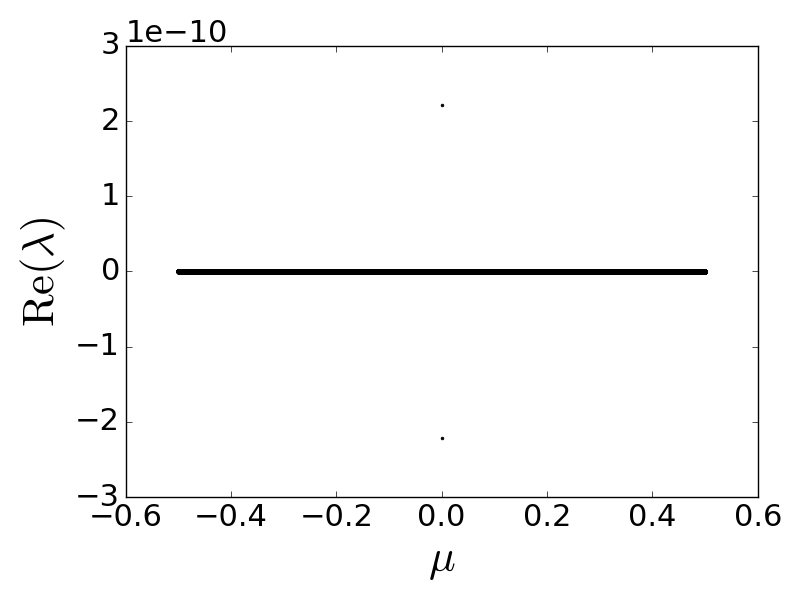}
	\caption{}
	\end{subfigure}%
    \begin{subfigure}[t]{0.33\linewidth}
    \centering
	\includegraphics[scale=0.28]{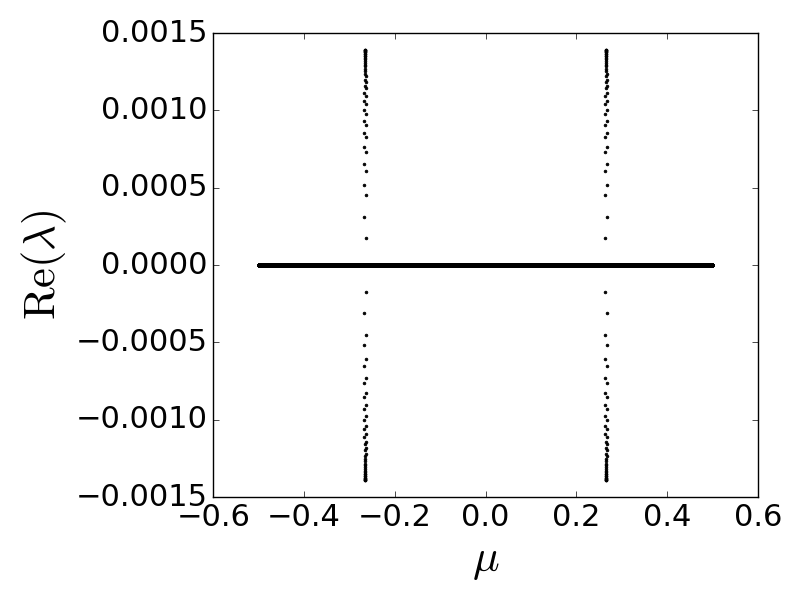}
	\caption{}
	\end{subfigure}%
	\begin{subfigure}[t]{0.33\linewidth}
	\centering
	\includegraphics[scale=0.28]{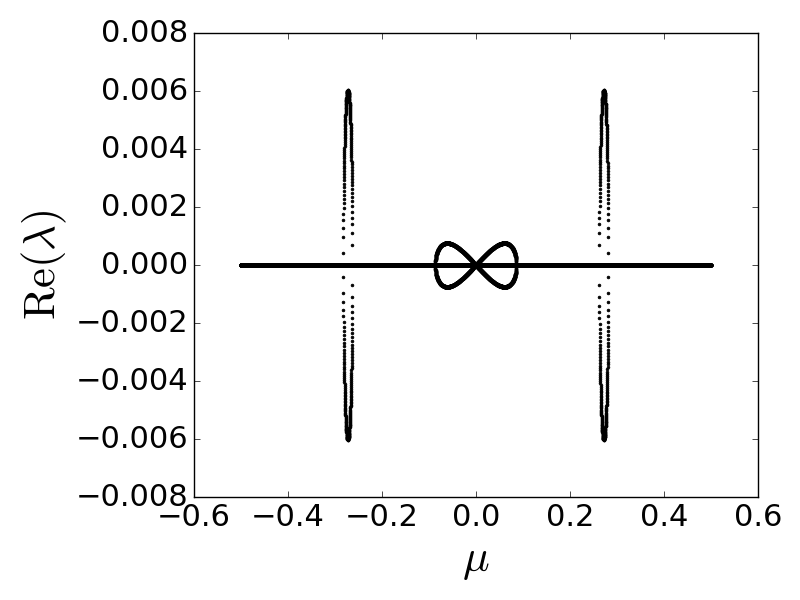}
	\caption{}
	\end{subfigure}%
	
    \begin{subfigure}[t]{0.33\linewidth}
    \centering
	\includegraphics[scale=0.28]{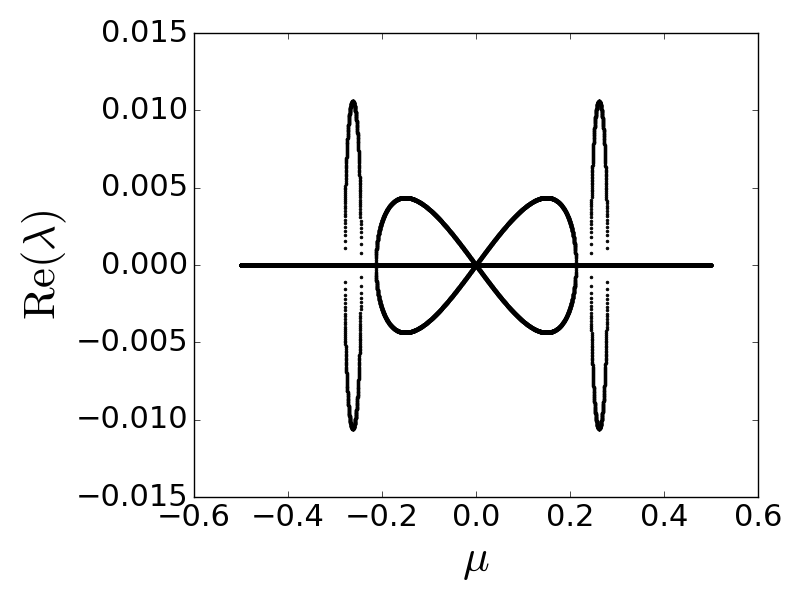}
	\caption{}
	\end{subfigure}%
	\hspace{1em}
    \begin{subfigure}[t]{0.33\linewidth}
    \centering
	\includegraphics[scale=0.28]{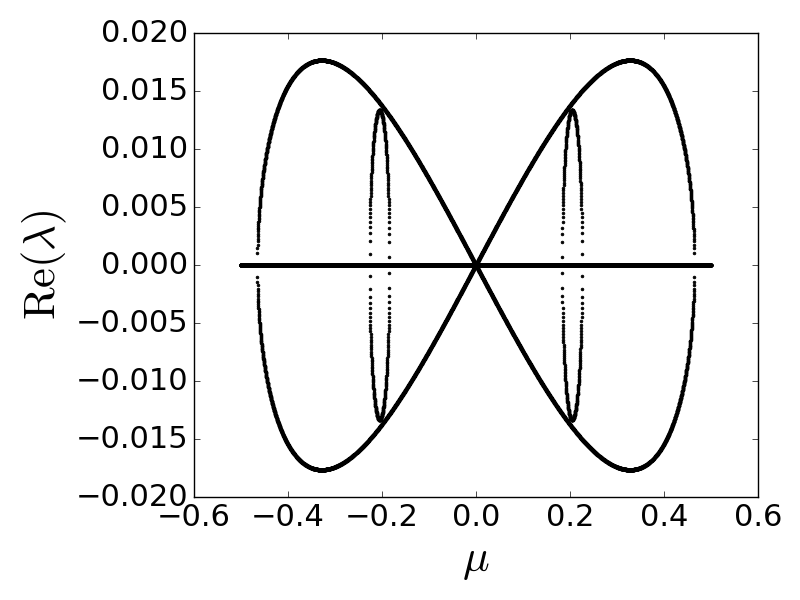}
	\caption{}
	\end{subfigure}

	\caption{(a)-(e) $\Re(\lambda)$ vs. $\mu$ for each of the spectral plots (a)-(e), respectively, of Boussinesq-Whitham shown in Figure \ref{F:BoussinesqWhitham_spectrum_plots}.}    
    \label{F:BoussinesqWhitham_Re_vs_mu}
\end{figure}

To study the stability of these periodic traveling wave solutions, we again follow the method of Appendix \ref{A:general_FFHM}
and define the following bi-infinite matrices whose eigenvalues, after finite-dimensional truncation, approximate the spectrum of the corresponding linearization:
\begin{align*}
\widehat{A}^\mu(m,l) &= ic\kappa(\mu+l)\delta_{m,l} \\
\widehat{B}^\mu(m,l) &= \delta_{m,l} \\
\widehat{C}^\mu(m,l) &= -2\kappa^2 (\mu+m)^2 \widehat{\phi}(m-l) - \kappa^2 (\mu+l)^2 \widehat{\CalK}(\kappa(\mu+l))\delta_{m,l} \\
\widehat{D}^\mu(m,l) &= ic\kappa(\mu+l) \delta_{m,l}.
\end{align*}
Plots of the spectrum, along with the growth rates $\Re(\lambda)$ with respect to the Bloch parameter $\mu$, are shown in Figures \ref{F:BoussinesqWhitham_spectrum_plots} and \ref{F:BoussinesqWhitham_Re_vs_mu}, respectively.  
Based on these numerics, we see that, similar to the other models considered here, all waves of sufficiently large waveheight appear to exhibit
both modulational and high-frequency instabilities, while high-frequency instabilities for asymptotically small waves are unobserved in our experiments.  
Furthermore, all computed waves appear to be spectrally stable with respect to co-periodic perturbations.

\section{Conclusion}\label{S:conclusion}

Using well-conditioned numerical methods, we have presented numerically-computed 
waveheight vs. wavespeed global bifurcation diagrams for three fully-dispersive bidirectional Whitham models and studied the spectral stability of such waves
by numerically approximating the spectrum of their associated linearized operators in both small and large amplitude regimes.  Our results confirm a number
of analytical results concerning the stability of asymptotically small waves in these models and provide new insight into the existence and stability of
large amplitude waves, including highest singular waves.  We note that while these models compare similarly with regard to existence and stability, we find evidence
that the model \eqref{E:HPShallow_model} does not have a highest singular wave.  Furthermore, we provide numerical evidence that the conditional
well-posedness result in \cite{EPW17} is in fact \emph{sharp}, in the sense that local evolution of initial data with height profiles with negative minima appears to be ill-posed.
This has led us to construct and analyze new wavetrain solutions of \eqref{E:EJ} with strictly positive height profile, ensuring the local dynamics about such waves
is indeed well-posed.

\appendix

\section{General Framework for Computing the Spectrum of the Linearization}  \label{A:general_FFHM}

In this appendix we seek to generalize the method for shown in Section \ref{SS:analysis_EJ} for setting up a bi-infinite matrix representation of a linearized operator, whose spectrum is suitably approximated by truncating to finite dimensions.

Recall that we are dealing with nonlinear traveling wave models of the form \eqref{E:general_model_travwave}.  To linearize the operator associated with such a model about an equilibrium solution $(u(x,t),\eta(x,t)) = (\phi(x),\psi(x))$, we first express solutions $u$ and $\eta$ as localized perturbations from equilibrium:
\begin{equation}  \label{E:equilibrium_plus_perturbation}
\begin{array}{rcl}
u(x,t) &=& \phi(x) + \varepsilon v(x,t) + \CalO(\varepsilon^2) \\
\eta(x,t) &=& \psi(x) + \varepsilon w(x,t) + \CalO(\varepsilon^2)
\end{array}
\end{equation}
with $v(\cdot,t),w(\cdot,t)\in L^2(\mathbb{R})$ for each $t>0$ for which they are defined.
Substituting \eqref{E:equilibrium_plus_perturbation} into \eqref{E:general_model}, Taylor expanding, and taking $\varepsilon \to 0$ yields
\begin{equation}  \label{E:linearization}
\left\{
\begin{array}{rcl}
v_t &=& \displaystyle \sum_{j=0}^n \left[ F_{\partial_x^j u}(\vec{\phi},\vec{\psi}) \partial_x^j v + F_{\partial_x^j u}(\vec{\phi},\vec{\psi}) \partial_x^j w \right] \\
w_t &=& \displaystyle \sum_{j=0}^n \left[ G_{\partial_x^j u}(\vec{\phi},\vec{\psi}) \partial_x^j v + G_{\partial_x^j u}(\vec{\phi},\vec{\psi}) \partial_x^j w \right]
\end{array}
\right.
\end{equation}
where
\[
\vec{\phi} := (\phi, \phi', \phi'', \ldots \phi^{(n)}), \quad \vec{\psi} := (\psi, \psi', \psi'', \ldots \psi^{(n)}).
\]
As was done in \cite{SKCK14} for the Whitham equation, we will apply Bloch theory to study the $2\pi$-periodic spectrum of the linearization by writing
\begin{equation}  \label{E:separation_of_variables}
v(x,t) = e^{\lambda t} V(x) + c.c., \quad w(x,t) = e^{\lambda t} W(x) + c.c.,
\end{equation}
where $c.c.$ denotes complex conjugate of the preceding expression (in order to ensure the reality of the functions $v$ and $w$), and
\begin{equation}  \label{E:bloch_form}
V(x) = \sum_{l \in \Z} \widetilde{V}(l) e^{i\kappa(\mu+l)x}, \quad
W(x) = \sum_{l \in \Z} \widetilde{W}(l) e^{i\kappa(\mu+l)x}
\end{equation}
with \emph{Bloch coefficients} $\widetilde{V}(l), \widetilde{W}(l)$ for $l \in \Z$ and \emph{Bloch parameter} $\mu \in [0,1)$.  Per equations \eqref{E:linearization}, we will need to understand how the operators $F_{\partial_x^j u}$, $F_{\partial_x^j \eta}$, $G_{\partial_x^j u}$, and $G_{\partial_x^j \eta}$ act on the Bloch decompositions \eqref{E:bloch_form}.  In the models of interest, these operators will be linear combinations of $2\pi$-periodic multiplication operators and the (pseudo)differential operators $\CalK, \partial_x$.  We derive the Bloch transforms of such operators here.

If $f$ is a $2\pi/\kappa$-periodic function (multiplication operator), then we may represent $f$ in Fourier series as
\[
f(x) = \sum_{m \in \Z} \widehat{f}(m) e^{im\kappa x}.
\]
Then
\begin{align*}
f(x)V(x) &= \left( \sum_{m \in \Z} \widehat{f}(m) e^{im\kappa x} \right) \left( \sum_{l \in \Z} \widetilde{V}(l) e^{i\kappa(\mu+l)x} \right) \\
&= \sum_{m \in \Z} \sum_{l \in \Z} \widehat{f}(m) \widetilde{V}(l) e^{i\kappa(\mu+m+l)x} \\
&= \sum_{m \in \Z} \left( \sum_{l \in \Z} \widehat{f}(m-l)\widetilde{V}(l) \right) e^{i\kappa(\mu+m)x} \quad \text{by taking $m \mapsto m-l$.}
\end{align*}
Hence the Bloch transform of the product $fV$ is given by
\begin{equation}  \label{E:bloch_product}
\widetilde{fV}(m) = \sum_{l \in \Z} \widehat{f}(m-l)\widetilde{V}(l), \quad m \in \Z.
\end{equation}
Moreover, the Bloch transform with parameter $\mu$ acts on the derivative operator $\partial_x$ as
\[
\widetilde{\partial_x f}(m) = \widehat{\partial_x}(\kappa(\mu+m))\widetilde{f}(m) = i\kappa(\mu+m) \widetilde{f}(m), \quad m \in \Z ,
\]
and, by analogy, the Bloch transform acts on the pseudodifferential operator $\CalK$ as
\begin{equation}
\widetilde{\CalK f}(m) := \widehat{\CalK}(\kappa(\mu+m)) \widetilde{f}(m) = \frac{\tanh(\kappa(\mu+m))}{\kappa(\mu+m)} \widetilde{f}(m), \quad m \in \Z.
\end{equation}
Substituting \eqref{E:separation_of_variables} into \eqref{E:linearization} and taking the Bloch transform of both sides, we have the simultaneous eigenvalue problems
\begin{align}
\lambda \widetilde{V}(m) &= \sum_{j=0}^n \left( F_{\partial_x^j u}(\vec{\phi},\vec{\psi}; c) V^{(j)} + F_{\partial_x^j \eta}(\vec{\phi},\vec{\psi}; c) W^{(j)} \right)^{\sim}(m)  \nonumber \\
&=: \sum_{l \in \Z} \left[ \widehat{A}^\mu(m,l)\widetilde{V}(l) + \widehat{B}^\mu(m,l)\widetilde{W}(l) \right] \label{E:VHat_eigenvalue_problem2} \\
\text{and} \quad \lambda \widetilde{W}(m) &= \sum_{j=0}^n \left( G_{\partial_x^j u}(\vec{\phi},\vec{\psi}; c) V^{(j)} + G_{\partial_x^j \eta}(\vec{\phi},\vec{\psi}; c) W^{(j)} \right)^{\sim}(m)  \nonumber \\
&=: \sum_{l \in \Z} \left[ \widehat{C}^\mu(m,l)\widetilde{V}(l) + \widehat{D}^\mu(m,l)\widetilde{W}(l) \right], \label{E:WHat_eigenvalue_problem2}
\end{align}
where $\widehat{A}^\mu(m,l)$, $\widehat{B}^\mu(m,l)$, $\widehat{C}^\mu(m,l)$, and $\widehat{D}^\mu(m,l)$ are the coefficients of the Bloch coefficients upon expanding the Bloch transform.  Define the following bi-infinite matrices entry-wise for row-$m$ and column-$l$ with $m,l \in \Z$ as
\begin{equation}
\widehat{A}^\mu := [ \widehat{A}^\mu(m,l) ]_{m,l \in \Z}, \quad \widehat{B}^\mu := [ \widehat{B}^\mu(m,l) ]_{m,l \in \Z}, \quad \widehat{C}^\mu := [ \widehat{C}^\mu(m,l) ]_{m,l \in \Z}, \quad \widehat{D}^\mu := [ \widehat{D}^\mu(m,l) ]_{m,l \in \Z}.
\end{equation}
Writing \eqref{E:VHat_eigenvalue_problem2} and \eqref{E:WHat_eigenvalue_problem2} jointly in block bi-infinte matrix form, we have
\begin{equation}  \label{E:FFHM_eigenvalue_problem2}
\def\arraystretch{1.5}  
\lambda \left[ \begin{array}{c}
\widetilde{V} \\
\hline
\widetilde{W}
\end{array} \right]
= \left[
\begin{array}{c|c}
\widehat{A}^\mu & \widehat{B}^\mu \\
\hline
\widehat{C}^\mu & \widehat{D}^\mu 
\end{array}
\right]
\left[ \begin{array}{c}
\widetilde{V} \\
\hline
\widetilde{W}
\end{array} \right]
=: \widehat{\CalL}^\mu 
\left[ \begin{array}{c}
\widetilde{V} \\
\hline
\widetilde{W}
\end{array}
\right]
\end{equation}
where
\begin{align*}
\widetilde{V} &:= \begin{bmatrix} \ldots & \widetilde{V}(-2) & \widetilde{V}(-1) & \widetilde{V}(0) & \widetilde{V}(1) & \widetilde{V}(2) & \ldots \end{bmatrix}^T \\
\widetilde{W} &:= \begin{bmatrix} \ldots & \widetilde{W}(-2) & \widetilde{W}(-1) & \widetilde{W}(0) & \widetilde{W}(1) & \widetilde{W}(2) & \ldots \end{bmatrix}^T.
\end{align*}
To numerically approximate the bi-infinite eigenvalue problem \eqref{E:FFHM_eigenvalue_problem2}, we truncate each of the bi-infinite matrices $\widehat{A}^\mu$, $\widehat{B}^\mu$, $\widehat{C}^\mu$, and $\widehat{D}^\mu$ to have dimension $(2N+1) \times (2N+1)$, i.e.
\[
\widehat{A}_N^\mu(m,l) := [\widehat{A}^\mu(m,l)]_{-N \leq m,l \leq N} \in \C^{(2N+1) \times (2N+1)}
\]
and similarly define $\widehat{B}_N^\mu$, $\widehat{C}_N^\mu$, and $\widehat{D}_N^\mu$, which together form $\widehat{\CalL}_N^\mu \in \C^{(4N+2) \times (4N+2)}$ as shown in \eqref{E:FFHM_discrete_eigenvalue_problem} below.  Moreover, we also truncate $\widetilde{V}$ and $\widetilde{W}$ as
\[
\widetilde{V}_N = \begin{bmatrix} \widetilde{V}(-N) & \ldots & \widetilde{V}(0) & \ldots & \widetilde{V}(N) \end{bmatrix}^T, \quad \widetilde{W}_N = \begin{bmatrix} \widetilde{W}(-N) & \ldots & \widetilde{W}(0) & \ldots & \widetilde{W}(N) \end{bmatrix}^T.
\]
Then, using a standard matrix eigenvalue solver (e.g. Matlab's \texttt{eig} or Python's \texttt{numpy.linalg.eig} \cite{SciPy}), we solve for $\lambda_N^\mu \in \C$ such that
\begin{equation}  \label{E:FFHM_discrete_eigenvalue_problem}
\def\arraystretch{1.5}  
\lambda_N^\mu \left[ \begin{array}{c}
\widetilde{V}_N \\
\hline
\widetilde{W}_N
\end{array} \right]
= \left[
\begin{array}{c|c}
\widehat{A}_N^\mu & \widehat{B}_N^\mu \\
\hline
\widehat{C}_N^\mu & \widehat{D}_N^\mu 
\end{array}
\right]
\left[ \begin{array}{c}
\widetilde{V}_N \\
\hline
\widetilde{W}_N
\end{array} \right]
=: \widehat{\CalL}^\mu_N 
\left[ \begin{array}{c}
\widetilde{V}_N \\
\hline
\widetilde{W}_N
\end{array}
\right]
\end{equation}
We solve this eigenvalue problem for $\mu$ in a discrete subset of $[0,1)$.  For our computations, we used a uniform mesh $\{\mu_j = j \Delta \mu\}_j$ with subintervals of constant width $\Delta \mu$.

\begin{remark} \label{Rmk:FFHM_convergence}
The above method is inspired by the Fourier-Floquet-Hill Method (FFHM) given in \cite{DK05}, which is used to compute the spectrum of a linear, locally-acting operator having periodic coefficients.  As was done in \cite{DK05}, we define convergence of the spectral approximation in the sense that any eigenvalue $\lambda_N^\mu$ of $\widehat{\CalL}_N^\mu$ above converges to some eigenvalue $\lambda$ of non-truncated linearization $\CalL$.  Moreover, the incorporation of all Bloch parameters $\mu \in [0,1)$ yields the entire spectrum, i.e.
\[
\lim_{N \to \infty} \bigcup_{\mu \in [0,1)} \sigma(\widehat{\CalL}_N^\mu) = \sigma(\CalL).
\]
\end{remark}
To obtain the approximate Fourier coefficients $\widehat{\phi}$ and $\widehat{\psi}$ (in the exponential basis) needed to construct the bi-infinite matrix $\widehat{\CalL}_N^\mu$, for a given wavespeed $c$, we use the approximate values of $\phi(x_m)$ at each collocation point $x_m = \frac{(2m-1)\pi}{2\kappa N}$ obtained in the numerical bifurcation (see Section \ref{S:bifurcation_methods}) and compute the approximate coefficients via their integral definitions using the midpoint quadrature:
\begin{align*}
\widehat{\phi}(n) &= \frac{\kappa}{2\pi} \int_{-\pi/\kappa}^{\pi/\kappa} \phi(x) e^{-i n \kappa x} \, dx \\
&= \frac{\kappa}{\pi} \int_0^{\pi/\kappa} \phi(x)\cos(n \kappa x) \, dx \quad \text{since $\phi$ is even} \\
&\approx \frac{1}{N} \sum_{m=1}^N \phi(x_m)\cos(n\kappa x_m) \quad \text{by \eqref{E:fourier_integral_discretization}.}
\end{align*}
Similarly,
\[
\widehat{\psi}(n) \approx \frac{1}{N} \sum_{m=1}^N \psi(x_m)\cos(n\kappa x_m).
\]
\begin{remark}
The approximation of $\widehat{\psi}(n)$ above requires that $\psi$ be an even function.  This is always the case since we are considering models for which $\psi$ can be resolved in terms of $\phi$, hence $\psi = \psi(\phi(x))$ is also an even function.

As an example, here we show the relevant calculations in detail for the model \eqref{E:EJ_traveling_wave_model} as a concrete demonstration of the more abstract framework established in this appendix.  Recall that this model is given by
\[
\left\{
\begin{array}{rcl}
u_t &=& cu_x - \eta_x - uu_x =: F(u,u_x,\eta,\eta_x) \\
\eta_t &=& c\eta_x - \CalK u_x - (\eta u)_x =: G(u,u_x,\eta,\eta_x),
\end{array}
\right.
\]
with equilibrium solutions $u(x,t) = \phi(x)$, $\eta(x,t) = \psi(x)$.
Beginning with \eqref{E:VHat_eigenvalue_problem}, we first compute
\[
F_u(\phi,\phi',\psi,\psi') = -\phi', \quad
F_{u_x}(\phi,\phi',\psi,\psi') = c-\phi, \quad
F_\eta(\phi,\phi',\psi,\psi')  = 0, \quad
F_{\eta_x}(\phi,\phi',\psi,\psi') = -1.
\]
All of the above are $2\pi/\kappa$-periodic multiplication operators, hence we apply \eqref{E:bloch_product} in \eqref{E:VHat_eigenvalue_problem2}, here with $n=1$, to achieve the following formulas for $m \in \Z$:
\begin{align*}
\widetilde{F_u V}(m) &= \sum_{l \in \Z} \widehat{-\phi'}(m-l)\widetilde{V}(l) = \sum_{l \in \Z} -i\kappa(m-l)\widehat{\phi}(m-l)\widetilde{V}(l) \\
\widetilde{F_{u_x} V'}(m) &= \sum_{l \in \Z} \widehat{(c-\phi)}(m-l) \ \widetilde{\partial_x V}(l) = \sum_{l \in \Z} \left( c \delta_{m,l} - \widehat{\phi}(m-l) \right) \cdot  i\kappa(\mu+l)\widetilde{V}(l) \\
\widetilde{F_\eta W}(m) &= \sum_{l \in \Z} 0 \, \widetilde{W}(l) \\
\widetilde{F_{\eta_x} W'}(m) &= \sum_{l \in \Z} \widehat{(-1)}(m-l) \ \widetilde{\partial_x W}(l) = \sum_{l \in \Z} -\delta_{m,l} \cdot i\kappa(\mu+l) \widetilde{W}(l),
\end{align*}
where $\delta_{m,l}$ is the Kronecker delta.  Summing the above expressions and grouping terms involving $\widetilde{V}(l)$ and $\widetilde{W}(l)$ per \eqref{E:VHat_eigenvalue_problem2}, we obtain
\begin{align}
\widehat{A}^\mu(m,l) &= -i\kappa(m-l)\widehat{\phi}(m-l) + \left( c \delta_{m,l} - \widehat{\phi}(m-l) \right) \cdot  i\kappa(\mu+l) \nonumber \\
&= ic\kappa(\mu+l)\delta_{m,l} - i\kappa(\mu+m)\widehat{\phi}(m-l) \label{E:appendix_block_matrix_A} \\
\widehat{B}^\mu(m,l) &= -i\kappa(\mu+l)\delta_{m,l}, \label{E:appendix_block_matrix_B}
\end{align}
Similarly, per \eqref{E:WHat_eigenvalue_problem2}, we compute
\[
G_u(\phi,\phi',\psi,\psi') = -\psi', \quad
G_{u_x}(\phi,\phi',\psi,\psi') = -\CalK - \psi, \quad
G_\eta(\phi,\phi',\psi,\psi') = -\phi', \quad
G_{\eta_x}(\phi,\phi',\psi,\psi') = c-\phi.
\]
Then
\begin{align*}
\widetilde{G_u V}(m) &= \sum_{l \in \Z} \widehat{-\psi'}(m-l) \widetilde{V}(l) = -i\kappa(m-l)\widehat{\psi}(m-l) \widetilde{V}(l) \\
\widetilde{G_{u_x} V'}(m) &= \left[ (-\CalK-\psi)\partial_x V \right]^{\sim}(m) \\
&= -\widetilde{\CalK \partial_x V}(m) - \widetilde{\psi \, \partial_x V}(m) \\
&= -i\kappa(\mu+m) \widehat{\CalK}(\kappa(\mu+m)) \widetilde{V}(m) - \sum_{l \in \Z} \widehat{\psi}(m-l) \, i\kappa(\mu+l) \widetilde{V}(l) \\
&= \sum_{l \in \Z} -i\kappa(\mu+l) \left[ \widehat{\CalK}(\kappa(\mu+l))\delta_{m,l} + \widehat{\psi}(m-l) \right] \widetilde{V}(l) \\
\widetilde{G_\eta W}(m) &= \sum_{l \in \Z} \widehat{-\phi'}(m-l) \widetilde{W}(l) = \sum_{l \in \Z} -i\kappa (m-l)\widehat{\phi}(m-l) \widetilde{W}(l) \\
\widetilde{G_{\eta_x} W'}(m) &= \sum_{l \in \Z} \widehat{(c-\phi)}(m-l) \ \widetilde{\partial_x W}(l) = \sum_{l \in \Z} \left( c\delta_{m,l} - \widehat{\phi}(m-l) \right) \cdot i\kappa(\mu+l)\widetilde{W}(l).
\end{align*}
Summing the above expressions and grouping terms involving $\widetilde{V}(l)$ and $\widetilde{W}(l)$ per \eqref{E:WHat_eigenvalue_problem2}, we obtain
\begin{align}
\widehat{C}^\mu(m,l) &= -i\kappa(\mu+l) \widehat{\CalK}(\kappa(\mu+l))\delta_{m,l} - i\kappa(\mu+m)\widehat{\psi}(m-l)  \label{E:appendix_block_matrix_C} \\
\widehat{D}^\mu(m,l) &= ic\kappa(\mu+l)\delta_{m,l} - i\kappa(\mu+m)\widehat{\phi}(m-l). \label{E:appendix_block_matrix_D}
\end{align}
Note that the formulas for $\widehat{A}^\mu(m,l)$, $\widehat{B}^\mu(m,l)$, $\widehat{C}^\mu(m,l)$, $\widehat{D}^\mu(m,l)$ given above in \eqref{E:appendix_block_matrix_A}, \eqref{E:appendix_block_matrix_B}, \eqref{E:appendix_block_matrix_C}, and \eqref{E:appendix_block_matrix_D} agree with those found in Section \ref{SS:analysis_EJ} circa \eqref{E:VHat_eigenvalue_problem}, \eqref{E:WHat_eigenvalue_problem} by manipulating Fourier series.
\end{remark}

\section{Numerical Parameters} \label{A:numerical_parameters}

\subsection{Bifurcation and Spectral Figures}  \label{AA:bifurcation_spectrum_parameters}

In Table \ref{T:bif_and_FFHM_parameters}, we provide the numerical parameters used to compute bifurcations via the pseudoarclength method and spectra via the Fourier-Floquet-Hill-Method (FFHM).  A brief reminder of the symbols involved:

\

\noindent Pseudoarclength Parameters:
\begin{itemize}
\item $\kappa$ is the wavenumber, which yields waves of period $2\pi/\kappa$.
\item $N$ is the number of collocation points $x_i = \frac{(2i-1)\pi}{2\kappa N}$, $i=1,\ldots,N$ on the half-period cell $[0,\pi/\kappa]$, as well as the number of modes used in the profile's truncated series $\phi(x) = \sum_{n=0}^{N-1} \widehat{\phi}(n) \cos(n\kappa x)$.  See Section \ref{SS:cosine_collocation_method}.
\item $h$ is the length of a step along the tangent direction in each iteration of the pseudoarclength method.  See Figure \ref{F:pseudoarclength_method} and the surrounding discussion.
\item $\varepsilon_0$ determines a suitable starting point (see \eqref{E:y0_star}) near the global bifurcation curve based on the model's local bifurcation formulas (e.g. \eqref{E:EJ_local_bif_formula_phi}, \eqref{E:EJ_local_bif_formula_phi}).
\end{itemize}

\noindent FFHM Parameters:
\begin{itemize}
\item $\Delta \mu$ is the width of a subinterval in the uniform mesh $\{ \mu_j = j\Delta \mu\}_j \subset [0,1)$ of Bloch parameters used to discretize the $L^2(\R)$ spectrum.  See \eqref{E:bloch_form} and Remark \ref{Rmk:FFHM_convergence}.
\item $N$ determines the dimension of the truncated bi-infinite matrix $\widehat{\CalL}_N^\mu \in \C^{(4N+2) \times (4N+2)}$. See  \eqref{E:FFHM_discrete_eigenvalue_problem} and the surrounding discussion.
\end{itemize}

\begin{table}[h]
\caption{Parameters for figures generated by the pseudoarclength method (bifurcation) and Fourier-Floquet-Hill method (spectrum).}
\label{T:bif_and_FFHM_parameters}
\begin{tabular}{l|l|l}
Figure(s) & Psuedoarclength Parameters & FFHM Parameters (if applicable) \\
\hline
\ref{F:EJ_global_bif_diagram} (a) & $\kappa=1$, $N=256$, $h=0.001$, $\varepsilon_0=1 \times 10^{-5}$ \\
\ref{F:EJ_global_bif_diagram} (b) & $\kappa \in \{ 0.8,1.0,1.3,1.6 \}$, $N=256$, $h=0.001$, $\varepsilon_0=1 \times 10^{-5}$ \\
\ref{F:EJ_profiles} (a), (b) & $\kappa=1$, $N=256$, $h=0.001$, $\varepsilon_0=1 \times 10^{-5}$ \\
\ref{F:EJ_spectrum_plots} (a) &  $\kappa=1$, $N=256$, $h=0.001$, $\varepsilon_0=1 \times 10^{-5}$ & $N=50$, $\Delta\mu = 1/10000$ \\
\ref{F:EJ_spectrum_plots} (b), (c) & $\kappa=1$, $N=256$, $h=0.001$, $\varepsilon_0=1 \times 10^{-5}$ & $N=50$, $\Delta\mu = 1/5000$ \\
\ref{F:EJ_spectrum_plots} (d), (e) & $\kappa=1$, $N=256$, $h=0.001$, $\varepsilon_0=1 \times 10^{-5}$ & $N=128$, $\Delta\mu = 1/5000$ \\
\ref{F:EJ_Re_vs_mu} (a)-(e) & Same data as Figure \ref{F:EJ_spectrum_plots} (a)-(e) & Same data as Figure \ref{F:EJ_spectrum_plots} (a)-(e) \\
\ref{F:EJ_critical_kappa_comparison} (a) & $\kappa=1.005$, $N=64$, $h=0.001$, $\varepsilon_0=1 \times 10^{-5}$ & $N=50$, $\Delta\mu = 1/15000$ \\
\ref{F:EJ_critical_kappa_comparison} (b) & $\kappa=1.008$, $N=64$, $h=0.001$, $\varepsilon_0=1 \times 10^{-5}$ & $N=64$, $\Delta\mu = 1/15000$ \\
\ref{F:EJ_positive_global_bif_diagram} (b) & $\kappa \in \{ 0.5, 1.0, 1.3, 1.6 \}$, $N=256$, $h=0.001$, $\varepsilon_0=1 \times 10^{-5}$ & \\
\ref{F:EJ_positive_bif_profiles} (a), (b) & $\kappa=1$, $N=2048$, $h=0.001$, $\varepsilon_0=1 \times 10^{-5}$ & \\
\ref{F:EJ_positive_spectrum_plots} (a) & $\kappa=1$, $N=2048$, $h=0.001$, $\varepsilon_0=1 \times 10^{-5}$ & $N=50$, $\Delta\mu = 1/50000$ \\
\ref{F:EJ_positive_spectrum_plots} (b) & $\kappa=1$, $N=2048$, $h=0.001$, $\varepsilon_0=1 \times 10^{-5}$ & $N=50$, $\Delta\mu = 1/10000$ \\
\ref{F:EJ_positive_spectrum_plots} (c) & $\kappa=1$, $N=2048$, $h=0.001$, $\varepsilon_0=1 \times 10^{-5}$ & $N=128$, $\Delta\mu = 1/5000$ \\
\ref{F:HPShallow_supcrit_bifurcation_and_profiles} (a)-(c) & $\kappa=1.611$, $N=256$, $h=0.001$, $\varepsilon_0=1 \times 10^{-5}$ \\
\ref{F:HPShallow_supcrit_bifurcation_and_profiles} (d) & $\kappa \in \{0.8,1,1.3,1.6\}$, $N=256$, $h=0.001$, $\varepsilon_0=1 \times 10^{-5}$ \\
\ref{F:HPShallow_secondary_turning_point} & $\kappa=1$, $N \in \{64,128,256,512\}$, $h=0.001$, $\varepsilon_0=1 \times 10^{-5}$ \\
\ref{F:HPShallow_supcrit_spectrum_plots} (a), (b) & $\kappa=1.611$, $N=256$, $h=0.001$, $\varepsilon_0=1 \times 10^{-5}$ & $N=50$, $\Delta\mu = 1/10000$ \\
\ref{F:HPShallow_supcrit_spectrum_plots} (c), (d) & $\kappa=1.611$, $N=256$, $h=0.001$, $\varepsilon_0=1 \times 10^{-5}$ & $N=50$, $\Delta\mu = 1/5000$ \\
\ref{F:HPShallow_supcrit_spectrum_plots} (e) & $\kappa=1.611$, $N=256$, $h=0.001$, $\varepsilon_0=1 \times 10^{-5}$ & $N=128$, $\Delta\mu = 1/5000$ \\
\ref{F:HPShallow_subcrit_spectrum_plots} (a), (b) & $\kappa=1.609$, $N=256$, $h=0.001$, $\varepsilon_0=1 \times 10^{-5}$ & $N=50$, $\Delta\mu = 1/10000$ \\
\ref{F:HPShallow_supcrit_Res} (a)-(e) & Same data as Figure \ref{F:HPShallow_supcrit_spectrum_plots} & Same data as Figure \ref{F:HPShallow_supcrit_spectrum_plots} \\
\ref{F:BoussinesqWhitham_bif} (a)-(c) & $\kappa=1$, $N=256$, $h=0.001$, $\varepsilon_0=1 \times 10^{-5}$ & \\
\ref{F:BoussinesqWhitham_bif} (d) & $\kappa \in \{ 0.8, 1.0, 1.3, 1.6 \}$, $N=256$, $h=0.001$, $\varepsilon_0=1 \times 10^{-5}$ \\
\ref{F:BoussinesqWhitham_spectrum_plots} (a)-(d) & $\kappa=1$, $N=256$, $h=0.001$, $\varepsilon_0=1 \times 10^{-5}$ & $N=50$, $\Delta\mu = 1/10000$ \\
\ref{F:BoussinesqWhitham_spectrum_plots} (e) & $\kappa=1$, $N=256$, $h=0.001$, $\varepsilon_0=1 \times 10^{-5}$ & $N=128$, $\Delta\mu = 1/10000$ \\
\ref{F:BoussinesqWhitham_Re_vs_mu} & Same data as Figure \ref{F:BoussinesqWhitham_spectrum_plots} & Same data as Figure \ref{F:BoussinesqWhitham_spectrum_plots}
\end{tabular}
\end{table}

\vspace{1em}

\subsection{Time Evolution} \label{AA:time_evolution_parameters}

We use a sixth-order pseudospectral operator splitting method \cite{Yoshida1990} to evolve a traveling wave solution of \eqref{E:EJ} by successively evolving the linear and nonlinear parts of the system.  The linear part of \eqref{E:EJ} is
\[
\left\{
\begin{array}{rcl}
u_t &=& -\eta_x \\
\eta_t &=& -\CalK u_x,
\end{array}
\right.
\]
which can be transformed into decoupled wave-like second-order in time IVPs
\[
\left\{
\begin{array}{rcl}
\eta_{tt}(x,t) &=& \CalK\eta_{xx}(x,t) \\
\eta(x,0) &=& \psi(x) \\
\eta_t(x,0) &=& -\CalK\phi'(x)
\end{array}
\right.
\quad
\left\{
\begin{array}{rcl}
u_{tt}(x,t) &=& \CalK u_{xx}(x,t) \\
u(x,0) &=& \phi(x) \\
u_t(x,0) &=& -\psi'(x),
\end{array}
\right.
\]
where $\psi(x)$, $\phi(x)$ are equilibria as discussed in Section \ref{SS:analysis_EJ}.
Taking the Fourier transform of the above equations in $x$ yields ODE IVPs in $t$ for each fixed $n \in \Z$:
\begin{equation} \label{E:EJ_time_evo_formulas}
\left\{
\begin{array}{rcl}
\widehat{\eta}_{tt}(n,t) &=& -\kappa^2 n^2 \widehat{\CalK}(n)\widehat{\eta}(n,t) \\
\widehat{\eta}(n,0) &=& \widehat{\psi}(n) \\
\widehat{\eta}_t(n,0) &=& -in\widehat{\CalK}(n) \widehat{\phi}(n)
\end{array}
\right.
\quad
\left\{
\begin{array}{rcl}
\widehat{u}_{tt}(n,t) &=& -\kappa^2 n^2 \widehat{\CalK}(n)\widehat{u}(n,t) \\
\widehat{u}(n,0) &=& \widehat{\phi}(n) \\
\widehat{u}_t(n,0) &=& -in\widehat{\psi}(n).
\end{array}
\right.
\end{equation}
The equations in \eqref{E:EJ_time_evo_formulas} are well-known to model simple harmonic motion and can be solved explicitly.  For the nonlinear part of \eqref{E:EJ}, the typical fourth-order Runge-Kutta scheme (RK4) is used.

\begin{itemize}
\item In Figure \ref{F:EJ_time_evo}, a small $2\pi$-periodic wave of height $\max\psi - \min\psi \approx 0.00109$ generated by the pseudoarclength method with parameters $\kappa=1$, $h=0.01$, $\varepsilon_0=1 \times 10^{-5}$, and  $N=2048$ collocation points on the half-periodic cell $[0,\pi]$ ($\Delta x = \frac{\pi}{2048} \approx 0.00153$) is used as initial data in \eqref{E:EJ_time_evo_formulas}.  A time step of $\Delta t = 0.001$ is used in each iteration of the 6th-order operator splitting method to integrate to time $t=2.55$, which is about 35.4\% of a temporal period for this wave.

\item In Figure \ref{F:EJ_positive_time_evo_15per}, a large $2\pi$-periodic \emph{positive} wave of height $\max\psi - \min\psi \approx 0.387$ generated by the pseudoarclength method with parameters $\kappa=1$, $h=0.01$, $\varepsilon_0=1 \times 10^{-5}$, and  $N=2048$ collocation points on the half-periodic cell $[0,\pi]$ ($\Delta x = \frac{\pi}{2048} \approx 0.00153$) is used as initial data in \eqref{E:EJ_time_evo_formulas}.  A time step of $\Delta t = 0.001$ was used to integrate for 15 temporal periods (to time $t \approx 80.5706$).
\end{itemize}

The ease of time-evolving large positive waves (which are proven to be locally well-posed; see \cite{EPW17}), as opposed to the difficulty of evolving non-positive waves in the same model with the same $\Delta t/\Delta x \approx 0.6536$ suggests that the non-positive waves are not locally well-posed in time.



\bibliographystyle{plain}
\bibliography{bwhitham.bib}

\end{document}